\documentclass[12pt,english,reqno]{amsart}
\usepackage{amssymb,amsmath,amstext,amsthm,amsfonts}
\usepackage[dvips]{graphicx}%%%PARA FAZER GRAFICOS EPS%%%
\usepackage{enumerate}
\usepackage[ansinew]{inputenc}
\usepackage{a4wide}
\newtheorem{Theorem}{Theorem}

\newtheorem{maintheorem}{Theorem}

\newtheorem{T}{Theorem}[section]
\newtheorem{Corollary}[T]{Corollary}
\newtheorem{Proposition}[T]{Proposition}
\newtheorem{Lemma}[T]{Lemma}

\newtheorem{Notation}[T]{Notation}
\newtheorem{Remark}[T]{Remark}

\newtheorem{Definition}[T]{Definition}
\newtheorem{Example}[T]{Example}
\newtheorem*{claim}{Claim}
\newtheorem{Claim}{Claim}

\def \DD {{\mathbb D}}
\def \RR {{\mathbb R}}
\def \ZZ {{\mathbb Z}}
\def \NN {{\mathbb N}}

\def \diam {\mbox{diameter}}
\def \interior {\mbox{int}}
\def \tz {{\mathrm{Z}}}
\def \th {{\mathrm{H}}}

\def \cm {\mathcal{M}}
\def \fz {\mathfrak{z}}

\def \cp {\mathcal{P}}
\def \cc {\mathcal{C}}
\def \ch {\mathcal{H}}
\def \ce {\mathcal{E}}

\def \cu {\mathcal{U}}
\def \co {\mathcal{O}}
\def \cn {\mathcal{N}}
\def \ck {\mathcal{K}}
\def \cj {\mathcal{J}}
\def \ca {\mathcal{A}}

\def \cz {\mathcal{Z}}
\def \ci {\mathcal{I}}
\def \cw {\mathcal{W}}

\newcommand{\dem}{\begin{proof}}
\newcommand{\cqd}{\end{proof}}

\newcommand{\leb}{\operatorname{Leb}}
\newcommand{\dist}{\operatorname{dist}}
\newcommand{\ord}{\operatorname{ord}}
\newcommand{\supp}{\operatorname{supp}}

\newcommand{\diameter}{\operatorname{diameter}}

\begin{document}

\author{Vilton  Pinheiro}
\address{Departamento de Matem\'atica, Universidade Federal da Bahia\\
Av. Ademar de Barros s/n, 40170-110 Salvador, Brazil.}
\email{viltonj@ufba.br}

\date{\today}

\thanks{Work carried out at the  Federal University of
Bahia, and IMPA. Partially supported by CNPq, PROCAD/CAPES, IMPA and PRONEX - Dynamical Systems}

\title{Expanding Measures}

\maketitle

\begin{abstract}
We prove that any $C^{1+\alpha}$ transformation, possibly with a
(non-flat) critical or singular region, admits an invariant
probability measure absolutely continuous with respect to any
expanding measure whose Jacobian satisfies a mild distortion
condition. This is an extension to arbitrary dimension of a
famous theorem of Keller \cite{Ke} for maps of the interval
with negative Schwarzian derivative.

We also show how to construct an induced Markov map $F$
such that every expanding probability of the initial transformation
lifts to an invariant probability of $F$. The induced time is bounded at each point by the corresponding first
hyperbolic time (the first time the dynamics exhibits hyperbolic behavior).
In particular, $F$ may be used to study decay of correlations and others statistical properties
of the initial map, relative to any expanding probability.
\end{abstract}

\tableofcontents

%%%%%%%%%%%%%%%%%%%%%%%%%%%%%%%%%%%%%%%%%%%%%%%%%%%%%%%%%%%%%%%%%%%%%%%%%%%%%%

\section{Introduction}

In this work we propose a general construction of Markov structures
for non-uniformly expanding transformations. A distinctive feature is that these
Markov structures capture all trajectories with expanding behavior.

In particular, we are able to use them to prove existence of ergodic invariant measures absolutely continuous with respect to any expanding reference measure with Hölder continuous Jacobian. In the special case when Lebesgue measure is the reference, this yields the physical measures of the transformation. Our Markov structures open the way for further development of the ergodic theory of this class of systems. In this direction, we construct Markov transformations induced from the original one, and we prove that any expanding invariant measure of the initial map lifts to invariant measure of these Markov transformations.

Markov partitions were the principal tool for analyzing the qualitative behavior of uniformly hyperbolic (Axiom A) or even uniformly expanding systems (see \cite{Shub}). For uniformly hyperbolic dynamics, the  systematic introduction of these partitions was due to Sinai \cite{Sin68a,Sin68b,Si} and Bowen \cite{Bow70a,Bow70b} and became a key technical tool in the ergodic theory of uniformly hyperbolic/expanding systems (see \cite{Bow75}).
Sinai, Ruelle, Bowen used Markov partitions to associate
these dynamical systems with symbolic ones, prove existence and uniqueness of equilibrium states, and several other properties, in a neighborhood of every transitive hyperbolic set.
Recall that a Markov partition for a map $f:\Lambda\to \Lambda$ is a cover $\cp=\{P_1,\cdots,P_s\}$ of $\Lambda$ satisfying
\begin{itemize}
\item[(a)] $\interior P_i\cap\interior P_j=\emptyset$ if $i\ne j$;
\item[(b)] if $f(P_j)\cap\interior P_i\ne\emptyset$ then $f(P_j)\supset P_i$.
\end{itemize}

Our setting is much more general than the classical family of uniformly expanding maps. Indeed we assume our systems to be non-uniformly expanding.
In this setting one can not expect, in general, the existence of a classical finite Markov partition as there exist parts of the system that spend arbitrarily large time to present some expanding behavior. Nevertheless, we shall prove the existence of a quite similar partition that will be called an {\em induced Markov partition}.
An {\em induced Markov partition}  is an at most countable cover $\cp=\{P_1,P_2,P_3,\cdots\}$ of $\Lambda$ satisfying
\begin{itemize}
\item[(a)] $\interior P_i\cap\interior P_j=\emptyset$ if $i\ne j$;
\item [(b)] for each $P_j$ there is a $R_j\ge1$ such that
\begin{itemize}
\item[(b.1)] if $\ell<R_j$ and $f^{\ell}(P_j)\cap\interior P_i\ne\emptyset$ then $f^{\ell}(P_j)\subset P_i$;
\item[(b.2)] if $f^{R_j}(P_j)\cap\interior P_i\ne\emptyset$ then $f^{R_j}(P_j)\supset P_i$.
\end{itemize}
\end{itemize}

Let us be more precise about the kind of systems we will deal with in this paper.
Formal statements will appear later.
Let $f:M\to M$ be a $C^{1+\alpha}$ transformation outside some critical/singular set $\cc\subset M$ (the case $\cc=\emptyset$ is a possibility).
A positively invariant set $\ch\subset M$ is called {\em expanding} if every point $x\in\ch$ satisfies
\begin{equation}\label{EquationExpanding0}
\limsup_{n\to\infty}\frac{1}{n}\sum_{i=0}^{n-1}
\log \|(Df(f^{i}(x)))^{-1}\|^{-1}>0
\end{equation}
and if $\ch$ satisfies the condition of slow approximation to the critical set, i.e.,
for each $\varepsilon>0$ there is a $\delta>0$
such that
\begin{equation} \label{EquationFaraway1}
\limsup_{n\to+\infty}
\frac{1}{n} \sum_{j=0}^{n-1}-\log \mbox{dist}_{\delta}(f^j(x),\cc)
\le\varepsilon
\end{equation}
for every $x\in\ch$,
where $\dist_{\delta}(x,\cc)$ denote the $\delta$-\emph{truncated} distance from $x$ to $\cc$ defined as
$\dist_{\delta}(x,\cc)=\dist(x,\cc)$ if $\dist(x,\cc) \leq \delta$
and $\dist_{\delta}(x,\cc) =1$ otherwise.

A probability measure is called {\em expanding} if there is an expanding set $\ch$ such that $\mu(\ch)=1$.
If $f$ is a $C^{1+\alpha}$ endomorphism then any invariant measure satisfying (\ref{EquationExpanding0}) almost everywhere is automatically an expanding measure (Corollary~\ref{CorollaryNonFlatToSlowRecurrence}).

Given any expanding set $\ch$, we construct an induced Markov partition $\cp$ of $\ch$ with respect to $f$ (or an iterate of it). Associated to this partition there is an induced map
$$F:\Lambda\to\Lambda,\,\,\,\,F(x)=f^{R(x)}(x),$$
which is Markov, with an appropriate upper bound on the inducing time.

Given any reference measure $\nu$ which gives positive weight to $\ch$, we can use the induced Markov map to construct $f$-invariant probabilities absolutely continuous with respect to $\nu$, and study decay of correlations and others statistical properties.

A crucial point to be noted is that every $f$-invariant measure $\mu$ that gives positive weight to $\ch$ can be lifted to the level of the induced map (the induced map does not depend on the measure $\mu$).

We also give several examples of expanding measures and applications of these results.

\subsection{Statement of main results}
\label{SettingAndStatementOfMainResults}

Let $M$ be a compact Riemannian manifold of dimension $d\ge1$ and $f:M\to M$ a map defined on $M$.

The map $f$ is called {\em non-flat} if it is a local
$C^{1+}$ ( i.e., $C^{1+\alpha}$ with $\alpha>0$ ) diffeomorphism in the whole manifold  except  in  a {\em
non-degenerate critical/singular set}  $\mathcal{C}\subset {M}$. We
say that $\mathcal{C}\subset {M}$ is a {\em non-degenerate critical/singular set} if $\exists \beta,B>0$ such that
the following two conditions hold.
\begin{enumerate}
\item[(C.1)]
\quad $\displaystyle{\frac{1}{B}dist(x,\mathcal{C})^{\beta}\leq
\frac {\|Df(x)v\|}{\|v\|}\leq B\,dist(x,\mathcal{C})^{-\beta} }$ for all $v\in T_x {M}$.
\end{enumerate}
For every $x,y\in {M}\setminus\mathcal{C}$ with
$dist(x,y)<dist(x,\mathcal{C})/2$ we have
\begin{enumerate}
\item[(C.2)] \quad $\displaystyle{\left|\log\|Df(x)^{-1}\|-
\log\|Df(y)^{-1}\|\:\right|\leq
\frac{B}{dist(x,\mathcal{C})^{\beta}}dist(x,y)}$.
\end{enumerate}

If $dim({M})=1$ and $f$
satisfies the usual one dimensional definition of non-flatness (see \cite{MvS}),
then it also satisfies the definition given above.

In the whole paper, a measure will be a countable additive measure defined on the Borel sets.
A measure $\mu$ is called {\em $f$-non-singular} if
$f_*\mu\ll\mu$, where $f_*\mu$ ($=\mu\circ f^{-1}$) is the push-forward of
$\mu$ by $f$.

Let $f$ be a non-flat map with critical/singular set $\cc\subset M$.
A finite measure $\mu$ is called {\em $f$-non-flat} if it is $f$-non-singular,
$\mu(\cc)=0$, $J_{\mu}f(x)$ is
well defined and positive for $\mu$-almost every $x\in M$, and
for $\mu$-almost every $x,y\in {M}\setminus\mathcal{C}$ with
$dist(x,y)<dist(x,\mathcal{C})/2$ we have
$$
\displaystyle{\left|\log\frac{J_\mu f(x)}{J_\mu f(y)}\:\right|\leq
\frac{B}{\dist(x,\mathcal{C})^{\beta}}\dist(x,y)}.
$$

\subsubsection{Expanding sets and measures}\label{SussectionEXPsetsANDmeas}
%Let $\dist_{\delta}(x,\cc)$ denote the $\delta$-\emph{truncated} distance from $x$ to $\cc$ defined as
%$\dist_{\delta}(x,\cc)=\dist(x,\cc)$ if $\dist(x,\cc) \leq \delta$
%and $\dist_{\delta}(x,\cc) =1$ otherwise.

\begin{Definition}\label{DefinitionExpandingSet}A positively invariant set $\ch\subset{M}$ (i.e., $f(\ch)\subset\ch$) is called {\em $\lambda$-expanding}, $\lambda\ge 0$, if
\begin{equation}\label{EquationExpanding}
\limsup_{n\to\infty}\frac{1}{n}\sum_{i=0}^{n-1}
\log \|(Df(f^{i}(x)))^{-1}\|^{-1}>\lambda,
\end{equation}
for every $x\in \ch$, and $\ch$ satisfies the slow approximation condition, i.e.,
for each $\varepsilon>0$ there is a $\delta>0$
such that (\ref{EquationFaraway1}) holds for every $x\in\ch$.
\end{Definition}

An expanding set is a positively invariant set but, in general, it is not a compact one.
In the one-dimensional case (\ref{EquationExpanding})
reduces to the Lyapunov exponent of $f$ on $x$ to be bigger than $\lambda$, i.e.,
$$
\limsup_{n\to\infty}\frac{1}{n}\sum_{i=0}^{n-1}
\log |f'(f^{n}(x))|= \limsup_{n\to\infty} |(f^{n})'(x)|>\lambda.
$$

\begin{Definition}[Expanding measures]
We call a measure $\mu$ (non necessarily invariant) a {\em $\lambda$-expanding} measure (with respect to $f$) if $\mu$ is $f$-non-singular and there exists a $\lambda$-expanding set
$\ch$ such that $\mu(M\setminus \ch)=0$.
\end{Definition}

\begin{maintheorem}[Existence of absolutely continuous invariant  measures]\label{TheoremSRB}Let $f:M\to M$ be a non-flat map.
If $\mu$ is a $f$-non-flat $\lambda$-expanding measure, $\lambda>0$,
then there exists a finite collection of $\mu$ absolutely continuous ergodic $f$-invariant probabilities such that $\mu$-almost every
point in $M$ belongs to the basin of one of these probabilities.
\end{maintheorem}

Recall that the {\em basin of measure} $\eta$ is the set $\mathcal{B}(\eta)$ of the points $x\in M$ such that
$$\lim_{n\to +\infty}\frac{1}{n}\sum_{j=0}^{n-1}\varphi\circ f^j(x)=\int\varphi d\eta,$$
for every continuous function $\varphi:M\to\RR$.

\subsubsection{Markov partitions}

Let $f:U\to U$ a measurable map defined on a Borel set $U$ of a compact, separable metric space $X$. A countable collection $\cp=\{P_1,P_2,P_3...\}$ of
Borel subsets of $U$ is called a {\em Markov partition} if
\begin{enumerate}
\item $\interior(P_i)\cap \interior(P_j)=\emptyset$ if $i\ne j$;
\item if $f(P_i)\cap\interior(P_j)\ne\emptyset$ then $f(P_i)\supset\interior(P_j)$;
\item $\#\{f(P_i)$ ${;}$ $i\in\NN\}<\infty$;
\item $f|_{P_i}$ is a homeomorphism and it can be extended to a homeomorphism sending  $\overline{P_i}$ onto $\overline{f(P_i)}$;
\item $\lim_n\mbox{diameter}(\cp_n(x))=0$ $\forall x\in\bigcap_{n\ge0}f^{-n}\big(\bigcup_{i}P_i\big)$,
\end{enumerate}
where $\cp_n(x)$ $=$
$\{y$ ${;}$ $\cp(f^j(y))=\cp(f^j(x))$ $\forall\,0\le j\le n\}$ and $\cp(x)$ denotes the element of $\cp$ that  contains $x$.

\begin{Definition}[Induced Markov partition]
A countable collection $\cp=\{P_1,P_2,P_3...\}$ of
Borel subsets of $U$ is called a {\em induced Markov partition} if it satisfies all conditions of a Markov partition except the second one which has to be substituted by the following
\begin{enumerate}
\item[(2)] for each $P_i\in\cp$ there is a $R_i\ge1$ such that
\begin{enumerate}
\item[(2.1)]  if $\ell<R_i$ and $\interior(f^{n}(P_i))\cap\interior(P_j)\ne\emptyset$ then $\interior(f^{n}(P_i))$ $\subset$ $\interior(P_j)$ or $\interior(f^{n}(P_i))$ $\supset$ $\interior(P_j)$;
\item[(2.2)]  if $f^{R_i}(P_i)\cap\interior(P_j)\ne\emptyset$ then $f^{R_i}(P_i)\supset \interior(P_j)$.
\end{enumerate}
\end{enumerate}
\end{Definition}

\begin{Definition}[Markov map] The pair $(F,\cp)$, where $\cp$ is a Markov partition of $F:U\to U$, is called a {\em Markov map} defined on $U$. If $F(P)=U$ $\forall\,P\in\cp$,
$(F,\cp)$ is called a {\em full Markov map}.
\end{Definition}

Note that if $(F,\cp)$ is a full Markov map defined on an open set $U$ then the elements of $\cp$ are open sets (because $F(P)=U$ and $F|_P$ is a homeomorphism $\forall\,P\in\cp$).

Consider a measurable map $f:M\to M$ from $M$ to $M$ (or, more in general, from the metric space $X$ to $X$).

\begin{Definition}[Induced Markov map]\label{DefinitionInducdMMap}
A Markov map  $(F,\cp)$ defined on $U$ is called a {\em induced Markov map} for $f$ on $U$ if
is there is a function $R:U\to\NN=\{0,1,2,3,...\}$ (called {\em inducing time}) such that
$\{R\ge 1\}=\bigcup_{P\in\cp}P$, $R|_P$ is constant $\forall\,P\in\cp$
and $F(x)=f^{R(x)}(x)$ $\forall\,x\in U$.
\end{Definition}

If an induced Markov map $(F,\cp)$ is a full Markov map, we call $(F,\cp)$ an {\em induced full Markov map}.

Given an induced Markov map $(F,\cp)$, an ergodic $f$-invariant probability $\mu$ is said {\em liftable} to $F$ if there exists $F$-invariant finite measure $\nu\ll\mu$ such that $$\mu=\sum_{P\in\cp}\sum_{j=0}^{R(P)-1}f_{\ast}^j\left(\nu|_{P}\right),$$ where $R$ is the inducing time of $F$, $\nu|_{P}$  denotes the measure given by $\nu|_{P}(A)=\nu(A\cap P)$ and $f_*^j$ is the
push-forward by $f^j$.

\begin{Definition}[Markov structure]
A {\em Markov structure} for a set $U\subset M$ (or $X$) is an at most countable collection $\mathfrak{F}=\{(F_i,\cp_i)\}_i$ of induced Markov maps such that
if $\mu$ is an ergodic $f$-invariant probability with $\mu(U)=1$ then $\exists(F_i,\cp_i)\in\mathfrak{F}$ such that $\mu$ is liftable to $F_i$.
\end{Definition}

\begin{maintheorem}[Markov structure for an expanding set]
 \label{TeoremMarkovStructureForExpanding}
 Every $\lambda$-expanding set with $\lambda>0$ admits a finite Markov structure $\mathfrak{F}=\{(F_1,\cp_1),...,(F_s,\cp_s)\}$. Furthermore, each $(F_i,\cp_i)\in\mathfrak{F}$ is a full Markov map defined on some topological open ball $U_i$ (also the elements of $\cp_i$ are topological open balls).
\end{maintheorem}

See Theorem~\ref{TeoremExpandingMarkovStructure} for a more complete version of this Theorem~\ref{TeoremMarkovStructureForExpanding}.

As an expanding set $\ch$ is positively invariant, it follows from Theorem~\ref{TeoremMarkovStructureForExpanding} that every ergodic $f$-invariant probability having $\mu(\ch)>0$ is liftable to one of the full induced Markov map $F_1,...,F_s$ given by the Theorem~\ref{TeoremMarkovStructureForExpanding}.

\label{DefinitionEssentiallyOpenSets}
We say that a Borel set $A$ is an {\em essentially open} set if\, $\overline{\mbox{interior}(A)}\supset A$, that is, the closure of the interior of $A$ contains $A$.
The theorem below shows that any expanding set $\ch$ admits a Markov structure composed by a single induced Markov map $(F,\cp)$ with $\cp$ being a collection of essentially open sets.

\begin{maintheorem}[Global expanding induced Markov map]\label{GlobalExpandingInducedMarkovMap}
If $\ch$ is a $\lambda$-expanding set with $\lambda>0$ then there exist $\ell\ge1$ and $\varepsilon_0>0$ such that for every $0<\varepsilon<\varepsilon_0$ there are an induced Markov map $(F,\cp)$, with inducing time $R$, and a finite partition $\cp_0$ of $M$ by essentially open sets satisfying the following conditions.
\begin{enumerate}[(i)]
\item\label{ItemE-1} $\sup\{\diameter(P)$ $;$ $P\in\cp\}\le\max\{\diameter(P)$ $;$ $P\in\cp_0\}<\varepsilon$.
\item\label{ItemE0} For each $Q\in\cp$ there exists $P\in\cp_0$ such that $\interior(Q)\subset P$.
\item\label{ItemE1} $F(P)\in\cp_0$ $\forall\,P\in\cp$(in particular, the elements of $P$ are essentially open sets).
\item\label{ItemE2} $\dist(F(x),F(y))\ge 8\dist(x,y)$ $\forall\,x,y\in P$ and $\forall\,P\in\cp$;
\item\label{ItemE3} For all $x,y\in P$, $P\in\cp$ and $0\le n\le R(P)$,
     $$\dist(f^n(x),f^n(y))\le\big(e^{-\lambda/8}\big)^{(R(P)-n)}\dist(F(x),F(y)).$$
\item\label{ItemE4} If $\mu$ is $f$-non-flat measure (not necessarily
        invariant) with $\mu(M\setminus\ch)=0$ then
        $\exists\,\rho>0$ such that
        $$\left|\log\frac{J_\mu F(x)}{J_\mu F(y)}\right| \le \rho \dist(F(x),
        F(y)),$$ for $\mu$ almost every $x, y\in P$ and $\forall\,P\in\cp$.
\item\label{ItemE5} If $\mu$ is an ergodic $f$-invariant measure with $\mu(\ch)>0$ then  $\mu\big|_{\bigcap_{j\ge0}f^{-\ell\,j}(\{R>0\})}$ is an invariant measure with respect to $f^{\ell}$.
\item\label{ItemE6}There is a good relation between the tail of the partition and the tail of the {\em hyperbolic times} (see Section~\ref{ApplicationsExpandingMeasures} for the definition), i.e.,
$$
\{R>n\}\cap\ch\subset\ch\setminus\bigcup_{j=1}^n\th_{\ell\,j},
$$
where $\th_i$ denotes the set of points having $i$ as a $(e^{-\lambda/4},\varepsilon)$-hyperbolic time.
\item\label{ItemE7} Every ergodic $f$-invariant probability $\mu$
with $\mu(\ch)>0$ is liftable to $F$.
\item\label{ItemUnifExpanding} If
$\ch$ is an uniformly expanding set
(i.e., $\overline{\ch}=\ch$, $\ch\cap\cc=\emptyset$ and $\exists\,n\ge1$ such that $\big\|\big(Df^n(x)\big)^{-1}\big\|^{-1}>1$ $\forall\,x\in\ch$)
then the partition $\cp$ is finite.
\end{enumerate}
\end{maintheorem}

We want to observe that in Section~\ref{SubsectionZoomingSetsAndMeasures} we
introduce the zooming sets (that generalizes the expanding sets)
and the theorems above are corollaries of
Theorems~\ref{TheoremZoomingInvariant},~\ref{TeoremMarkovStructure}
and \ref{GlobalMarkovStructureForZooming} for zooming sets.

A second note is that, if $f$ is a $C^{1+}$ endomorphism, $\cp$ is an induced Markov partition of $\ch$, with respect to $f^{\ell}$, and the estimate of the item (\ref{ItemE6}) of Theorem~\ref{GlobalExpandingInducedMarkovMap}
can be rewritten in a more direct
dependence on the expansion (\ref{EquationExpanding}) and on the
recurrence to the critical region (\ref{EquationFaraway1}) (see
Theorem~\ref{GlobalExpandingInducedMarkovMapForMapsWithBV}).

\subsection{Overview of the paper}

In Section~\ref{SectionNestedSets} we introduced
the notion of {\em nested sets} adapted to the kind of pre-images we want to deal
with (for example, pre-images with some contraction).

In Section~\ref{SectionErgodicComponents} we study the ergodic components for non (necessarily) invariant measures for maps on metric spaces.

In Section~\ref{SectionIndMarMap} we obtain a statistical characterization of the liftable measures for a given
induced map.

Although we are basically interested in expanding measures
(Section~\ref{SettingAndStatementOfMainResults}), we weakened the
expansion condition to permit more flexibility in the applications.
For this we introduce the {\em zooming measures} in
Section~\ref{SubsectionZoomingSetsAndMeasures}.

In Section~\ref{SectionConstructingLocalInducingMap} and \ref{SectionGlobalInducedMarkovMap} we show
most of the results
for zooming sets and measures. In particular the existence of induced Markovian maps for zooming sets (Theorem~\ref{TeoremMarkovStructure} and \ref{GlobalMarkovStructureForZooming})
and the existence of an invariant measure $\nu{\ll}\mu$ that is absolutely continuous with respect to a given zooming measure with some distortion control
(Theorem~\ref{TheoremZoomingInvariant}).

Section~\ref{ApplicationsExpandingMeasures} is dedicated to the definition and properties of expanding measures,
as well as to establish the connection between these measures with the zooming ones.
The existence of an absolutely continuous
invariant measure for a given expanding measure, the induced Markovian maps for expanding sets and so on are consequences of the analogous result for zooming measures and in this section we use the zooming results to get the expanding ones.

In Section~\ref{SectionExamplesAndAplications} we give  many examples of expanding and zooming sets and measures.
We give also some applications of the results of the previous sections. In particular, we study the decay of correlations for general expanding measures.

{\em Acknowledgments.} We thank  V. Araújo, P. Brandão A. Castro and K. Oliveira for comments and for useful conversations. We thanks also IMPA, Brazil, and Universidade do Porto, Portugal, for the hospitality (specially J. Alves, M. F. Carvalho and J. Rocha) and the opportunity to present and discuss this work. We are specially grateful to P. Varandas, M. Todd, M. Viana and J. Palis not only for useful conversations but also for the incentive.

%%%%%%%%%%%%%%%%%%%%%%%%%%%%%%%%%%%%%%%%%%%%%%%%%%%%%%%%%%%%%%%%%%%%%%%%%%%%%%

\section{Nested sets}\label{SectionNestedSets}

The notion of nice interval, introduced by Martens in \cite{Ma},
is a useful tool in the theory of real and complex one-dimensional dynamical
systems (see, for instance, \cite{MvS,PrRLe}).
A \emph{nice interval} is an open interval $I$ such that
the forward orbit $\mathcal O^+ (\partial I)$ of the boundary of
$I$ does not return to I, i.e., $\mathcal O^+ (\partial I) \cap I
=\emptyset$. Note that nice intervals are natural and easy to
construct for interval maps. For instance, two consecutive points
of a periodic orbit define a nice interval. Its main property is
that there are no linked pre-images of a nice interval, that is,
if $I_1$ and $I_2$ are sent homeomorphically onto an open nice
interval $I$ by $f^{n_1}$ and $f^{n_2}$ respectively then either
$I_1 \cap I_2 =\emptyset$, $I_1 \subset I_2$ or $I_2 \subset I_1$.

In the multidimensional case, the boundary of topological open balls are
connect topological manifolds and if a chaotic transitive dynamic is not much symmetric,
it is natural to expect that this dynamic will spread these boundaries to the whole manifold,
forbidding any ``nice ball''. In general, the same seems true for any set whose the boundary is not totally disconnected.

In this section
we present the abstract construction of nested sets. This reformulates and generalizes the concept of {\em nice interval}. In
Section~\ref{SubsectionZoomingSetsAndMeasures} we show
their abundance in the presence of some expansion
(see Lemma~\ref{LemmaZoomingNestedBall}).

Let $f:X\to X$ be a map defined on a complete, separable metric space $X$.
Fixed some $K\subset {X}$, a set $P\subset X$ is called a {\em regular pre-image of
order} $n\in\NN$ of $K$ if $f^n$ sends $P$ homeomorphically onto
$K$. Denote the order of $P$ (with respect to $K$) by $\ord(P)$.

Let us fix in all Section~\ref{SectionNestedSets} a collection $\ce_0$
of connected open subsets of ${X}$ (for instance, $\ce_0$ can be the collection
$\{f^n(V_n(x))$ ${;}$ $x\in {\th}_n$ and $n\in\NN\}$ of all hyperbolic balls of ${X}$,
see Proposition~\ref{PropositionHyperbolicBalls}). For each $n\in\NN$ and $V\in\ce_0$
consider some collection $\ce_n(V)$ of regular pre-images of order $n$ of $K$.
Set $\ce_n=(\ce_n(V))_{V\in\ce_0}$.
We call the sequence $\ce=(\ce_n)_n$ a {\em dynamically
closed  family of (regular) pre-images} if
$f^\ell(E)\in\ce_{n-\ell}$ $\forall\,E\in\ce_n$ and $\forall\,0\le
\ell\le n$. Given $Q\in\ce_n$  we denote $f^n|_Q$ by $f^Q$ and we
denote the $\ce$-inverse branch of  associated to $Q$,
$(f^n|_Q)^{-1}$, by $f^{-Q}$.

Let $\ce=(\ce_n)_n$ be a dynamically closed  family of pre-images.
A set $P$ is called an
$\ce$-pre-image of a set $W\subset {X}$ if  there is $n\in\NN$
and $Q\in\ce_n$ such that $W\subset f^n(Q)$ and $P=f^{-Q}(W)$.

\begin{Remark}\label{RemarkPreImagensNaoLinkadas}
Two distinct $\ce$-pre-images $\mathcal{X}_1$ and $\mathcal{X}_2$ of some set
$\mathcal{X}\subset {X}$ having the same order cannot intersect.
Indeed, write $n=\ord(\mathcal{X}_1)=\ord(\mathcal{X}_2)$ and for each $i\in\{1,2\}$
write $\mathcal{X}_i=f^{-Q_i}(\mathcal{X})$, with $Q_i\in\ce_{n}$. Let $P_j=f^{-Q_j}(Q_1\cap Q_2)$,
for $j=1,2$. It follows that $P_1\cap P_2\supset \mathcal{X}_1\cap \mathcal{X}_2\ne\emptyset$.
Of course $P_1\ne P_2$, otherwise $\mathcal{X}_2$ $=$ $f^{-Q_2}(\mathcal{X})$ $=$ $(f^n|_{P_2})^{-1}(\mathcal{X})$ $=$ $(f^n|_{P_1})^{-1}(\mathcal{X})$ $=$ $f^{-Q_1}(\mathcal{X})$ $=$ $\mathcal{X}_1$.
Thus $P_1\cap\partial P_2\ne\emptyset$ or $P_2\cap\partial P_1\ne\emptyset$.
Assume that $P_1\cap\partial P_2\ne\emptyset$. So, $\emptyset\ne f^n(P_1\cap\partial P_2)\subset f^n(P_1)\cap\partial \big(f^n(P_2)\big)\subset (Q_1\cap Q_2)\cap\partial \big(Q_1\cap Q_2\big)=\emptyset$. An absurd.
\end{Remark}

\begin{Definition}[Linked sets] We  say that two open sets $U_1$ and $U_2$ are {\em linked}
if both $U_1\setminus U_2$ and $U_2\setminus U_1$ are not empty sets.
\end{Definition}

Note that two connected open sets $U_1$ and $U_2$ are
linked if and only if $\partial U_1\cap U_2$ and $U_1\cap\partial U_2$ are not empty sets.

\label{SubsectionNestedSetsAndCollections}

\begin{Definition}[$\ce$-nested set] A set $V$ is called {\em $\ce$-nested} if $V$ is an open set and $V$ is not linked
with any $\ce$-pre-image of $V$.
\end{Definition}

The fundamental property of a nested set is that
any  $\ce$-pre-images $P_1$ and $P_2$ of it are not linked (see Corollary~\ref{CorollaryMainNestedProperty}).

We can extend the concept of $\ce$-nested set to a collection of sets in the following way.
\begin{Definition}[$\ce$-nested collection of sets]\label{DefinitionNestedCollection}
A  collection $\ca$ of open sets
is called an {\em $\ce$-nested collection of sets} if every
$A\in\ca$ is not linked with any $\ce$-pre-image of an element of $\ca$ with order bigger than zero.
Precisely, if $A_1\in\ca$ and $P$ is an $\ce$-pre-image of some $A_2\in\ca$,
then either $A_1$ and $P$ are not linked or $P=A_2$.
\end{Definition}

It follows from the definition of an $\ce$-nested collection of sets that every sub-collection of
an $\ce$-nested collection is also an $\ce$-nested collection. In particular, each element of an
$\ce$-nested collection is an $\ce$-nested set.

\begin{Lemma}[Main property of a nested collection]\label{LemmaPreMainNestedProperty}
If $\ca$ is an $\ce$-nested collection of open sets and $P_1$ and $P_2$
are $\ce$-pre-images of two elements of $\ca$ with $\ord(P_1)\ne\ord(P_2)$
then $P_1$ and $P_2$ are not linked.
\end{Lemma}
\dem Let $\ell_j=\ord(P_j)$ for $j=1,2$.
We may assume that
$\ell_1<\ell_2$ and, by contradiction, assume that $P_1$  and $P_2$
are linked. Let, for $i=1,2$, $p_j\in P_j\cap \partial P_{3-j}$,
$Q_i\in\ce_{\ell_i}$  and $A_i\in\ca$ be such that $P_i=f^{-Q_i}(A_i)$.  As $\ce$ is a
dynamically closed family of pre-images of elements of $\ce_0$,
$Q=f^{\ell_1}(Q_2)\in\ce_{\ell_2-\ell_1}$ and
$P=f^{\ell_1}(P_2)=f^{-Q}(A_2)$ is an $\ce$-pre-image of $A_2$. On the other hand $f^{\ell_1}(P_1)=A_1\in\ca$.
As
$f^{\ell_1}(p_1)\in f^{\ell_1}(P_1)\cap \partial\big(f^{\ell_1}(P_{2})\big)=A_1\cap\partial P$
and
$f^{\ell_1}(p_2)\in f^{\ell_1}(P_2)\cap \partial\big(f^{\ell_1}(P_{1})\big)=P\cap\partial A_1$,
it follows that $P$ and $A_1$ are linked, but this is impossible because $\ca$ is $\ce$-nested.

\cqd

\begin{Corollary}[Main property of a nested set]\label{CorollaryMainNestedProperty}
If $V$ is an $\ce$-nested set and $P_1$ and $P_2$
are $\ce$-pre-images  of $V$ then $P_1$ and $P_2$ are not linked.
Furthermore,
\begin{enumerate}
\item if $P_1\cap P_2\ne\emptyset$ then $\ord(P_1)\ne\ord(P_2)$;
\item if $P_1\subsetneqq P_2$ with $\ord(P_1)<\ord(P_2)$ then
$V$ is contained in an $\ce$-pre-image of itself with order bigger than zero,
$$f^{\,\ord(P_2)-\ord(P_1)}(V)\subset V.$$
\end{enumerate}
\end{Corollary}
\dem Lets suppose that  $P_1\ne P_2$ are $\ce$-pre-images of $V$ and set
$\ell_j=\ord(P_j)$ for $j=1,2$.
By Remark~\ref{RemarkPreImagensNaoLinkadas}, $\ell_1\ne\ell_2$. Thus, we may assume that $\ell_1<\ell_2$.
By Lemma~\ref{LemmaPreMainNestedProperty} it follows that $P_1$ and $P_2$ are not linked.

Now, suppose in addition that $P_1\subset P_2$.
Then $V=f^{\ell_1}(P_1)\subset
f^{\ell_1}(P_2)$ ($f^{\ell_1}(P_2)$ is an $\ce$-pre-image of $V$) and this will imply
that $f^{\,\ell_2-\ell_1}(V)\subset f^{\ell_2}(P_2)=V$.
\cqd

\subsection{Constructing nested sets}\label{SubsectionConstrNestedSets}

In this section (Section~\ref{SubsectionConstrNestedSets}) let $\ca$ be a collection of connected open sets such that the elements of $\ca$ are not contained in any $\ce$-pre-image of order bigger than zero of an element of $\ca$.

%%%%%%%%%%%%%%%%%%%%%%%%%%%%%%%%%%%%%%%%%%%%%%
\begin{figure}
  \includegraphics{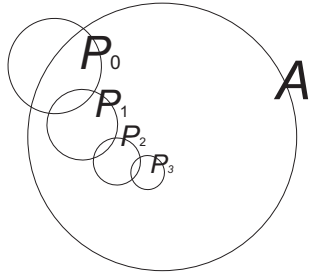}\\
  \caption{A chain $(P_0,P_1,P_2,P_3)$ of pre-images beginning in $A$.}\label{FigureChains}
\end{figure}
%%%%%%%%%%%%%%%%%%%%%%%%%%%%%%%%%%%%%%%%%%%%%%

A finite sequence $\ck=(P_0,P_1,\cdots,P_n)$ of $\ce$-pre-images of
elements of $\ca$ is called a {\em chain of pre-images of $\ca$} beginning in $A\in\ca$
(Figure~\ref{FigureChains}) if
\begin{enumerate}
\item $0<\ord(P_0)\le\cdots\le\ord(P_{n-1})\le\ord(P_{n})$;
\item $A$ and $P_0$ are linked;
\item $P_{j-1}$ and $P_{j}$ are linked $\forall 1\le j\le n$;
\item $P_i\ne P_j$ $\forall\,i\ne j$.
\end{enumerate}
Denote by $ch_{\ce}(A)$ the collection of all chain of pre-images of $\ca$ beginning in $A\in\ca$.
As the elements of $\ca$ are connected and open, it is easy to check the following remark.
\begin{Remark}\label{RemarkConnectedChains}
If $(P_0,P_1,\cdots,P_n)\in ch_{\ce}(A)$, with $A\in\ca$, then $\bigcup_{j=n_0}^{n_1}P_j$ is
a connected open set $\forall\,0\le n_0\le n_1\le n$.
\end{Remark}

For each $A\in\ca$ define the open set
\begin{equation}\label{EquationStar}A^{\star}=A\setminus\,
\overline{\bigcup_{(P_j)_j\in ch_{\ce}(A)}\bigcup_j
P_j}.\end{equation}

%%%%%%%%%%%%%%%%%%%%%%%%%%%%%%%%%%%%%%%%%%%%%%
\begin{figure}
\includegraphics{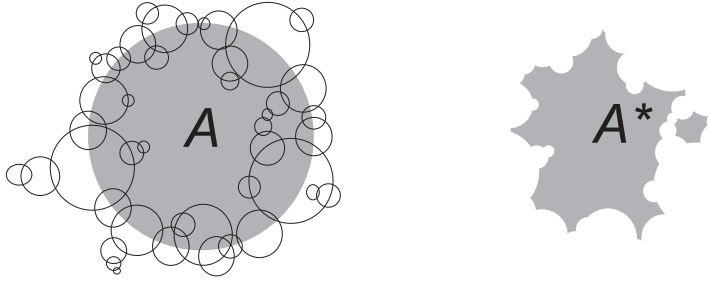}\\
\caption{On the left side it is shown a ball $A$ (in grey) and the boundaries of the pre-images of $A$ that belong to the chains. On the right side $A^{\star}$ is shown.}\label{FigureNestedSet}
\end{figure}
%%%%%%%%%%%%%%%%%%%%%%%%%%%%%%%%%%%%%%%%%%%%%%

\begin{Proposition}[An abstract construction of a nested collection]\label{PropositionExisteNested2}
For each $A\in\ca$ such that $A^\star\ne\emptyset$ choose a connected
component $A'$ of $A^\star$. If $\ca'=\{A'$ ${;}$ $A\in\ca$ and $A^{\star}\ne\emptyset\}$ is
not an empty collection  then $\ca'$ is an $\ce$-nested collection of sets.
\end{Proposition}
\dem
Suppose that $\ca'\ne\emptyset$. By contradiction, assume that there exist $A_1,A_2\in\ca$ and a $\ce$-pre-image $P$ of $A_2'$ such that $A_1'$ and $P$ are linked. So, as $A_1'$ and $P$ are connected sets,
$\exists\,p\in P\cap\partial A_1'$. Let $\wp=\ord(P)$ and let $E\in\ce_\wp$ be such that
$P=f^{-E}(A_2')$. Setting $Q=f^{-E}(A_2)$, we get $P\subset Q$.

\begin{claim}
$Q\subset A_1$.
\end{claim}
\dem
First note that $Q\cap A_1\supset Q\cap A_1'\supset P\cap A_1'\ne\emptyset$.
On the other hand, if $Q\cap\partial A_1\ne\emptyset$, the unitary
sequence $(Q)$ will be a chain of $\ce$-pre-images beginning in $A_1$, i.e.,
$(Q)\in ch_{\ce}(A_1)$. But this is a contradiction to the definition of
$A_1^{\star}$ because $Q\cap A_1^\star\supset Q\cap A_1'\ne\emptyset$. Thus,
$Q\cap\partial A_1=\emptyset$. As $Q$ and $A_1$ are connected sets and $Q\cap A_1\supset Q\ne\emptyset$,
we get $Q\supset A_1$ or $Q\subset A_1$. The first option is not possible because (by hypothesis) the elements of $\ca$ are not contained in any $\ce$-pre-image of order bigger than zero of an element of $\ca$.
Therefore, $Q\subset A_1$.
\cqd

%%%%%%%%%%%%%%%%%%%%%%%%%%%%%%%%%%%%%%%%%%%%%%
\begin{figure}
  \includegraphics{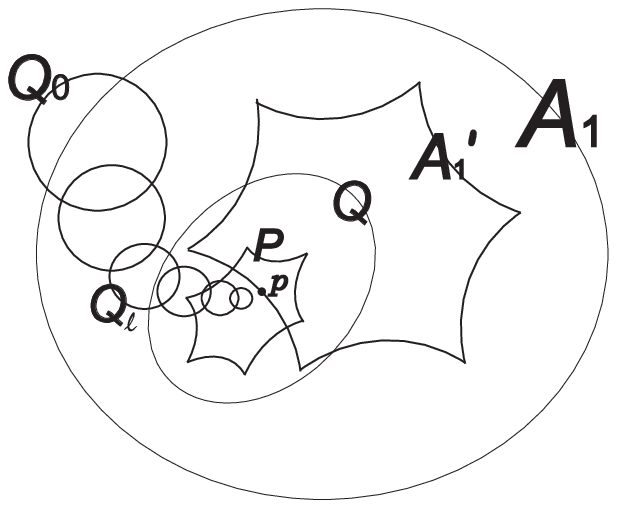}\\
  \caption{}\label{FigureProposicaoSTAR}
\end{figure}
%%%%%%%%%%%%%%%%%%%%%%%%%%%%%%%%%%%%%%%%%%%%%%

As $p\in\partial A_1'$, for a given $\varepsilon>0$ there exists a chain
$(Q_0,\cdots,Q_n)\in ch_{\ce}(A_1)$ such that $\dist(p,\cup_{j=0}^n
Q_j)$ $<$ $\varepsilon$. On the other hand, as $P$ and $Q$ are open sets and $p\in P\subset Q$, taking
$\varepsilon$ small enough, $P\cap(\cup_{j=0}^n Q_j)\ne\emptyset$
and so,
\begin{equation}\label{EquationE222}Q_m\cap Q\supset Q_m\cap
P\ne\emptyset,\end{equation} for some $1\le m\le n$. As
$Q_0\cup\cdots\cup Q_m$ is a connected set (Remark~\ref{RemarkConnectedChains}) and $Q_0\cap ({X}\setminus
Q)\supset Q_0\cap ({X}\setminus
A_1)\ne\emptyset$ (because $Q_0$ and $A_1$ are linked),
there exists $0\le j\le m$ such that $Q_j\cap\partial
Q\ne\emptyset$. Let $\ell=\min\{0\le j\le m$ ${;}$ $Q_j\cap Q\ne\emptyset\}$.

We have two cases, either $\ord(Q_\ell)\le\ord(Q)$
or $\ord(Q_\ell)>\ord(Q)$.
Suppose first that $\ord(Q_\ell)\le\ord(Q)$.
By the minimality of $\ell$, $Q\ne Q_j$ $\forall\,0\le j\le\ell$. Thus, it is easy to check that
$\ck=(Q_0,\cdots,Q_\ell,Q)\in ch_{\ce}(A_1)$.
As $Q\cap A_1^{\star}\supset Q\cap A_1'\ne\emptyset$, the existence of the chain $\ck$ is a contradiction to
(\ref{EquationStar}) and so, this case cannot occur. For the second
case ($\ord(Q_\ell)>\ord(Q)$),
consider the sequence
$\ck=(f^\wp(Q_\ell),\cdots,f^\wp(Q_m))$. It is also easy to check that
$\ck\in ch_{\ce}(A_2)$ (note that, as $f^\wp(Q)=A_2$,
$f^\wp(Q_\ell)\cap\partial A_2=f^\wp(Q_\ell\cap\partial
Q)\ne\emptyset$). But, as
$f^\wp(P)=A_2'\subset A_2^{\star}$, it follows from (\ref{EquationE222}) that $f^\wp(Q_m)\cap
A_2^{\star}\supset f^\wp(Q_m\cap P)\ne\emptyset$, contradicting (\ref{EquationStar}) again and concluding the proof.
\cqd

An easy way to assure the existence of nested sets (or collections) is
to show that the chains have small diameter,
where the diameter of a chain $(P_j)_j$ is defined as the diameter of $\bigcup_j P_j$.

\begin{Corollary}\label{CorollarySmallDiameter}
Let $\varepsilon\in(0,1/2)$ and let $A=B_r(p)$ be a connected open ball with radius $r$ centered in $p\in X$ such that $f^n(A)\not\subset A$ $\forall\,n>0$. If
every chain of $\ce$-pre-images of $A$ has diameter smaller than
$2\varepsilon r$ then the set $A^{\star}$, given by (\ref{EquationStar}) contains the ball $B_{r(1-2\varepsilon)}(p)$. Moreover, the connected component $A'$ of $A^{\star}$ that contains $p$ is a $\ce$-nested set containing $B_{r(1-2\varepsilon)}(p)$.
\end{Corollary}
\dem Set $\ca=\{A\}$. As $f^n(A)\not\subset A$ $\forall\,n>0$,
it follows that $A$ is not contained in any $\ce$-pre-image of itself (with order bigger than zero).
Let $\Gamma$ be the collection of all chains of
$\ce$-pre-images of $A$. If $(P_j)_j\in\Gamma$ then
$\bigcup_j P_j$ is a connected open set intersecting $\partial A$
with diameter smaller than
$2\varepsilon r$. Thus, $\bigcup_j P_j\subset B_{2\varepsilon r}(\partial A)$,
$\forall\,(P_j)_j\in\Gamma$. As a consequence, $A^{\star}=A\setminus\,
\overline{\bigcup_{(P_j)_j\in\Gamma}\bigcup_j P_j}\supset A\setminus
\overline{B_\varepsilon(\partial A)}\supset B_{r(1-2\varepsilon)}(p)$ is a non-empty open set.
Taking $A'$ as the connected component  of $A^{\star}$ that contains $p$ (and so, contains $B_{r(1-2\varepsilon)}(p)$),
it follows from Proposition~\ref{PropositionExisteNested2} that $A'$ is a $\ce$-nested set.\cqd

%%%%%%%%%%%%%%%%%%%%%%%%%%%%%%%%%%%%%%%%%%%%%%%%%%%%%%%%%%%%%%%%%%%%%%%%%%%%%%

\section{Ergodic components}\label{SectionErgodicComponents}

Before constructing the Markov partition using the  adapted {\em nested sets},
we need also some preliminary knowledge of the so called ergodic components for non (necessarily) invariant measures.
This knowledge is important to assure good
statistical properties for these nested sets with respect to the class of measures
that we are working on.

Let $\mu$ be a finite
measure defined on the Borel sets of the compact,
separable metric space ${X}$ and let $f:X\to X$ be a measurable map. A
subset $U\subset {X}$ is called an invariant set (with respect to $f$) if $f^{-1}(U)=U$, and it is called a positively invariant set if $f(U)\subset U$.
\begin{Definition}[Ergodic components] An invariant set $U$ with $\mu
(U)>0$ is called an {\em ergodic component} (indeed, a  $\mu$
{\em  ergodic component} with respect to $f$), if it does not admit any
smaller invariant subset with positive measure, that is, if
$V\subset U$ is invariant, $f^{-1}(V)=V$, then either $\mu(V)$ or
$\mu(U\setminus V)$ is zero. The measure $\mu$ is called {\em ergodic} if
$X$ is an ergodic component.
\end{Definition}

We stress that in the definition of ergodic measure and ergodic
components we are not assuming the invariance of the measure $\mu$ with respect to $f$.
Let us give some examples of non-invariant ergodic measures.
\begin{Example}\label{ExampleErg1}
Given any $p\in X\setminus\mbox{Fix}(f)$, the Dirac measure $\delta_p$ is ergodic and non-invariant ($\mbox{Fix}(f)$ is the set of fixed points of $f$). More in general, given a finite subset $\cu\subset\co_f^-(p)$ of the pre-orbit of a point $p\in X$, let
$\mu=\frac{1}{\#\cu}\sum_{q\in\cu}\delta_{q}$. If $f^{-1}(\cu)\ne\cu$ then $\mu$ is an ergodic probability but not invariant.
\end{Example}

\begin{Example}\label{ExampleErg2}
Given an ergodic (not necessarily invariant) measure $\mu$, let $Y\subset X$ be such that $\mu(Y) \mu(f^{-1}(Y)\setminus Y)>0$. Then $\mu|_Y$, the restriction of $\mu$ to $Y$, is non-invariant and ergodic.
\end{Example}

\begin{Example}\label{ExampleErg3}
By Martens~\cite{Ma}, the Lebesgue measure is ergodic and non-invariant for every non-flat $S$-unimodal map $f$ without a periodic attractor. In particular when $f$ is an infinitely renormalizable map the Lebesgue measure is ergodic but there is no absolutely continuous invariant measure
(for multimodal maps, see Blokh-Lyubich~\cite{BL89,BL91} and Vargas-van Strien~\cite{vSV}).
\end{Example}

Following Milnor's definition of attractor (indeed, minimal
attractor \cite{Mi}), a compact positively invariant set $A$
will be called  a $\mu$-attractor, or for short, an attractor, if
its basin of attraction $\mathcal{B}_f({A})=\{x\in
{X}\,{;}\,{{\omega}_f}(x)\subset A\}$ has positive measure and, in
contrast, the basin of every positively invariant compact subset
$A'\subsetneqq{A}$ has zero measure. Here,
${{\omega}_f}(x)$ denotes the $\omega$-limit set of $x\in {X}$.

A collection $\mathcal{P}$ of sets with
positive measure  is called a {\em partition mod $\mu$} of $U\subset {X}$ if this
collection covers $U$ almost everywhere
($\mu(U\setminus\bigcup_{P\in\mathcal{P}}P)=0$) and $\mu(P\cap Q)=0$ for every $P,Q\in\cp$ with $P\ne Q$.
The {\em diameter} of a partition $\cp$ is defined by $\diam (\cp)$ $=$ $\sup\{\diam(P)$ ${;}$ $P\in\cp\}.$

\begin{Proposition}[Ergodic attractors]\label{PropositionErgodicComponents2}
Given an ergodic component $U\subset {X}$, there exists a unique
attractor $A\subset {X}$ that attracts almost every point of $U$.
Moreover, ${{\omega}_f}(x)=A$ for almost every point of $U$.
\end{Proposition}

\dem Let $\mathcal{P}_1$ be any finite partition (mod $\mu$) of $U$ formed
by open subsets and with $\diam(\cp_1)<1$.
We will construct by induction a sequence of partitions $\cp_1<\cp_2<...$ of $U$
in the measure-theoretical sense. Thus, suppose that the collection
$\mathcal{P}_{n-1}$ has already been constructed. Set, for each
$P\in\mathcal{P}_{n-1}$, $U_P=\{x\in U\,{;}\,{{\omega}_f}(x)\cap
P\ne\emptyset\}$. As $U$ is an ergodic component and
$f^{-1}(U_P)=U_P$ (because ${{\omega}_f}(x)={\omega}(f(x))$ $\forall\,x$), either $U_P$ or $U\setminus U_P$ is a zero
measure set.

Given $P\in\cp_{n-1}$, we define a partition $\mathbb P_P$ (mod $\mu$) of $P$ as follows. If $\mu(U_P)=0$, we set $\mathbb P_P$ as the trivial refinement, i.e., $\mathbb P_P=\{P\}$. On the other hand, if $\mu(U_P)>0$, we choose any $\mathbb P_P$ in the collection of  finite partitions (mod $\mu$) of $P\in\cp_{n-1}$ formed by open subset of $P$ with diameter smaller that $\frac{1}{2}\diam (P)$ (see Lemma~\ref{LemmaSmallPartition}). Now, define
$$P_n=\{Q\in\mathbb P_P\,{;}\,P\in\cp_{n-1}\}.$$

For each $n\in\mathbb{N}$, set
$\mathcal{P}_n^*=\{P\in\mathcal{P}_n\,{;}\,U\setminus U_P$ is a
zero measure set$\}$ and $K_n=\bigcup_{P\in\mathcal{P}_n^*}\overline{P}$. As
$K_1\supset K_2\supset ...\supset K_n\supset...$ is a nested
sequence of non-empty compact sets, ${A}=\bigcap_n K_n$ is
also a non-empty compact set. By construction, for almost every
point $x\in U$ and $\forall n\in\mathbb{N}$, ${{\omega}_f}(x)\subset K_n$
and ${{\omega}_f}(x)\cap P$ $\forall P\in\mathcal{P}_n^*$. Moreover, as
$\diam (\overline{P})<2^{-n}$ $\forall P\in\mathcal{P}_n^*$,
it follows that $\sup\{\dist(y,\co_f(x))\,{;}\,y\in A\}\le 2^{-n}$
and ${{\omega}_f}(x)\subset K_n\subset
\overline{B_{2^{-n}}(A)}=\{p\in {X}\,{;}\,\dist(p,
A)\le 2^{-n}\}$ for every $n\in\NN$ and
almost every point $x\in U$. Thus, ${{\omega}_f}(x)=A$ for
$\mu$-almost every point $x\in U$.
\cqd

Consider for each point $x$ of a positively invariant set $U\subset {X}$,
a subset $\cu(x)\subset\co^+(x)$ of the positive orbit of $x$.
\begin{Definition}\label{DefinitionAsymInv}The collection $\cu=(\cu(x))_{x\in U}$ is called {\em asymptotically
invariant} if for every $x\in U$,
\begin{enumerate}
\item $\#\{j\in\NN\,{;}\,f^j(x)\in\cu(x)\}=\infty$ and
\item $\cu(x)\cap\co^+(f^{n}(x))=\cu(f(x))\cap\co^+(f^{n}(x))$ for every big $n\in\NN$.
\end{enumerate}
\end{Definition}

\begin{Definition}[${\omega}_{f,\cu}$]\label{DefinitionOmegaU}
Given an asymptotically invariant collection $\cu=(\cu(x))_{x\in U}$,
define for each $x$ the {\em omega-$\,\cu$ limit set} of $x$ ({\em omega-$\,\cu$} of $x$, for short),
denoted by ${\omega}_{f,\cu}(x)$, as the set of accumulation points of $\cu(x)$
and, that is, the set of points $p\in {X}$ such that there is a sequence $n_j\to+\infty$
satisfying $\cu(x)\ni f^{n_j}(x)\to p$.
\end{Definition}

It is easy to check that
${\omega}_\cu(x)$ is a non-empty compact set but not necessarily invariant.

We say that the asymptotically invariant collection $\cu=(\cu(x))_{x\in U}$
has {\em positive frequency} if $\limsup\frac{1}{n}\#\{1\le j\le n\,{;}\, f^j(x)\in\cu(x)\}>0$,
for every $x\in U$.

\begin{Definition}[${\omega}_{_{+},f,\cu}$]\label{DefinitionOmegaUTheta}
If $\cu$ is an asymptotically invariant collection with positive frequency,
define
${\omega}_{_{+},f,\cu}(x)$, the {\em set of $\cu$-frequently visited points} of $x$ orbit,
as the set of points $p\in {X}$ such that $\limsup\frac{1}{n}\#\{1\le j\le n\,{;}\,
f^j(x)\in\cu(x)\cap V\}>0$ for every neighborhood $V$ of $p$.
\end{Definition}

\begin{Lemma}\label{LemmaOmegaLimitSets}
Let $\cu=(\cu(x))_{x\in U}$ be an asymptotically invariant collection
defined in an ergodic component $U$ and let $A\subset {X}$ be the attractor
associated to $U$.
There is a compact set $A_{\cu}\subset A$
such that ${\omega}_{f,\cu}(x)=A_{\cu}$ for $\mu$-almost every
$x\in U$. Furthermore, if $\cu$ has positive frequency then
there is also a compact set $A_{_{+},\cu}\subset A_{\cu}$
such that ${\omega}_{_{+},f,\cu}(x)= A_{_{+},\cu}$ for $\mu$-almost every
$x\in U$.
\end{Lemma}

\dem We construct the compact sets ${A}_{\cu}$ and ${A}_{_{+},\cu}$
in the same way
we did for ${A}$ in the proof of
Proposition~\ref{PropositionErgodicComponents2}. For ${A}_{\cu}$
the only difference is that we have to change
${{\omega}_f}(x)$ by ${\omega}_{f,\cu}(x)$ in the proof.
Note that the key property of ${{\omega}_f}(x)$ used there is that ${{\omega}_f}(x)={\omega_f}(f(x))$
and we also have the same property for ${\omega}_{f,\cu}$,
i.e., ${\omega}_{f,\cu}(x)={\omega}_{f,\cu}(f(x))$.

For ${A}_{_{+},\cu}$ we have, of course, to change in the proof
${{\omega}_f}$ by ${\omega}_{_{+},f,\cu}$ (again ${\omega}_{_{+},f,\cu}(x)={\omega}_{_{+},f,\cu}(f(x))$ $\forall\,x$)
and we have also to change the definition
of the set $U_P$. For this we proceed as follows. Given a point $x\in
U$ and a set $K\subset {X}$ denote the $\cu$-visit frequency of $x$ to $K$
by $\phi_K(x)=\limsup\frac{1}{n}\#\{0\le j< n\,{;}\, f^j(x)\in K\cap\cu(x)\}$.
Set, for each $P\in\mathcal{P}_{n-1}$, $U_P=\{x\in
U\,{;}\,\phi_P(x)>0\}$. As we are using $\limsup$ in the definition
of $\phi_K$, we get $\phi_K(x)>0$ or $\phi_{{X}\setminus K}(x)>0$.
This is important to ensure that $K_n\ne\emptyset$ $\forall n$ (see
proof of Proposition~\ref{PropositionErgodicComponents2}).

To finish the proof, we remark that every point of
${A}_{_{+},\cu}=\bigcap_n K_n$ is accumulated by the sequence
$\{f^n(x)\,{;}\,n\in\NN$ and $f^n(x)\in\cu(x)\}$ for
almost every point $x\in U$ and so, ${A}_{_{+},\cu}$ is contained in
$A_\cu$ which is contained in ${A}$.
Moreover, if $B$ is an open set with
$B\cap{A}_{_{+},\cu}\ne\emptyset$ then for any big $n$ there will
be some element $P$ of $\mathcal{P}^*_n$ contained in $B$.
Therefore, by construction, $\limsup\frac{1}{n}\#\{0\le j< n\,{;}\,
f^j(x)\in B\cap\cu(x)\}\ge \limsup\frac{1}{n}\#\{0\le j< n\,{;}\, f^j(x)\in
P\cap\cu(x)\}>0$. \cqd

%%%%%%%%%%%%%%%%%%%%%%%%%%%%%%%%%%%%%%%%%%%%%%
\begin{figure}
  \includegraphics{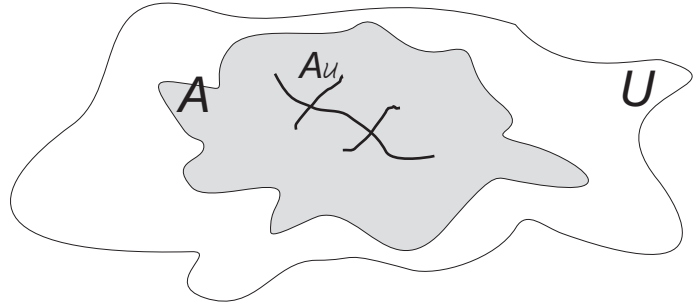}\\
  \caption{$U$ is an ergodic component with its attractor $A$ and its omega-$\,\cu$ set $A_\cu$.}
\end{figure}
%%%%%%%%%%%%%%%%%%%%%%%%%%%%%%%%%%%%%%%%%%%%%%

As defined in Section~\ref{SubsectionZoomingSetsAndMeasures}, a measure $\mu$ is $f$-non-singular if
the pre-image by $f$ of any set with zero measure has also
zero measure ($f_*\mu\ll\mu$). The ergodic measures that appears in Example~\ref{ExampleErg1}~and~\ref{ExampleErg2} are not in general $f$-non-singular. The lemma below gives a way to construct new $f$-non-singular ergodic measures.

\begin{Lemma}\label{LemmaContructingNon-SingErgMeasure}
If $\mu$ is a $f$-non-singular ergodic measure (not necessarily invariant) then $\frac{1}{\mu(E)}\mu|_{_E}$ is a $f$-non-singular ergodic probability whenever $E\subset X$ is a positively invariant Borel set with positive measure (i.e., $f(E)\subset E$ and $\mu(E)>0$).
\end{Lemma}
\dem
As in Example~\ref{ExampleErg2}, $\frac{1}{\mu(E)}\mu|_{_E}$ is an  ergodic probability.
We need only to show that this probability is $f$-non-singular. Given $Y\subset X$, we have
$\mu(f^{-1}(Y)\cap E)$ $\le$ $\mu(f^{-1}(Y)\cap f^{-1}(f(E)))$
$=$ $\mu(f^{-1}(Y\cap f(E)))$. Thus,
if $\mu|_{_E}(Y)=\mu(Y\cap E)=0$ then $0\le \mu|_{_E}(f^{-1}(Y))\le\mu(f^{-1}(Y\cap f(E)))=0$
(because $\mu$ is $f$-non-singular).
As a consequence,  $\frac{1}{\mu(E)}\mu|_{_E}$ is $f$-non-singular.

\cqd

\begin{Lemma}\label{LemmaCriteionForErgodicity}
Let $\mu$ be a finite measure.
If there exists some $\delta>0$ such that every
invariant set has $\mu$ measure either zero or bigger
than $\delta$, then ${X}$ can be decomposed into a finite number of $\mu$
ergodic components.
\end{Lemma}
\dem
Let $W_1\subset {X}$ be any invariant subset of ${X}$ (for example,
$W_1={X}$) with non zero $\mu$ measure and let $\mathcal{F}(W_1)$ be
the collection of all invariant subsets  $U\subset W_1$ with
$\mu$ measure bigger than zero. Note that  $\mathcal{F}(W_1)$ is
non-empty, because $W_1\in\mathcal{F}(W_1)$. Let us consider the
inclusion (mod $\mu$) as a partial  order on $\mathcal{F}(W_1)$.

\begin{claim}
Every  totally ordered subset $\Gamma\subset\mathcal{F}(W_1)$ is finite.
In particular, it has an upper bound.
\end{claim}
\dem
Otherwise there is an infinite sequence
$\gamma_0\supset \gamma_1\supset\gamma_3\supset\cdots$ with
$\mu(\gamma_{k}\setminus\gamma_{k+1})>0$ $\forall k$.
But as $\sum_k\mu(\gamma_{k}\setminus\gamma_{k+1})=\mu(\gamma_0)<\infty$,
it follows that $\mu(\gamma_{k}\setminus\gamma_{k+1})<\delta$ for $k$ big and this contradicts
our hypothesis as every $\gamma_{k}\setminus\gamma_{k+1}$ is an invariant set.
\cqd

From Zorn's lemma, there exists a maximal
element $U_1\in\mathcal{F}(W_1)$ and this is necessarily an ergodic
component.

As $W_2={X}\setminus U_1$ is
an invariant set, either it has zero $\mu$ measure or we can use
the argument above to $W_2$ and obtain a new ergodic component $U_2$ inside
${X}\setminus U_1$. Inductively, we can construct a collection of
ergodic components $U_1,...,U_r$ while $\mu({X}\setminus U_1\cup
...\cup U_r)>0$. But, as $\mu(U_j)>\delta$, this processes will stop
and we will get the decomposition of ${X}$ into $\mu$ ergodic
components as desired.

\cqd

\begin{Proposition}[A criterion for ergodicity] \label{PropositionFatErgodicAttractors}
Let $\mu$ be a $f$-non-singular finite measure.
If there exists some $\delta>0$ such that every
positively invariant set has $\mu$ measure either zero or bigger
than $\delta$, then ${X}$ can be decomposed into a finite number of $\mu$
ergodic components. Moreover, the attractor associated to each ergodic
component has positive $\mu$ measure.
\end{Proposition}
\dem
As every invariant set is positively invariant, it follows from Lemma~\ref{LemmaCriteionForErgodicity}
that ${X}$ can be decomposed into a finite number of $\mu$
ergodic components.

From Proposition~\ref{PropositionErgodicComponents2} each ergodic component $U$ of ${X}$ is
the basin of some attractor ${A}$.  Let us, for instance,
suppose that $\mu({A})=0$. In this case, one can choose an
open neighborhood $V$ of ${A}$ such that $\mu(V)<\delta$ and
an integer $n_0$ such that $\mu(U')>0$, where $U'=\{x\in U\ {;}\
f^n(x)\in V\ \forall n\ge n_0\}$. Note that $\mu(f^{n_0}(U'))>0$ because $\mu$ is $f$-non-singular.
As $U'$ is positively invariant, $f^{n_0}(U')$ is a positively invariant set with
$0<\mu(f^{n_0}(U'))<\mu(V)<\delta$, but this is impossible by ours hypothesis.
So, $\mu({A})>0$ (indeed, $\mu({A})>\delta$).
\cqd

We end this section relating the number of $\mu$ ergodic components with respect to $f$
to the number of  $\mu$ ergodic components with respect to $f^k$.

\begin{Lemma}\label{LemmaReletingErgNumber} Let $\mu$ be a $f$-non-singular finite measure.
If $U$ is an ergodic component with respect to $f$ then $U$ can be partitioned in at most $k$
ergodic components with respect to $f^k$. Furthermore, if $U_1,U_2\subset U$ are ergodic components
with respect to $f^k$ then $U_2=f^{-j}(U_1)$ (mod $\mu$) for some $0\le j<k$.
\end{Lemma}
\dem
First we will prove by induction that $U$ can be partitioned (mod $\mu$) in a finite number of ergodic components with respect to $f^k$.
Of course this claim is true for $k=1$. Thus, suppose by induction that for every $1\le j\le k-1$ we can decompose $U$ (mod $\mu$) in a finite number of ergodic components with respect to $f^j$.

If $U$ is ergodic with respect to $f^k$ there is nothing to prove. Thus, we may assume that there is an invariant set $Y\subset U$ (that is, $f^{-k}(Y)=Y$) with $0<\mu(Y)<\mu(U)$.

Let $\{j_1,...,j_s\}$ be a maximal subset of $\{1,...,k\}$ (with respect to the inclusion)
such that $\mu(Y\cap f^{-j_1}(Y)\cap...\cap f^{-j_s}(Y))>0$.
Set ${Y_1}=Y\cap f^{-j_1}(Y)\cap...\cap f^{-j_s}(Y)$. Note that $f^{-k}({Y_1})={Y_1}$.
Furthermore, by maximality,
if $\mu(f^{-\ell}({Y_1})\cap {Y_1})>0$ then $f^{-\ell}({Y_1})={Y_1}$ (mod $\mu$).
Let $a_1=\min\{1\le\ell\le k$ ${;}$ $f^{-\ell}({Y_1})={Y_1}\}$.
Of course, $f^{-1}(\bigcup_{j=0}^{a_1-1}f^{-j}({Y_1}))=\bigcup_{j=0}^{a_1-1}f^{-j}({Y_1})$ (mod $\mu$).
As $U$ is ergodic component for $f$, we get $$U=\bigcup_{j=0}^{a_1-1}f^{-j}({Y_1})\mbox{ (mod $\mu$)}.$$

\begin{claim}
${Y_1}$ is an ergodic component for $f^{a_1}$.
\end{claim}
\dem[Proof of the claim]
Suppose that ${Y_1}'\subset {Y_1}$ is $f^{a_1}$ invariant and $\mu({Y_1}\setminus {Y_1}')>0$.
As $f^{-a_1}({Y_1}\setminus {Y_1}')={Y_1}\setminus {Y_1}'$, we get
$f^{-1}(\bigcup_{j=0}^{a_1-1}f^{-j}({Y_1}\setminus {Y_1}'))=\bigcup_{j=0}^{a_1-1}f^{-j}({Y_1}\setminus {Y_1}')$ and,
as $\widetilde{U}$ is ergodic component for $f$, $\widetilde{U}=\bigcup_{j=0}^{a_1-1}f^{-j}({Y_1}\setminus {Y_1}')$ (mod $\mu$).
Thus
\begin{equation}\label{EquationReletingErgNumber2}
\sum_{j=0}^{a_1-1}\mu(f^{-j}({Y_1}))=\mu(\widetilde{U})=\sum_{j=0}^{a_1-1}\mu(f^{-j}({Y_1}\setminus {Y_1}')),
\end{equation}
because $\mu(f^{-i}({Y_1})\cap f^{-j}({Y_1}))=0$ $\forall\,0\le i<j\le a_1-1$ (here we are using that $\mu$ is $f$-non-singular).
As $\mu(f^{-j}({Y_1}))\ge\mu(f^{-j}({Y_1}\setminus {Y_1}'))$ $\forall\,j$, it follows
from (\ref{EquationReletingErgNumber2}) that $\mu(f^{-j}({Y_1}))=\mu(f^{-j}({Y_1}\setminus {Y_1}'))$ $\forall\,j$
and so, $\mu({Y_1}')=0$. \cqd

Denote by $\cu$ the collection of all ergodic component $\widetilde{U}\subset U$ with respect to some iterate $f^j$, $j=1,\cdots,k-1$.
By induction $\cu$ is finite and so, $\delta=\min\{\mu(\widetilde{U})$ $;$ $U\in\cu\}>0$.

From the claim above follows that if $U$ is not an ergodic component with respect to $f^k$ then every $f^k$-invariant set $Y\subset U$ with $0<\mu(Y)<\mu(U)$ contains some element of $\cu$. Thus, every positively invariant subset of $U$ has $\mu$ measure either zero or bigger than $\delta$. Applying Lemma~\ref{LemmaCriteionForErgodicity} to $\mu$ (indeed to $\widetilde{\mu}=\mu|_U$), it follows that $U$ can be decomposed into a finite number of $\mu$ ergodic components with respect to $f^k$.

To finish the proof of the lemma, let $W\subset U$ be an ergodic component with respect to $f^k$. As $f^{-k}(W)=W$, $f^{-1}(\bigcup_{j=0}^{k-1}f^{-j}(W))=\bigcup_{j=0}^{k-1}f^{-j}(W)$. Thus, by the ergodicity of $U$, $U=\bigcup_{j=0}^{k-1}f^{-j}(W)$ (mod $\mu$).
Note that, if $\widetilde{W}\subset U$ is an ergodic component with respect to $f^k$ and $\mu(\widetilde{W}\cap f^{-j}(W))>0$, then $\widetilde{W}= f^{-j}(W)$ (mod $\mu$), because $f^{-k}(\widetilde{W}\cap f^{-j}(W))=\widetilde{W}\cap f^{-j}(W)$ and $\widetilde{W}$ is ergodic with respect to $f^k$.
As $U=\bigcup_{j=0}^{k-1}f^{-j}(W)$ (mod $\mu$), we can conclude that any ergodic component $\widetilde{W}\subset U$ with respect to $f^k$ is (mod $\mu$) an element of $\{W,f^{-1}(W),...,f^{-(k-1)}(W)\}$.

\cqd

%%%%%%%%%%%%%%%%%%%%%%%%%%%%%%%%%%%%%%%%%%%%%%%%%%%%%%%%%%%%%%%%%%%%%%%%%%%%%%

\section{Characterizing the liftable measures}
\label{SectionIndMarMap}

In this section we obtain a statistical characterization of the liftable measures for a given
induced map (see Corollary~\ref{CorollaryRogers}). Differently of Zweimüller's results \cite{Zw},  this characterization is given by a statistical condition, condition (\ref{EquationStatisticalCondition}), not by the integrability of the induced time with respect to the reference measure (the one that we want to lift). This is important to avoid an additional condition of integrability of the induced time with respect to the reference measure (in our context this is not a natural condition).

Let $X$ be a compact separable metric space and $f:X\to X$ a measurable map defined on $X$.

\begin{Definition}[Markov map compatible with a measure]
\label{DefinitionMarkovMapCompatibleWithAMeasure}
We say that a Markov map $(F,\cp)$ defined on an open set $Y\subset X$ is {\em compatible} with a measure $\mu$ if
\begin{enumerate}
\item $\mu(Y)>0$;
\item $\mu$ is $F$-non-singular;
\item $\mu(\bigcup_{P\in\cp}P)=\mu(Y)$ (in particular, $\mu(\partial P)=0$ $\forall P\in\cp$).
\end{enumerate}
\end{Definition}

We say that a measure $\mu$ has a {\em Jacobian with respect to the map $f:X\to X$} if there is a function $J_{\mu}f\in L^1(\mu)$ such that $$\mu(f(A))=\int_A J_{\mu}f d\mu$$ for every measurable set $A$ such that $f|_A$ is injective. When the Jacobian exists, it is essentially unique. In general, the Jacobian may not exist, but if, for instance, $\mu$ is a $f$-invariant measure and $f$ is a countable to one map then the Jacobian of $\mu$ with respect to $f$ is well defined (see \cite{Par}).

\begin{Definition}[Markov map with $\mu$-bounded distortion]\label{DefinitionMU-boundDist}
We say that a Markov map $(F,\cp)$ defined on an open set $Y\subset X$ has {\em bounded distortion with respect to a measure $\mu$} (for short, has {\em $\mu$-bounded distortion}) if $(F,\cp)$ {\em compatible} with  $\mu$, $\mu$ has a Jacobian with respect to $F$ and
$\exists\,K>0$ such that
$$\left|\log\frac{J_\mu F(x)}{J_\mu F(y)}\right| \leq K \dist(F(x),
F(y)),$$ for $\mu$ almost every $\,x, y\in P$ and for all $\,P\in\cp$.
\end{Definition}

The remark below is a well known fact about projections of invariant measures of induced maps, see for instance Lemma~$3.1$ in Chapter~V of \cite{MvS}.

\begin{Remark}\label{RemarkProjection}
Let $(F,\cp)$ be an induced Markov map for $f$ defined on some $Y\subset X$ and let $R$ be its induced time. If $\nu$ is a $F$-invariant finite measure such that $\int R d\nu<\infty$ then
$$\eta =\sum_{P\in\cp}\sum_{j=0}^{R(P)-1}f_{\ast}^j\left(\nu|_{P}\right)\,
\bigg(=\sum_{j=0}^{+\infty}f^j_*(\nu|_{\{R>j\}})\bigg)$$
is a $f$-invariant finite measure.
\end{Remark}

Note that, if $(F,\cp)$ is compatible with a measure $\mu$,
the $\sigma$-algebra generated by $\{F^{-n}(P)$ ${;}$ $P\in\cp$ and $n\ge 0\}$ is
equal to the Borel sets of $U$ (mod $\mu$). Thus, using for example Lemma~$4.4.1$ of \cite{Aa}, it
is easy to obtain the following result.

\begin{Proposition}[Folklore Theorem]\label{FolkloreTheorem}
Let $\mu$ be $f$-non-singular measure. If $(F,\cp)$ is an induced full Markov map for $f$ with $\mu$-bounded distortion
then there exists an ergodic $F$ invariant probability $\nu\ll\mu$ whose density belongs to
$L^{\infty}(\mu)$. Indeed, $\log\frac{d\nu}{d\mu}\in L^{\infty}(\mu|_{\{\frac{d\nu}{d\mu}>0\}})$.

Moreover, if the inducing time  $R$ of $F$ is $\nu$-integrable, then
$\eta =\sum_{P\in\cp}\sum_{j=0}^{R(P)-1}f_{\ast}^j\left(\nu|_{P}\right)$
is a $\mu$ absolutely continuous ergodic $f$-invariant finite measure.
\end{Proposition}

In Theorem~\ref{TheoremRogers} we obtain an absolutely continuous $F$-invariant measure $\nu$ replacing the condition of bounded distortion (that appears in Proposition~\ref{FolkloreTheorem}) by $\mu$ being $f$-invariant and the statistical condition (\ref{EquationStatisticalCondition}). Furthermore, this statistical condition assures that projecting $\nu$ by the dynamics of $f$ we recover $\mu$. That is, every invariant measure satisfying (\ref{EquationStatisticalCondition}) can be lifted (indeed this is necessary and sufficient condition, see Corollary~\ref{CorollaryRogers}).

\begin{Theorem}\label{TheoremRogers}
Let $(F,\cp)$ be an induced full Markov map for $f$ defined on an open set $B\subset X$. Let $R$ be the inducing time of $F$ and $\mu$ be an ergodic $f$-invariant probability such that $\mu(\{R=0\})=0$.
If there exist $\Theta>0$ such that
\begin{equation}\label{EquationStatisticalCondition}
\limsup\frac{1}{n}\#\{0\le j<n\,;\,f^j(x)\in\co_F^+(x)\}\ge\Theta
\end{equation}
for $\mu$ almost every $x\in B$, where $\co_F^+(x)=\{F^j(x)\,;\,j\ge0\}$ is the positive orbit of $x$ by $F$, then there is a non trivial ($\not\equiv0$) finite $F$-invariant measure
$\nu$ such that $\nu(Y)\le\mu(Y)$ for all Borel set $Y\subset {B}$ and such that $\int Rd\nu\le \Theta^{-1}$.
\end{Theorem}
\dem Let $\mathfrak{B}=\{ x\in {B}$ $;$ $F^j(x)\in\bigcup_{P\in\cp}P$ $\forall\,j\ge0\}$.
Of course, $\mathfrak{B}$ is a metric space with the distance of $X$ and $\mathfrak{B}={B}$ ($\mu$ mod).

Let $\cw $ be a collection of subsets of $\mathfrak{B}$ formed by the empty set $\emptyset$ and all $Y\subset \mathfrak{B}$ such that $Y=(F|_{P_1})^{-1}\circ\dots\circ(F|_{P_s})^{-1}(\mathfrak{B})$ for some sequence of $P_1,...,P_s\in\cp$. That is, the elements of $\cw $ are the empty set and all homeomorphic $F$ pre-image of $\mathfrak{B}$. Note that $\cw $ is a collection of open sets of $\mathfrak{B}$. Given $Y\subset\mathfrak{B}$  and $r>0$, let $\cw (r,Y)$ be the set of all countable covers $\{I_i\}$ of $Y$ by elements of $\cw $ with $\diameter(I_i)\le r$ $\forall\,i$. It is clear that $\cw (r,Y)\ne\emptyset$ $\forall\,Y\subset\mathfrak{B}$ and $\forall\,r>0$.

Given a Borel set $Y\subset \mathfrak{B}$, let $\tau(Y)\in[0,1]$ be such that
$$\tau(Y)=\limsup_{n\to+\infty}\frac{1}{n}\{0\le j< n\,;\,f^j(x)\in Y\cap\co_F^+(x)\},$$
for $\mu$ almost every $x\in X$.

\begin{Claim}\label{ClaimRogersPropertiesTau}
The function $\tau$ has the following properties.
\begin{enumerate}
\item $\tau(\emptyset)=0;$
\item $\tau(\mathfrak{B})\ge\Theta>0$;
\item $\tau(Y_1)\le\tau(Y_2)$ whenever $Y_1\subset Y_2$ are Borel subsets of $\mathfrak{B}$;
\item $\tau(\bigcup_{i=1}^\infty Y_i)\le \bigcup_{i=1}^\infty \tau(Y_i)$ $\forall\, Y_1, Y_2, Y_3,
\dots$ Borel subsets of $\mathfrak{B}$;
\item $\tau(Y)\le\mu(Y)$ for all Borel set $Y\subset\mathfrak{B}$;
\item $\tau(F^{-1}(Y))=\tau(Y)$ for all Borel set $Y\subset\mathfrak{B}$.
\end{enumerate}
\end{Claim}
\dem[Proof of Claim~\ref{ClaimRogersPropertiesTau}]
The first four items follows from (\ref{EquationStatisticalCondition}) and the definition of $\tau$.
From Birkhoff Theorem follows the fifth item. Indeed, $\mu(Y)=\lim\frac{1}{n}\#\{0\le j< n$ $;$ $f^j(x)\in Y\}$ for every Borel set $Y\subset \mathfrak{B}$ and $\mu$ almost every $x$. Thus, $\tau(Y)\le\mu(Y)$ for every Borel set $Y\subset \mathfrak{B}$. To check the last item considers a Borel set $Y\subset\mathfrak{B}$.
As $F^j(x)\in Y$ $\Leftrightarrow$ $F^{j-1}(x)\in F^{-1}(Y)$ $\forall\,j\ge1$ and
$\forall\,x\in \mathfrak{B}$, we get $\tau(F^{-1}(Y))=\tau(Y)$.\cqd

Following the definition of pre-measure of Rogers~\cite{Ro}, $\tau$ restricted to $\cw $ is a pre-measure (Definition~5 of \cite{Ro}). Given $Y\subset\mathfrak{B}$, define
$$\nu(Y)=\sup_{r>0}\nu_r(Y)\,\,\bigg(=\lim_{r\searrow0}\nu_r(Y)\bigg),$$ where
$\nu_r(Y)=\inf_{_{\ci\in\cw(r,Y)}}\sum_{I\in\ci}\tau(I)$ and $\cw (r,Y)$ is the set of all countable covers $\ci=\{I_i\}$ of $Y$ by elements of $\cw $ with $\diameter(I_i)\le r$ $\forall\,i$.
The function $\nu$, defined on the class of all subset of $\mathfrak{B}$, is called in \cite{Ro} the {\em metric measure constructe from the pre-measured $\tau$ by Method~II} (Theorem~15 of \cite{Ro}).

As $(F,\cp)$ is a full Markov map,
\begin{equation}\label{EquationROGERS1}
\mbox{either $I_1\subset I_2$ or $I_2\subset I_1$ or $I_1\cap I_2=\emptyset$, $\forall\,I_1, I_2\in\cw $.}
\end{equation}
Thus,
$$\nu_r(Y)=\inf_{_{\ci\in\widetilde{\cw}(r,Y)}}\sum_{I\in\ci}\tau(I),$$
where $\widetilde{\cw }(r,Y)=\big\{\{I_i\}\in\cw (r,Y)$ $;$ $I_i\cap I_j=\emptyset$ $\forall\,i\ne j\big\}$.

\begin{Claim}\label{ClaimROGERS2}
$\nu(Y)\le\mu(Y)$ for every Borel set $Y\subset \mathfrak{B}$.
\end{Claim}
\dem[Proof of Claim~\ref{ClaimROGERS2}]
Let $Y\subset\mathfrak{B}$. As we are working only with countable additive measures defined on the Borel sets (see Section~\ref{SettingAndStatementOfMainResults}), $\mu$ is a regular measure. So, $\mu(Y)$ $=$ $\inf_{_{\ci\in\widetilde{\cw }(r,{Y})}}\mu\big(\bigcup_{I\in\ci}I\big)$  $=$ $\inf_{_{\ci\in\widetilde{\cw }(r,{Y})}}\sum_{I\in\ci}\mu(I)$ $\ge$ $\inf_{_{\ci\in\widetilde{\cw }(r,Y)}}\sum_{I\in\ci}\tau(I)$ $=$ $\nu_r(Y)$ for every $r>0$. Thus $\nu(Y)\le\mu(Y)$.
\cqd

It follows from Claim~\ref{ClaimRogersPropertiesTau}~and~\ref{ClaimROGERS2} that $\nu$ restricted to the Borel subsets $\mathfrak{B}$ is finite and non trivial, i.e., $\nu\not\equiv0$. Indeed, $\nu(\emptyset)$ $=$ $0$ $<$ $\theta$ $\le$ $\tau(\mathfrak{B})\le\nu(\mathfrak{B})$ $\le$ $\mu({B})$. Therefore, Theorem~19~and~3 of \cite{Ro} assures that $\nu$ restrict to the Borel subsets of $\mathfrak{B}$ is a countable additive measure.

Before we show that $\nu$ is $F$-invariant (Claim~\ref{ClaimROGERS3}) let us introduce some notation.
\begin{Notation}
Let $Y$ being Borel subset of $\mathfrak{B}$ and $r, r'>0$.
\begin{itemize}
\item Given $\ci\in\widetilde{\cw}(r,Y)$, set $F^*\ci=\{(F|_P)^{-1}(I)\}_{I\in\ci,P\in\cp}.$
\item Given $\ci\in\widetilde{\cw}(r,Y)$ and $x\in\mathfrak{B}$, let $\ci(x)$ be the element of $\ci$ that contains $x$, if $x\in\sum_{I\in\ci}I$. Otherwise, $\ci(x)=\emptyset$.
\item Given $\ci\in\widetilde{\cw}(r,Y)$ and $\ci'\in\widetilde{\cw}(r',Y)$, define $\ci\cap_Y \ci'=\{\ci(x)\cap \ci'(x)\,;\,x\in Y\}$.
\end{itemize}
\end{Notation}

Given $\ci\in\widetilde{\cw}(r,Y)$ and $\ci'\in\widetilde{\cw}(r',Y)$, follows easily from (\ref{EquationROGERS1}) that \begin{equation}\label{EquationROGERS3}
\ci\cap_Y \ci'\in\widetilde{\cw}(\min\{r,r'\},Y).
\end{equation}
Furthermore, (\ref{EquationROGERS1}) and Claim~\ref{ClaimRogersPropertiesTau} give
\begin{equation}\label{EquationROGERS4}
\sum_{I\in\ci\cap_Y\ci'}\tau(I)\le\min\{\sum_{I\in\ci}\tau(I),\sum_{I\in\ci'}\tau(I)\}.
\end{equation}

\begin{Claim}\label{ClaimROGERS3}
$\nu$ is $F$-invariant.
\end{Claim}
\dem[Proof of Claim~\ref{ClaimROGERS3}]
Let $Y$ be a Borel subset of $\mathfrak{B}$.
Let be a sequence of $a_1>a_2>\dots>a_j\searrow0$ a sequence of positive real numbers and $\ci_1,\ci_2,\dots$
a sequence of covers of $Y$ by elements of $\cw$ satisfying the follows properties.
\begin{enumerate}
\item[(P1) ] $\nu_{a_j}(Y)\le\nu_{a_{j+1}}(Y)$ $\forall j\ge1$;
\item[(P2) ] $\nu_{a_j}(Y)\le\nu(Y)<\nu_{a_j}(Y)+(1/j)$ $\forall j\ge1$;
\item[(P3) ] $\ci_j\in\widetilde{\cw}(a_j,Y)$ $\forall j\ge1$.
\end{enumerate}

Given $\varepsilon>0$, let $\delta>0$ be such that
\begin{equation}\label{EquationROGERS5}
\mu\bigg(\bigcup_{P\in\cp_1}P\bigg)<\frac{\varepsilon}{6},
\end{equation}
where $\cp_1=\{P\in\cp$ $;$ $\mu(P)<\delta\}$.

Set $\cp_0=\{P\in\cp$ $;$ $\mu(P)\ge\delta\}$. Of course, $n_0:=\#\cp_0<\infty$. For each $P\in\cp_0$, let $0<b_p<\diameter(P)$ be such that
\begin{equation}\label{EquationROGERS6}
\nu_{b_P}((F|_P)^{-1}(Y))\le\nu((F|_P)^{-1}(Y))\le\nu_{b_P}((F|_P)^{-1}(Y))+\frac{\varepsilon}{6 n_0},
\end{equation}
and let $\cj_P\in\cw(b_P,(F|_P)^{-1}(Y))$ such that
\begin{equation}\label{EquationROGERS7}
\nu_{b_p}((F|_P)^{-1}(Y))\le\sum_{J\in\cj_P}\tau(J)\le\nu_{b_P}((F|_P)^{-1}(Y))+\frac{\varepsilon}{6 n_0}.
\end{equation}

As $b_P<\diameter(P)$, it follows from (\ref{EquationROGERS1}) that $J\subset (F|_P)^{-1}(\mathfrak{B})=P\cap\mathfrak{B}$ $\forall J\in\cj_P$. Thus, for every $P\in\cp$ with $R(P)\le n_0$ we have
$$(F|_P)^{-1}(Y)\subset \bigcup_{J\in\cj_P}J \subset P$$

As $\bigcup_{J\in\cj_P}J \subset P$ and $F|_P$ is a homeomorphism, it follows that $\{F(J)\}_{_{J\in\cj_P}}$ $\in$ $\widetilde{\cw}(r_B,Y)$
$\forall\,P\in\cp_0$, where $r_B=\diameter(B)$. So, by (\ref{EquationROGERS3}),
$$\mbox{$\cj_0$ $:=$ ${{\large\mbox{$\bigcap_{_{Y}}$}}\atop {P\in\cp_0}}\{F(J)\}_{_{J\in\cj_P}}\in\widetilde{\cw}(r_B,Y)$}$$
$$\mbox{and}$$
$$\ci_j\cap_Y\cj_0\in\widetilde{\cw}(a_j,Y)\mbox{ for every }j\ge1.$$

Given $P\in\cp_0$, note that
\begin{equation}\label{EquationROGERS10}
\{(F|_P)^{-1}(I)\}_{_{I\in\ci_j\cap_{_{Y}}\cj_0}}=\ck_1 \cap_{Y_P} \cj_P\cap_{Y_P}\ck_2,
\end{equation}
where $Y_P=(F|_P)^{-1}(Y)$, $\ck_1=\{(F|_P)^{-1}(I)\}_{I\in\ci_j}$ and
$$\ck_2=\mbox{${{\large \mbox{$\bigcap_{_{Y_P}}$}} \atop P\ne Q\in\cp_0}\{(F|_P)^{-1}(I)\}_{I\in\cj_{Q}}$}.$$
It follows from (\ref{EquationROGERS3}) and (\ref{EquationROGERS10}) that
\begin{equation}\label{EquationROGERS11}
\{(F|_P)^{-1}(I)\}_{_{I\in\ci_j\cap_{_{Y}}\cj_0}}\in\widetilde{\cw}(b_P,(F|_P)^{-1}(Y)),\,\,\,\,\forall P\in\cp_0.
\end{equation}
Furthermore, by (\ref{EquationROGERS4}) and (\ref{EquationROGERS10}) we get
\begin{equation}\label{EquationROGERS11b}
\sum_{I\in\{(F|_P)^{-1}(I)\}_{_{I\in\ci_j\cap_{_{Y}}\cj_0}}}\tau(I)\le
\sum_{I\in\cj_P}\tau(I),\,\,\,\,\forall P\in\cp_0.
\end{equation}

Using the definition of $\nu_{b_p}$, (\ref{EquationROGERS11}), (\ref{EquationROGERS11b}) and (\ref{EquationROGERS7}), we obtain for all $P\in\cp_0$ that
$$\nu_{b_P}((F|_P)^{-1}(Y))\le\sum_{I\in\{(F|_P)^{-1}(I)\}_{_{I\in\ci_j\cap_{_{Y}}\cj_0}}}\tau(I)= $$
\begin{equation}\label{EquationROGERS12}
=\sum_{I\in\ci_j\cap_Y\cj_0}\tau((F|_P)^{-1}(I))\le\sum_{I\in\cj_P}\tau(I) \le\nu_{b_P}((F|_P)^{-1}(Y))+\frac{\varepsilon}{6 n_0}.
\end{equation}

Therefore

$$\bigg|\nu(F^{-1}(Y))-\sum_{I\in F^*(\ci_j\mbox{$\cap_Y$}\cj_0)}\tau(I)\bigg|=$$
$$=\bigg|\sum_{P\in\cp}\nu((F|_P)^{-1}(Y))-\sum_{P\in\cp}\sum_{I\in\ci_j\cap_Y\cj_0}\tau((F|_P)^{-1}(I))\bigg|\le$$
$$\le\sum_{P\in\cp_0}\bigg|\nu((F|_P)^{-1}(Y))-\sum_{I\in\ci_j\cap_Y\cj_0}\tau((F|_P)^{-1}(I))\bigg|+$$
$$+\sum_{P\in\cp_1}\nu((F|_P)^{-1}(Y))
+\sum_{P\in\cp_1}\sum_{I\in\ci_j\cap_Y\cj_0}\tau((F|_P)^{-1}(I))<$$
$$<\sum_{P\in\cp_0}\bigg|\nu((F|_P)^{-1}(Y))-\sum_{I\in\ci_j\cap_Y\cj_0}\tau((F|_P)^{-1}(I))\bigg|+$$
$$+\underbrace{\sum_{P\in\cp_1}\mu((F|_P)^{-1}(Y))}_{*}
+\underbrace{\sum_{P\in\cp_1}\sum_{I\in\ci_j\cap_Y\cj_0}\mu((F|_P)^{-1}(I))}_{**}.$$

As $*$ $\le$ $\mu(\bigcup_{P\in\cp_1}P)$ and also $**$ $\le$ $\sum_{P\in\cp_1}\mu((F|_P)^{-1}(\bigcup_{I\in\ci_j\cap_Y\cj_0}I))$ $\le$ $\sum_{P\in\cp_1}\mu(P)$ $=$ $\mu(\bigcup_{P\in\cp_1}P)$,
it follows from (\ref{EquationROGERS5}) that
$$\bigg|\nu(F^{-1}(Y))-\sum_{I\in F^*(\ci_j\mbox{$\cap_Y$}\cj_0)}\tau(I)\bigg|<$$
$$<\sum_{P\in\cp_0}\underbrace{\bigg|\nu((F|_P)^{-1}(Y))-\sum_{I\in\ci_j\cap_Y\cj_0}
\tau((F|_P)^{-1}(I))\bigg|}_{***}+\varepsilon/3.$$
By (\ref{EquationROGERS6}) and (\ref{EquationROGERS12}),
$$***\le \bigg|\nu((F|_P)^{-1}(Y))-\nu_{b_P}((F|_P)^{-1}(Y))\bigg|+$$
$$+\bigg|\nu_{b_P}((F|_P)^{-1}(Y))-\sum_{I\in\ci_j\cap_Y\cj_0}\tau((F|_P)^{-1}(I))\bigg|<\frac{\varepsilon}{3 n_0}.$$

Therefore,
\begin{equation}\label{EquationROGERS13}\bigg|\nu(F^{-1}(Y))-\sum_{I\in F^*(\ci_j\mbox{$\cap_Y$}\cj_0)}\tau(I)\bigg|<\frac{2}{3}\varepsilon.
\end{equation}

Let $j>3/\varepsilon$. Using the properties of $\tau$ (Claim~\ref{ClaimRogersPropertiesTau}), the fact that $\ci_j\cap_{_Y}\cj_0$ $\in$ $\widetilde{\cw}(a_j,Y)$, (P2) and (\ref{EquationROGERS13}), we get
$$\nu(Y)<\frac{1}{j}+\nu_{a_j}(Y)\le\frac{1}{j}+\sum_{I\in\ci_j\cap_{_Y}\cj_0}\tau(I)=$$
$$=\frac{1}{j}+\sum_{I\in\ci_j\cap_{_Y}\cj_0}\tau(F^{-1}(I))=
\frac{1}{j}+\sum_{I\in\ci_j\cap_{_Y}\cj_0}\tau\bigg(\sum_{P\in\cp}(F|_P)^{-1}(I)\bigg)\le$$
$$\le\frac{1}{j}+\sum_{I\in\ci_j\cap_{_Y}\cj_0}\sum_{P\in\cp}\tau\big((F|_P)^{-1}(I)\big)=
\frac{1}{j}+\sum_{I\in F^*(\ci_j\cap_{_Y}\cj_0)}\tau(I)\le$$
$$\le\frac{1}{j}+\nu(F^{-1}(Y))+\frac{2}{3}\varepsilon<\nu(F^{-1}(Y))+\varepsilon.$$

Thus, given a Borel set $Y\subset\mathfrak{B}$, we can conclude that $\nu(Y)<\nu(F^{-1}(Y))+\varepsilon$ for every
$\varepsilon>0$. That is, $$\mbox{$\nu(Y)\le\nu(F^{-1}(Y))$ for all Borel set $Y\subset\mathfrak{B}$.}$$

To conclude the proof of Claim~\ref{ClaimROGERS3}, let us assume the existence of a Borel set $L\subset\mathfrak{B}$ such that $\nu(L)<\nu(F^{-1}(L))$. As $\nu(\mathfrak{B}\setminus L)\le \nu(F^{-1}(\mathfrak{B}\setminus L))$, we obtain $\nu(\mathfrak{B})=\nu(L)+\nu(\mathfrak{B}\setminus L)<\nu(F^{-1}(L))+\nu(F^{-1}(\mathfrak{B}\setminus L))=\nu(\mathfrak{B})$, which is an absurd.\cqd

Now, suppose that $\int R d\nu\in(\gamma,+\infty]$, for some $\frac{1}{\Theta}<\gamma\in\RR$.
As $\nu$ is $F$ invariant and $R\ge0$, it follows from Birkhoff Theorem that $\exists\,\widetilde{\mathfrak{B}}\subset\mathfrak{B}$ with $\nu(\widetilde{\mathfrak{B}})>0$ such that for every $x\in\widetilde{\mathfrak{B}}$ there is some $n_x\in\NN$ satisfying $\sum_{k=0}^{n}R\circ F^k(x)>\gamma\,n$ $\forall\,n\ge n_x$. In this case, for every $n>\gamma\,n_x$ $(\ge n_x)$ and every $\frac{1}{\gamma}\,n\le j< n$ we get $\sum_{k=0}^{j}R\circ F^k(x)>\gamma\,j=\gamma\frac{j}{n}n\ge n$.
Thus, $$\sup\{j\ge0\,;\,\sum_{k=0}^{j}R\circ F^k(x)< n\}\le\frac{1}{\gamma}n<\Theta\,n,$$
for all $n\ge\gamma\, n_x$ and all $x\in\widetilde{\mathfrak{B}}$.

Because
$\{j\ge0$ $;$ $\sum_{k=0}^j R\circ F^k(x)<n\}$ $=$ $\{0\le j< n$ $;$ $f^j(x)\in\co_F^+(x)\}$
and $\sup\{j\ge0$ $;$ $\sum_{k=0}^{j}R\circ F^k(x)< n\}$ $=$ $\#\{j\ge0$ $;$ $\sum_{k=0}^{j}R\circ F^k(x)< n\}$, it follows that

\begin{equation}\label{EquationROGERS8.5}
\#\{0\le j< n\,;\,f^j(x)\in\co_F^+(x)\}=\sup\{j\ge0\,;\,\sum_{k=0}^j R\circ F^k(x)<n\}.
\end{equation}

So, for every $x\in\widetilde{\mathfrak{B}}$, we get
\begin{equation}\label{EquationROGERS9}
\limsup\frac{1}{n}\#\{0\le j< n\,;\,f^j(x)\in\co_F^+(x)\}<\Theta
\end{equation}
But this is
a contradiction. Indeed, as $\nu\ll\mu$, we have by hypothesis that $\nu(\{x\in\mathfrak{B}\,;\, (\ref{EquationROGERS9})\mbox{ holds}\})=\nu(\mathfrak{B}\setminus\{x\in\mathfrak{B}\,;\mbox{ (\ref{EquationStatisticalCondition}) holds}\})=0$. This proves that $\int R d\nu\le\Theta^{-1}$.
To finish the proof of the theorem, we extend $\nu$ to $B$ by setting $\nu(B\setminus\mathfrak{B})=0$.
\cqd

Using Theorem~\ref{TheoremRogers} we obtain the following characterization of the liftable
measures.

\begin{Corollary}\label{CorollaryRogers}Let $(F,\cp)$ be an induced full Markov map for $f$ defined on an open set $B\subset X$. Let $R$ be the inducing time of $F$ and $\mu$ be an ergodic $f$-invariant probability such that $\mu(\{R=0\})=0$.
The following statements are equivalent.
\begin{enumerate}
\item[(i)] There is a $F$-invariant finite measure $\nu\ll\mu$ such that $\mu=\sum_{j=0}^{+\infty}f^j_*(\nu|_{\{R>j\}})$.
\item[(ii)] For $\mu$ almost every $x\in B$, $\limsup_n\frac{1}{n}\#\{0\le j<n\,;\,f^j(x)\in\co_F^+(x)\}>0$.
\item[(iii)] For $\mu$ almost every $x\in B$, $\limsup_n\frac{1}{n}\sup_j\{j\ge0\,;\,\sum_{k=0}^j R\circ F^k(x)<n\}>0$.
\item[(iv)] There is a $F$-invariant finite measure $\nu\ll\mu$ such that $0<\int Rd\nu<\infty$.
\end{enumerate}
\end{Corollary}
\dem
By (\ref{EquationROGERS8.5}) follows that (ii)$\Leftrightarrow$(iii).
As $\mu=\sum_{j=0}^{+\infty}f^j_*(\nu|_{\{R>j\}})$ implies that $\int R d\nu$ $=$ $\sum_{j=0}^{+\infty}f^j_*(\nu|_{\{R>j\}})(X)$ $=$ $\mu(X)$, it follows that (i)$\Rightarrow$(iv). We get (i)$\Leftarrow$(iv) from Proposition~\ref{FolkloreTheorem}.
As (ii)$\Rightarrow$(i) follows from Theorem~\ref{TheoremRogers}, only (iv)$\Rightarrow$(iii)
remains to be proved.

Suppose that (iv) holds. For every $n\in\NN$ and each $x\in\mathfrak{B}:=\{ x\in {B}$ $;$ $F^j(x)\in\bigcup_{P\in\cp}P$ $\forall\,j\ge0\}$, let $i_x(n)=\sup_j\{j\ge0\,;\,\sum_{k=0}^j R\circ F^k(x)<n\}$. Thus, for every $x\in {B}$,
\begin{equation}\label{EquationCorollaryR}
\frac{1}{i_x(n)+1}\sum_{k=0}^{i_x(n)+1}R\circ F^k(x)\ge \frac{1}{i_x(n)+1}n=\frac{n}{i_x(n)}\bigg(\frac{i_x(n)}{i_x(n)+1}\bigg).
\end{equation}

If, by contradiction, $\limsup_n\frac{1}{n}\sup_j\{j\ge0\,;\,\sum_{k=0}^j R\circ F^k(x)<n\}=0$ for $\mu$ almost every $x\in {B}$, then $\lim_{n}n/i_x(n)=\infty$ for $\mu$ almost every $x\in {B}$. Using (\ref{EquationCorollaryR}) it follows that $\limsup_k\frac{1}{k}\sum_{j=0}^{k-1}R\circ F^j(x)=\infty$ for $\mu$ almost every $x\in {B}$. This contradicts (iv) as, by Birkhoff Theorem, $\limsup_k\frac{1}{k}\sum_{j=0}^{k-1}R\circ F^j(x)=\int R d\nu<\infty$ for $\nu$ for a.e. $x\in {B}$ and so, $\mu(\{x\in {B}$ $;$ $\limsup_k\frac{1}{k}\sum_{j=0}^{k-1}R\circ F^j(x)<\infty\})>0$.

\cqd

Lemma~\ref{LemmaIntegrabilidade} just below will be useful to bound the space average of the induced time by having some information about the time average of the induced time. This will be necessary for projecting an invariant measure of the induced map onto a $f$-invariant measure.

\begin{Lemma}
\label{LemmaIntegrabilidade}Let $\{G_j\}_{j\in\NN}$ be a collection of
subsets of $X$ such that $f^j(x)\in G_{n-j}$ $\forall\,0\le j<n$ $\forall\,x\in G_n$.
Let $B\subset X$ and let $x\in B$ be a point such that
$\#\{j\ge 0$ $;$ $x\in G_j$ and $f^j(x)\in B\}$ $=\infty$.
Let $T:\co^+(x)\cap B\to\co^+(x)\cap B$ be a map given by $T(y)=f^{g(y)}(y)$,
with  $1\le g(y)\le\min\{j\in\NN$ ${;}$ $y\in{G}_j$ and $f^j(y)\in B\}$. Then
$$\#\{1\le j\le n\,;\,x\in{G}_j\mbox{ and }f^j(x)\in B\}\le\#\{j\ge0\,;\,\sum_{k=0}^j g(T^k(x))\le n\}.$$
Moreover, If
$\limsup_n\frac{1}{n}\#\{1\le j\le n\,{;}\,x\in G_j\mbox{ and }f^j(x)\in B\}>\Theta>0$
then
$$\liminf_{n\to\infty}\frac{1}{n}\sum_{j=0}^{n-1}g\circ T^j(x)\le{\Theta}^{-1}.$$
\end{Lemma}
\dem
Given $n\in\NN$, set $\Gamma_n=\{1\le j\le n$ $|$ $x\in{G}_j$ and $f^j(x)\in B\}$
and $\Sigma_n=\{j\ge0$ ${;}$ $\sum_{k=0}^j g(T^k(x))\le n\}$.

As $\Gamma_0=\emptyset=\Sigma_0$,
we have $\#\Gamma_0\le \#\Sigma_0$. By induction, assume that $\#\Gamma_j\le \#\Sigma_j$ $\forall\,0\le j<n$.
To prove that $\#\Gamma_n\le \#\Sigma_n$ we may assume that $n\in\Gamma_n$, otherwise $\#\Gamma_{n-1}=\#\Gamma_n$ and so, $\#\Gamma_{n}=\#\Gamma_{n-1}\le\#\Sigma_{n-1}\le\#\Sigma_n$. Let $\ell=\max\{j$ ${;}$ $j\in\Sigma_{n-1}\}$ and $s=\sum_{k=0}^{\ell} g(T^k(x))$. As $s\le n-1$ and $x\in G_n$, we have $T^{\ell+1}(x)=f^{s}(x)\in G_{n-s}$. Moreover, we also know that $f^s(x)\in B$, $f^{n-s}(f^s(x))=f^n(x)\in B$ and so, $g(f^s(x))\le n-s$ and, as a consequence, $\sum_{k=0}^{\ell+1} g(T^k(x))$ $=$ $\sum_{k=0}^{\ell} g(T^k(x))$ $+$ $g(T^{\ell+1}(x))\le$ $s+(n-s)\le n$. Therefore, $\ell+1\in\Sigma_n\setminus\Sigma_{n-1}$ and so, $\#\Gamma_{n}$ $=$ $\#\Gamma_{n-1}+1$ $\le$ $\#\Sigma_{n-1}+1$ $\le$ $\#\Sigma_n$ (as $n\in\Gamma_n$, $\Gamma_n=\{n\}\cup\Gamma_{n-1}$), completing the induction.

Assume now that $\limsup_n\frac{1}{n}\#\{1\le j\le n\,{;}\,x\in G_j\mbox{ and }f^j(x)\in B\}>\Theta>0$. If $\liminf_n\frac{1}{n}\sum_{k=0}^{n-1}g\circ T^k(x)>{\Theta}^{-1}$, there is some $n_0$ such that $\sum_{k=0}^{n}g\circ T^k(x)>{\Theta}^{-1}\,n$ $\forall\,n\ge n_0$. In this case, if $n_0\le \Theta n\le j\le n$ then $\sum_{k=0}^{j}g\circ T^k(x)>{\Theta}^{-1}j={\Theta}^{-1}\frac{j}{n}n\ge n$. So, $\#\Sigma_n(x)\le \Theta\,n$ $\forall\,n\ge n_0$ and, as a consequence of $\#\Gamma_n\le \#\Sigma_n$ $\forall n$,
$$
\limsup_{n\to+\infty}\frac{1}{n}\#\{1 \le j\le n\,{;}\,x\in G_j\mbox{ and }f^j(x)\in B\}=\limsup_{n\to+\infty}\frac{1}{n}\#\Gamma_n\le\Theta,
$$
contradicting our hypotheses.
\cqd

%%%%%%%%%%%%%%%%%%%%%%%%%%%%%%%%%%%%%%%%%%%%%%%%%%%%%%%%%%%%%%%%%%%%%%%%%%%%%%

\section{Zooming sets and measures}\label{SubsectionZoomingSetsAndMeasures}

In this section we introduce the notion of {\em zooming times}.
This notion captures and weakens the geometric aspects of the {\em hyperbolic times} (Section~\ref{ApplicationsExpandingMeasures}),
allowing more flexibility in the applications and examples.

Let $f:X\to X$ be a measurable map defined on a connected, compact, separable metric space.

\begin{Definition}[Zooming contraction]\label{DefinitionZoomingContration} A sequence $\alpha=\{\alpha_n\}_{1\le n\in\NN}$ of
functions $\alpha_n:[0,+\infty)\to[0,+\infty)$ is called a {\em zooming contraction} if it satisfies the following conditions
\begin{itemize}
\item $\alpha_n(r)<r$ $\forall\,r>0$ and $\forall\,n\ge1$;
\item $\alpha_n\circ\alpha_m(r)\le\alpha_{n+m}(r)$ $\forall\,r>0$ and $\forall\,n,m\ge1$;
\item $\sup_{0\le r\le1}\big(\sum_{n=1}^{\infty}\alpha_n(r)\big)<\infty$.
\end{itemize}
\end{Definition}

For instance, an exponential contraction corresponds to a zooming contraction $\alpha_n(r)=\lambda^n r$ with $0<\lambda<1$.
We note that we can deal with polynomial contractions ($\alpha_n(r)=n^{-a}r$, $a>1$) and also with contractions that becomes in small scales as weak as we want ($\alpha_n(r):=(\frac{1}{1+n\sqrt{r}})^2 r$ defines a zooming contraction and $\lim_{r\to0}\frac{a_n(r)}{r}=1$, see Example~\ref{ShiftPolinomial}).

Let $\alpha=\{\alpha_n\}_{n}$ be a zooming contraction and $\delta>0$ be a positive constant.

\begin{Definition}[Zooming times]
We say that $n\ge 1$ is a {\em $(\alpha,\delta)$-zooming time}
for $p\in X$ (with respect to $f$) if there is a neighborhood $V_{n}(p)$ of $p$ satisfying
\begin{enumerate}
\item $f^{n}$ sends $\overline{V_n(p)}$  homeomorphically onto $\overline{B_{\delta}(f^{n}(p))}$;
\item $\dist(f^j(x),f^j(y))\le\alpha_{n-j}\big(\dist(f^{n}(x),f^{n}(y))\big)$, for every $x,y\in V_{n}(p)$
and every $0\le j<n$.
\end{enumerate}
\end{Definition}

The ball $B_{\delta}(f^{n}(p))$ is called a {\em zooming ball} and the set $V_{n}(p)$ is called a {\em zooming pre-ball}.
Denote by $\tz_n(\alpha,\delta,f)$ the set of points of $X$ for which $n$
is a $(\alpha,\delta)$-zooming time.

\begin{Definition}[Zooming sets]\label{DefinitionZoomingSets}
A positively invariant set $\Lambda\subset X$ is called a {\em zooming set} if (\ref{eqZZ}) holds for every $x\in\Lambda$.
\end{Definition}

\begin{Definition}[Zooming measures]\label{DefinitionZoomingMeasure}
A $f$-non-singular finite measure $\mu$ defined on the Borel set of $X$ is called a {\em weak zooming measure} if $\mu$ almost every point has infinitely many $(\alpha,\delta)$-zooming times. A weak zooming measure is called a {\em zooming measure} if
\begin{equation}\label{eqZZ}
\limsup\frac{1}{n}\#\{1\le j\le n\,{;}\, x\in {\tz_j(\alpha,\delta,f)}\}>0,
\end{equation}
for $\mu$ almost every $x\in X$.
\end{Definition}

\begin{Definition}[Bounded distortion]\label{DefinitionMeasureWithBoundedDistortion} We say that a weak zooming measure $\mu$ has bounded distortion if
$\exists\,\rho>0$ such that, $\forall\,n\in\NN$ and $\mu$ almost every $p\in{\tz_n(\alpha,\delta,f)}$,
the jacobian of $f^n$ with respect to $\mu$, $J_{\mu}f^n$,
is well defined on $V_{n}(p)$ and
$$\left|\log\frac{J_{\mu}f^{n}(x)}{J_{\mu}f^{n}(y)}
\right|\le\rho \dist(f^{n}(x),f^{n}(y)),$$
for $\mu$-almost every $x$ and $y\in V_{n}(p)$.
\end{Definition}

\begin{Remark}\label{RemarkCONNECTED} We use the connectivity (indeed, local connectivity is enough) only in the proof of Lemma \ref{LemmaZoomingNestedBall} (the local connectivity is necessary to apply Proposition~\ref{PropositionExisteNested2}, where $\mathcal{A}=\{B_r(x)\}$). This Lemma assures the existence of {\em nested sets} containing a given point $x\in X$. Thus, to obtain all the results of Sections~\ref{SubsectionZoomingSetsAndMeasures},~\ref{SectionConstructingLocalInducingMap}~and~\ref{SectionGlobalInducedMarkovMap} we can remove the additional hypotheses above if the existence of sets like $(B_r(x))^{\star}$  can be ensured in another way.
\end{Remark}

\begin{Lemma}
\label{LemmaConcatenar}The zooming times have the following properties.
\begin{enumerate}
\item If $p\in\tz_j(\alpha,\delta,f)$ then $f^{\ell}(p)\in\tz_{j-\ell}(\alpha,\delta,f)$ for all $0\le\ell<j$.
\item If $p\in\tz_j(\alpha,\delta,f)$ and $f^j(p)\in\tz_\ell(\alpha,\delta,f)$ then $p\in\tz_{j+\ell}(\alpha,\delta,f)$.
\item If $p\in\tz_{j \ell}(\{\alpha_n\}_n,\delta,f)$ then $p\in\tz_{j}(\{\alpha_{\ell\, n}\}_n,\delta,f^{\ell})$.
\end{enumerate}
\end{Lemma}
\dem Follows easily from the properties of zooming times.
\cqd

It follows from Lemma~\ref{LemmaConcatenar} that if
$x\in\tz_{k\,m+j}(\{\alpha_{n}\}_n,\delta,f)$, with $0\le j<k$,
then $f^j(x)\in\tz_{k\,m}(\{\alpha_{n}\}_n,\delta,f)\subset\tz_m(\{\alpha_{k n}\}_n,\delta,f^k)$. Thus,
$$\limsup_m\tz_m(\{\alpha_{n}\}_n,\delta,f)\subset
\bigcup_{j=0}^{k-1}f^{-j}\big(\limsup_m\tz_{k\,m}(\{\alpha_{n}\}_n,\delta,f)\big)\subset$$
\begin{equation}\label{EquationPeqZZ}\subset\bigcup_{j=0}^{k-1}f^{-j}\big(\limsup_m\tz_m(\{\alpha_{k n}\}_n,\delta,f^k)\big).
\end{equation}

Let $\cz$ be the set of all points of $X$ with positive frequency of $(\{\alpha_n\}_n,\delta)$-zooming times, that is, (\ref{eqZZ}) holds.

\begin{Notation}\label{NotationCollectionZoomingPreBalls}
Denote by  $\ce_\cz=(\ce_{\cz,_n})_n$ as the collection of all
$(\alpha,\delta)$-zooming pre-balls, where $\ce_{\cz,n}=\{V_n(x)$
${;}$ $x\in\tz_n(\alpha,\delta,f)\}$ is the collection of all $(\alpha,\delta)$-zooming pre-balls of order $n$.
\end{Notation}

One can check easily that the collection of all
$(\alpha,\delta)$-zooming pre-balls is a {\em dynamically closed family of pre-images} as defined in Section~\ref{SectionNestedSets}.

Given $x\in {X}$ and $0<r<\delta$ let $(B_r(x))^{\star}$ be the set defined by
(\ref{EquationStar}). If $x\in (B_r(x))^{\star}$, it follows from Proposition~\ref{PropositionExisteNested2} (taking $\ca=\{B_r(x)\}$) that the connected component of $(B_r(x))^{\star}$ which contains $x$ is an $\ce_\cz$-nested set.

\begin{Definition}[Zooming nested balls]\label{DefinitionZoomingNestedBalls}
If $x\in (B_r(x))^{\star}$, define {\em the $(\alpha,\delta)$-zooming nested ball (with respect to $f$) of radius $r$ and center on $x$}, denoted by $B_r^{\star}(x)$, as the connected component of $(B_r(x))^{\star}$ which contains $x$.
\end{Definition}

Note that, as we have contraction in any zooming time, $B_r(x)$
cannot be contained in any zooming pre-image (with order bigger than zero) of itself. So $\ca=\{B_r(x)\}$, in the definition above, is indeed a collection of open sets as desired on Section~\ref{SubsectionConstrNestedSets}.

\begin{Remark}\label{RemarkIncreasingOrder} As two distinct $\ce_\cz$-pre-images of the same set cannot intersect (Remark~\ref{RemarkPreImagensNaoLinkadas}), the order of the elements of a chain are strictly increasing. That is, if $(P_0,...,P_n)$ is a chain of $\ce_\cz$-pre-images of $B_r(x)$  then $0<\ord(P_0)<\ord(P_1)<\cdots<\ord(P_n)$. \end{Remark}

\begin{Definition}[Backward separated map] We say that $f$ is backward separated if
\begin{equation}\label{EquationHypothesis2}\inf\bigg\{\dist\bigg(x,\bigcup_{j=1}^nf^{-j}(x)\setminus\{x\}\bigg)\,{;}\,x\in X\bigg\}>0\,\,\forall n\ge1.
\end{equation}
\end{Definition}

For instance, every continuous map $f$ with bounded number of pre-images ($\sup\{\#f^{-1}(x)$ ${;}$ $x\in X\}<+\infty$) is backward separated.

\begin{Lemma}[Existence of zooming nested balls]\label{LemmaZoomingNestedBall} If for some $0<r<\delta/2$ we have $\sum_{n\ge 1}{\alpha}_n(r)<r/4$ then the zooming nested ball $B_r^{\star}(x)$ is well defined and $B_r^{\star}(x)\supset B_{r/2}(x)$, $\forall\,x\in {X}$.
Furthermore, if $f$ is backward separated and $\sup_{r>0}\big(\sum_{n\ge 1}{\alpha}_n(r)/r\big)<+\infty$ then there exists $0<r_0<\delta/2$ such that given $0<\gamma<1$ one can find
$0<r_\gamma<r_0$, depending only on  $\delta$, $\alpha$ and $\gamma$, such that
$B_r^{\star}(x)\supset  B_{\gamma r}(x)$ $\forall x\in {X}$ and $\forall 0<r\le r_{\gamma}$.
\end{Lemma}

\dem If $\sum_{n\ge 1}{\alpha}_n(r)<r/4$, $0<r<\delta/2$, as the order of the elements of a chain of $\ca=\{B_r(x)\}$, $0<r<\delta$, are strictly increasing (Remark~\ref{RemarkIncreasingOrder}), the diameter of any chain is smaller than $\sum_{n\ge 1}\alpha_n\big(\diameter(B_r(x))\big)<r/2$. Thus, using Corollary~\ref{CorollarySmallDiameter}, we get $B_r^{\star}(x)\supset B_{r/2}(x)$.

Let suppose now that $f$ is backward separated and $\sup_{r>0}\big(\sum_{n\ge 1}{\alpha}_n(r)/r\big)<+\infty$. Given $0<\gamma<1$, let $n_0\in\NN$ be such that $\sum_{n>n_0}\alpha_n(r)<(1-\gamma)r/2$.
Let $\varepsilon>0$ be such that $\inf_x\dist(x,\bigcup_{j=1}^{n_0} f^{-j}(x)\setminus\{x\})
>\varepsilon$ $\forall\,x\in {X}$, $r_{\gamma}=\frac{1}{3}\min\{\varepsilon,\delta\}$ and $0<r\le r_\gamma$.
Note that if $j<n_0$ then $B_r(x)\cap P=\emptyset$ $\forall\,P\in\ce_{\cz,j}$ (because
$P\cap\big(\bigcup_{j=1}^{n_0} f^{-j}(x)\big)\ne\emptyset$ and $\diam (P)<r< \varepsilon/2$).
Thus, every chain of $\ce_\cz$-pre-images of $B_r(x)$ begins with a pre-image of order bigger than $n_0$.
By Remark~\ref{RemarkIncreasingOrder}, the diameter of any chain is smaller than $\sum_{n>n_0}\alpha_n\big(\diameter(B_r(x))\big)<(1-\gamma)r$ and,
as a chain intersects the boundary of $B_r(x)$, we can conclude that a chain cannot intersect
$B_{\gamma r}(x)$. So, $(B_r(x))^{\star}$, and also $B_r^{\star}(x)$, contains $B_{\gamma r}(x)$.

\cqd

\begin{Notation}
Given any sequence of sets $\{U_n\}_n$, denote by $\limsup_n U_n$ the set of points that belong
to infinitely many elements of this sequence, i.e.,
$${\limsup}_n U_n=\bigcap_{n\ge1}\bigcup_{j\ge n}U_j.$$
\end{Notation}

Using the notation above, $f$-non-singular finite measure $\mu$ is weak zooming if
 $\mu({X}\setminus\limsup{\tz_m(\alpha,\delta,f)})=0$ (see Definition~\ref{DefinitionZoomingMeasure}).

If $x\in X$ has a zooming time, we can define the
{\em first zooming time} of $x$
as $\min\{n\,{;}\,x\in{\tz_n(\alpha,\delta,f)}\}$.
It is easy to show that,
if $\mu$ is a finite $f$-non-singular measure and the first zooming time is well defined for $\mu$-almost everywhere, then $\mu$ is a weak zooming measure.
That is, $$\mu\big({X}\setminus\bigcup_{j=1}^{\infty}\tz_j(\alpha,\delta,f)\big)=
0\,\,\Longrightarrow\,\,\mu({X}\setminus\limsup{\tz_m(\alpha,\delta,f)})=0.$$

\begin{Notation}[The zooming images set] Denote the collection of zooming images of $f$ by $\fz=(\fz(x))_{x\in\limsup\tz_m(\alpha,\delta,f)}$,
where $\fz(x)=\{f^m(x)\,{;}\,m\in\NN$ and $x\in{\tz_m(\alpha,\delta,f)}\}$ is the set of
zooming images of $x$ by $f$.
\end{Notation}
It is easy to see that if $x\in{\tz_m(\alpha,\delta,f)}$ then $f^{m-j}(x)\in\tz_{m-j}(\alpha,\delta,f)$,
$\forall\,0\le j<m$. Thus, using this information and Lemma~\ref{LemmaConcatenar},
one can prove that $\fz$ is an asymptotically invariant collection.
Indeed, if $x\in\tz$ and $m_0$ is the first zooming time for $x$
then
$\big\{f^m(x)\,{;}\,m\ge 2$ and $x\in{\tz_m(\alpha,\delta,f)}\big\}$ $=$ $\big\{f^{m}(f(x))\,
{;}\,m\ge\max\{m_0-1,1\}$ and $f(x)\in{\tz_{m}(\alpha,\delta,f)}\big\}$.

In the lemma below, let $\mu$ be a weak $(\alpha,\delta)$-zooming measure with bounded distortion (see Definition~\ref{DefinitionMeasureWithBoundedDistortion}), where $\alpha=\{\alpha_n\}_n$.

\begin{Lemma}\label{LemmaErgodicComponents}
Suppose that for some $0<r_0<\delta/2$ and $p\in X$ the $(\alpha,\delta)$-zooming nested open ball $B_{r_0}^{\star}(p)$ is well defined and contains $B_{r_0/2}(p)$. If $U\subset {X}$ is positively invariant, $\mu(U)>0$ and
$\mu(\{x\in U$ ${;}$ $B_{r_0/2}(p)\cap{\omega}_{\fz}(x)\ne\emptyset\})>0$ then $\mu(B_{r_0/2}(p)\cap U)=\mu(B_{r_0/2}(p))>0$.
\end{Lemma}

\dem

Let $\rho>0$ be the distortion constant that appear in Definition~\ref{DefinitionMeasureWithBoundedDistortion}
and let $K\subset \{x\in U$ ${;}$ $B_{r_0}(p)\cap{\omega}_{\fz}(x)\ne\emptyset\}$ be a compact set
with positive $\mu$ measure.

Given $\ell>0$  choose an open neighborhood $V\supset K$ of
$K$ such that $\mu(V\setminus K)<\mu(K)/\ell$. Choose for each $x\in K$ a
zooming time $n(x)$ such that
$V_{n(x)}(x)\subset V$ and $f^{n(x)}(x)\in B_{r_0/2}(p)$. As $V_{n(x)}(x)$ is
mapped diffeomorphically by $f^{n(x)}$ onto
$B_{\delta}(f^{n(x)}(x))$ and $B_{\delta}(f^{n(x)}(x))\supset B_{r_0}^{\star}(p)$
(because $r_0<\delta/2$), set, for each $x\in K$,
$$W(x)=(f^{n(x)}|_{V_{n(x)}(x)})^{-1}\big(B_{r_0}^{\star}(p)\big).$$
By compactness $K\subset W(x_1)\cup...\cup W(x_s)$ for some
$x_1,...,x_m\in K$. As $B_{r_0}^{\star}(p)$ is a nested set, we can assume that
$W(x_j)\cap W(x_i)=\emptyset$ whenever $j\ne i$.
Thus, at least for one
$j$ we have $\mu(W(x_j)\setminus K)<\mu(W(x_j))/\ell$.
Otherwise,
$\mu(V\setminus K)$ $\ge$ $\mu\Big(\big(\bigcup_j W(x_j)\big)\setminus K\Big)$
$=$ $\sum_j \mu(W(x_j)\setminus K)$
$\ge$ $\sum_j\mu(W(x_j))/\ell$ $=$ $\mu(\bigcup_j W(x_j))/\ell$ $\ge\mu(K)/\ell$.
Therefore, for each $\ell\in\NN$ we can find some pre-ball $W_\ell$ which is sent by some iterate $f^{n_{\ell}}$
of $f$ diffeomorphically onto $B_{r_0}^{\star}(p)$ and with the distortion bounded by  $\rho$.
Furthermore, $\mu(W_\ell\setminus K)<\mu(W_\ell)/\ell$ $\forall\,\ell$. By the bounded
distortion we get $$\frac{\mu(B_{r_0}^{\star}(p)\setminus U)}{\mu(B_{r_0}^{\star}(p))}\le
\frac{\mu(B_{r_0}^{\star}(p)\setminus f^{n_\ell}(K))}{\mu(B_{r_0}^{\star}(p))}\le\rho\frac{\mu(W_\ell\setminus K)}{\mu(W_\ell)}
<\frac{\rho}{\ell}\to 0.$$
As $B_{r_0}^{\star}(p)\supset B_{{r_0}/2}(p)$, we conclude the proof.

\cqd

\begin{Corollary}\label{CorollaryBigSets} If $\mu$ is a weak zooming
measure with compact support and bounded distortion,
then there is $\varepsilon>0$ such that every positively invariant set has either
$\mu$-measure bigger than $\varepsilon$ or equal to zero.
\end{Corollary}
\dem
Let $\mu$ be a weak $(\alpha,\delta)$-zooming
measure with compact support and bounded distortion, where $\alpha=\{\alpha_n\}_n$.
Let $U$ be a positively invariant set with $\mu(U)>0$.

First, assume that $\sum_{n\ge1}\alpha_n({r_0})<{r_0}/4$ for some $0<{r_0}<\delta/2$.
It follows from  Lemma~\ref{LemmaZoomingNestedBall} that for every $p\in X$ the
$(\alpha,\delta)$-zooming nested open ball $B_{r_0}^{\star}(r)$ is well defined and contains $B_{{r_0}/2}(p)$.
Let $p$ be any point on the support of $\mu$ such that
$\mu(\{x\in U\,{;}\, {\omega}_{f,\fz}(x)\cap B_{r_0/2}(p)\ne\emptyset\})>0$ (of course at least one of such point exist).
It follows from Lemma~\ref{LemmaErgodicComponents} that
$\mu(U)\ge\mu(B_{{r_0}/2}(p)\cap U)=\mu(B_{{r_0}/2}(p))>0$. Let $\varepsilon:=
\inf\{\mu(B_{{r_0}/2}(x))$ ${;}$ $x\in\supp\mu\}$. It is easy to see that $\varepsilon>0$,
$\varepsilon$ does not depend on $U$ and $\mu(U)\ge\varepsilon>0$.

In the general case, let ${r_0}=\frac{\delta}{3}$ and $\widetilde{f}=f^{k}$,
where $k\ge 1$ is such that $\sum_{n\ge 1}\alpha_{k\,n}({r_0})$ $<$ $\frac{{r_0}}{4}$.

By (\ref{EquationPeqZZ}), there is $0\le j<k$ such that
$\mu\big(f^{-j}\big(\limsup_m\tz_m(\{\alpha_{k n}\}_n,\delta,f^k)\big)$ $\cap$ $U\big)$ $>$ $0$
and, as $\mu\circ f^{-1}\ll\mu$ and $f(U)\subset U$, we get
\begin{equation}\label{4545}
\mu\big(\limsup_m\tz_m(\{\alpha_{k n}\}_n,\delta,f^k)\cap U\big)>0.
\end{equation}

Taking $\widetilde{\mu}=\mu|_{\limsup_m\tz_m(\{\alpha_{k n}\}_n,\delta,f^k)}$,
it is easy to see that
$\widetilde{\mu}$ is a weak $(\{\widetilde{\alpha}_{n}\}_n,\delta)$-zooming measure
with respect to $\widetilde{f}$, where $\widetilde{\alpha}_n=\alpha_{k n}$.
Moreover $\widetilde{\mu}$ has compact support,
bounded distortion and $\sum_{n\ge 1}\widetilde{\alpha}_{n}({r_0})/{r_0}$ $<$ $1/4$.
As $\widetilde{f}(U)\subset U$ and, by (\ref{4545}), $\widetilde{\mu}(U)>0$,
we can apply the particular case and get $\varepsilon>0$,
not depending on $U$, such that $\mu(U)\ge\widetilde{\mu}(U)>\varepsilon$.
\cqd

As a consequence of Proposition~\ref{PropositionErgodicComponents2},
Proposition~\ref{PropositionFatErgodicAttractors}, Corollary~\ref{CorollaryBigSets}
and Lemma~\ref{LemmaOmegaLimitSets} we have
the following result.

\begin{Theorem}
\label{PropositionFatZoomingErgodicAttractors}
If $\mu$ is a weak zooming measure with bounded distortion
then ${X}$ can be partitioned into a finite collection
of $\mu$-ergodic components. Inside each $\mu$-ergodic component $U$ there
exists a fat attractor ${A}$ (i.e, $\mu(A)>0$)
such that ${{\omega}_f}(x)={A}$ for $\mu$-almost every
point $x\in U$.

Furthermore, there is a compact set ${A}_\fz\subset{A}$ such that
${\omega}_{f,\fz}(x)={A}_\fz$ for $\mu$-almost every point $x\in U$ and,
if $\mu$ is a zooming measure, there is a compact set
${A}_{_{+},\fz}\subset{A}_\fz$ such that ${\omega}_{+,f,\fz}(x)={A}_{_{+},\fz}$ for
$\mu$-almost every point $x\in U$.
\end{Theorem}

%%%%%%%%%%%%%%%%%%%%%%%%%%%%%%%%%%%%%%%%%%%%%%%%%%%%%%%%%%%%%%%%%%%%%%%%%%%%%%

\section{Constructing a local inducing Markov map}
\label{SectionConstructingLocalInducingMap}

Section~\ref{SectionConstructingLocalInducingMap} and \ref{SectionGlobalInducedMarkovMap} are the kernel of this paper.
Most of the results for zooming sets and measures are proved in these sections, and from them we will obtain their analogues for expanding sets and measures. The existence of an invariant measure $\nu{\ll}\mu$ that is absolutely continuous with respect to a given zooming measure with some distortion control is given by Theorem~\ref{TheoremZoomingInvariant}. In Theorem~\ref{TeoremMarkovStructure} we prove the existence of Markov structures  for zooming sets. The existence of global induced Markov maps for zooming sets is given in Section~\ref{SectionGlobalInducedMarkovMap} by  Theorem~\ref{GlobalMarkovStructureForZooming}.
Note that our approach to construct induced Markov map for dynamics with some hyperbolic behavior has to be very different from the one of Alves, Luzzatto, Pinheiro \cite{ALP,ALP1}, Gouëzel~\cite{G}, Pinheiro~\cite{Pi} and Young~\cite{Y1}. That is because
this construction in those papers depends in an essential way on the good
relation between the diameter and the volume (Lebesgue measure) of balls and this is not true for general zooming (or expanding) measures.

Let $X$, $f$, $\delta$, $\alpha=\{\alpha_n\}_n$ and $\fz=(\fz(x))_{x\in\limsup\tz_n(\alpha,\delta,f)}$ be as in Section~\ref{SubsectionZoomingSetsAndMeasures}. Let $$\Lambda\subset\limsup_{n\to\infty}\tz_n(\alpha,\delta,f)\subset X$$ be a positively invariant set.

Let $\Delta$ be a {\em $(\alpha,\delta)$-zooming nested open ball}, that is, $\Delta$ is a topological open ball and also a $(\alpha,\delta)$-zooming nested set. Assume also that $\diameter(\Delta)<\delta/2$. For example, if $\sum_{n\ge1}\alpha_n(r)<r/4$ for some $0<r<\delta/4$ (or if $f$ is backward separated and $\sup_{r>0}\sum_{n\ge1}\alpha_n(r)/r<+\infty$) we can take $\Delta$ as any zooming nested ball $B^{\star}_r(q)$ given by
Lemma~\ref{LemmaZoomingNestedBall}.

It is sometimes useful not to use all the zooming times but a sub-collection of them in the construction of the induced Markov map (for instance, this is necessary in the proof of item~(4) of Theorem~\ref{GlobalExpandingInducedMarkovMap}).
This motivates the definitions below.

For each $x\in\Lambda$ consider a set $\widetilde{\fz}(x)\subset\fz(x)$.
We say that $n$ is a $\widetilde{\fz}$-time for $x$ if $f^n(x)\in\widetilde{\fz}(x)$. A zooming pre-ball $V_n(x)$ is called
a $\widetilde{\fz}$-pre-ball if $n$ is a $\widetilde{\fz}$-time
for $x$. Let $\widetilde{\ce}_{\cz}\subset\ce_{\cz}$
be the collection of all $\widetilde{\fz}$-pre-balls $V_n(x)$ for all $x\in\Lambda$ and all $\widetilde{\fz}$-time for $x$.

\begin{Definition}\label{DefinitionPropreZooColl}We say that $\widetilde{\fz}=\big(\widetilde{\fz}(x)\big)_{x\in\Lambda}$ is a {\em proper zooming sub-collection} if
\begin{enumerate}
\item $\widetilde{\fz}$ is asymptotically invariant;
\item $\widetilde{\fz}(x)\subset\fz(x)$ for all $x\in\Lambda$;
\item $\widetilde{\fz}$ has positive frequency whenever $\fz$ has positive frequency;
\item $\widetilde{\ce}_{\cz}$ is a dynamically closed family of pre-images.
\end{enumerate}
\end{Definition}

The zooming collection itself is an example of a proper zooming sub-collection. Another example of proper zooming sub-collection that we are interested in is the following. Fixed $\ell\ge1$, set $\widetilde{f}=f^{\ell}$ and $\widetilde{\alpha}=\{\widetilde{\alpha}_n\}_n$, where $\widetilde{\alpha}_n=\alpha_{\ell\,n}$.
For instance, denote the collection of $(\alpha,\delta)$-zooming images of $f$ by
$\fz_{_f}=(\fz_{_f}(x))_{x\in\limsup_{n}\tz_n(\alpha,\delta,f)}$ and the collection of
$(\widetilde{\alpha},\delta)$-zooming images of $\widetilde{f}$ by
$\fz_{_{\widetilde{f}}}$ $=$ $(\fz_{_{\widetilde{f}}}(x))_{x\in\limsup_{n}
\tz_n(\widetilde{\alpha}\,\delta,\widetilde{f})}$. It follows from Lemma~\ref{LemmaConcatenar} that $\limsup_{n}
\tz_{\ell\,n}(\alpha\,\delta,f)$
$\subset$ $\limsup_{n}
\tz_n(\widetilde{\alpha},\delta,\widetilde{f})$. Thus, taking  $\widetilde{\fz}_{_{\widetilde{f}}}(x)$ $=$ $\{\widetilde{f}^n(x)$ $;$ $f^{\ell\,n}(x)\in\fz_{_f}(x)\}$, the collection $\widetilde{\fz}$ $=$ $\widetilde{\fz}_{\widetilde{f}}$ $=$ $(\widetilde{\fz}_{_{\widetilde{f}}}(x))_{x\in\limsup_{n}\tz_{\ell\,n}(\alpha,\delta,f)}$ is a proper $(\widetilde{\alpha},\delta)$-zooming sub-collection for the map $\widetilde{f}$.
This sub-collection will be necessary in the proof of item~(\ref{ItemGMSFZ3}.2) of Teorem~\ref{GlobalMarkovStructureForZooming} (and, as a consequence, item~(\ref{ItemE3}) of Theorem~\ref{GlobalExpandingInducedMarkovMap}) to acquire more contraction on the pre-balls (changing $\alpha$ for $\widetilde{\alpha}$ and $f$ for $\widetilde{f}$) maintaining the distortion control even for each iterate of the original map.
We emphasize  that for all the other results  we do not really need to work with a sub-collection.

Now, let $\widetilde{\fz}=\big(\widetilde{\fz}(x)\big)_{x\in\Lambda}$ be a proper zooming sub-collection  and let $\widetilde{\ce}_{\cz}\subset\ce_{\cz}$ be the collection of all $\widetilde{\fz}$-pre-balls.
Given $x\in\Delta$, let ${\Omega}(x)$ be the collection of
$\widetilde{\ce}_{\cz}$-pre-images $V$ of $\Delta$ such that $x\in V$.

The set $\Omega(x)$ is not empty for every $x\in\Delta$ that has a $\widetilde{\fz}$-return to $\Delta$. Indeed, if  $x\in\Delta$ and $f^n(x)$ $\in$ $\Delta\cap\widetilde{\fz}(x)$ then $B_{\delta}(f^n(x))$ $=$ $f^n(V_n(x))$ $\supset$ $\Delta$ (because $\diameter(\Delta)$ $<$ $\delta/2$). Thus, for each $\widetilde{\fz}$-return time of a point $x\in\Delta$ we can associated the $\widetilde{\ce}_{\cz}$-pre-image $P$ $=$ $(f^n|_{V_n(x)})^{-1}(\Delta)$ of $\Delta$ with $x\in P$.

\begin{Definition}\label{DefinitionRetIndTime}The  {\em inducing time} on
$\Delta$ associated to {\em ``the first $\widetilde{\ce}_{\cz}$-return time to $\Delta$}'' is the function $R:\Delta\to\NN$ given by
\begin{equation}\label{EquationInducing}
R(x)=\begin{cases}
\min\{\ord(V)\,;\,V\in{\Omega}(x)\} & \text{ if }{\Omega}(x)\ne\emptyset\\
0 & \text{ if }{\Omega}(x)=\emptyset
\end{cases}.
\end{equation}
\end{Definition}

Note that $R(x)$ is smaller than or equal to the first $\widetilde{\fz}$-return time to
$\Delta$, i.e., $R(x)\le\min\{n\ge1$ ${;}$ $f^n(x)\in\widetilde{\fz}(x)\cap\Delta\}$.
\begin{Definition}\label{DefinitionIndMapToRetIndTime}
The induced map $F$ on
$\Delta$ associated to {\em ``the first $\widetilde{\ce}_{\cz}$-return time to $\Delta$}'' is the map $F:\Delta\to\Delta$ given by
\begin{equation}\label{EquationinducedMap}
F(x)=f^{R(x)}(x),\,\forall x\in \Delta.
\end{equation}
\end{Definition}

As the collection of sets ${\Omega}(x)$ is totally ordered by inclusion, it follows from Corollary~\ref{CorollaryMainNestedProperty} that there is a unique $I(x)\in{\Omega}(x)$ such that $\ord(I(x))=R(x)$, whenever ${\Omega}(x)\ne\emptyset$.

%%%%%%%%%%%%%%%%%%%%%%%%%%%%%%%%%%%%%%%%%%%%%%
\begin{figure}
  \includegraphics{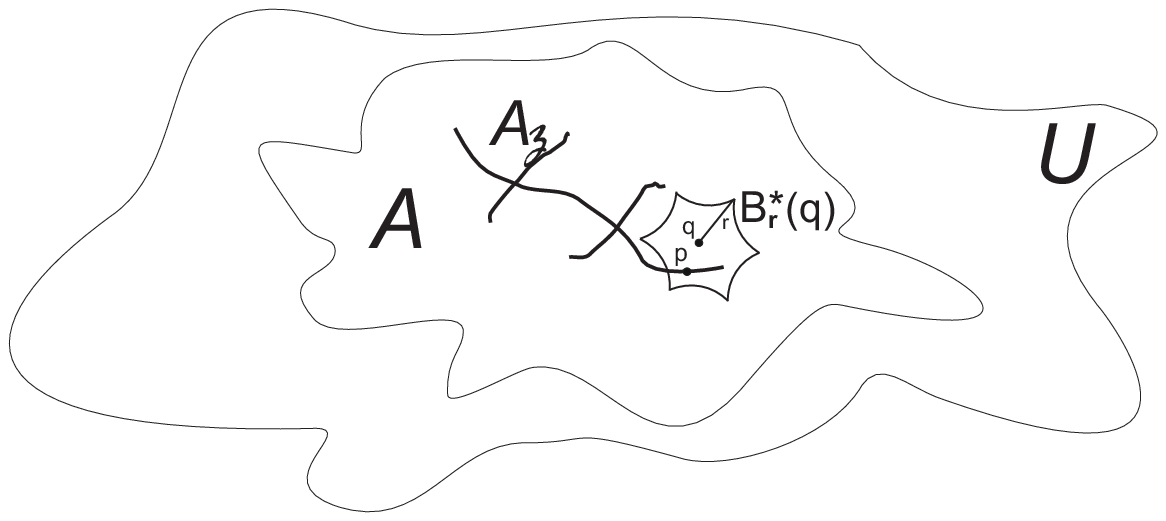}\\
  \caption{$\Delta=B_r^\star(x)$}
\end{figure}
%%%%%%%%%%%%%%%%%%%%%%%%%%%%%%%%%%%%%%%%%%%%%%

\begin{Lemma}\label{LemmaDisjointness}
If ${\Omega}(x)\ne\emptyset\ne{\Omega}(y)$ then either $I(x)\cap I(y)=\emptyset$ or $I(x)=I(y)$.
\end{Lemma}
\dem
We claim that, if ${\Omega}(x)\ne\emptyset$, $I(x)\supset V$ $\forall V\in{\Omega}(x)$. Indeed, if $I(x)\subsetneqq V$ with $V\in{\Omega}(x)$, as $\ord(I(x))<\ord(V)$, it follows from Corollary~\ref{CorollaryMainNestedProperty} that $\Delta$ is contained in an $\widetilde{\ce}_\cz$-pre-image of itself of order bigger than zero.  But this is impossible because we have contraction in the zooming times, i.e., the diameter of an $\widetilde{\ce}_\cz$-pre-image of $\Delta$ has diameter smaller than the diameter of $\Delta$.

Let $x,y\in X$ with ${\Omega}(x)\ne\emptyset\ne{\Omega}(y)$. As $I(x)$ and $I(y)$ are $\widetilde{\ce}_\cz$-pre-images of $\Delta$, if $I(x)\cap I(y)\ne\emptyset$ then $I(x)\supset I(y)$ or $I(x)\subset I(y)$. Thus, $I(x)\cap I(y)\ne\emptyset$ implies that $I(x)\in{\Omega}(y)$ or $I(y)\in{\Omega}(x)$. In any case, by uniqueness, $I(x)=I(y)$.
\cqd

\begin{Definition}\label{DefinitionMarkCollecFirstRetTime}
The Markov partition associated to {\em ``the first $\widetilde{\ce}_{\cz}$-return time to $\Delta$}'' is the collection of open sets $\cp$  given by
\begin{equation}\label{EquationIndMArkovCollection}
\cp=\{I(x)\,;\,x\in \Delta\mbox{ and }\,{\Omega}(x)\ne\emptyset\}.
\end{equation}
\end{Definition}

The Corollary below shows that $\cp$ is indeed a Markov partition of open sets.

\begin{Corollary}[Existence of a full induced Markov map for a zooming set]\label{CorollaryDisjointness} Let $F$ be given by (\ref{EquationinducedMap}), $R$ given by (\ref{EquationInducing}) and $\cp$ by (\ref{EquationIndMArkovCollection}). If $\cp\ne\emptyset$ then $(F,\cp)$ is an induced full Markov map for $f$ on $\Delta$.
\end{Corollary}
\dem
By construction the elements of $\cp$ are open sets. By Lemma~\ref{LemmaDisjointness}, $\cp$ satisfies the first condition of a Markov partition for $F$.
As $F(P)=\Delta\supset Q$ $\forall\,P,Q\in\cp$, $\cp$ also satisfies the second and third conditions of a Markov partition. On the other hand, as $F|_P=f^{\ord(P)}|_P$ and $P$ is a $\widetilde{\ce}_\cz$-pre-image of order $n=\ord(P)$,
there is a zooming pre ball $V_n(x)$, $x\in{\tz_n(\alpha,\delta,f)}$, containing $P$ and $F|_P$ can be extended to
a homeomorphism  between $\overline{P}$ and $\overline{\Delta}$ (because $f^n|_{\overline{V_n(x)}}$ is a homeomorphism).
Given $x\in\bigcap_{n\ge0}F^{-n}\big(\bigcup_{P\in\cp}P\big)$, set $P_j=\cp(F^j(x))$.
As
$\mbox{diameter}(\cp_n(x))=\mbox{diameter}({F}|_{P_1}^{-1}\circ{F}|_{P_2}^{-1}\circ\cdots\circ{F}|_{P_n}^{-1}({\Delta}))$
$<$ $\prod_{j=1}^{n}\alpha_{_{\ord(P_j)}}(\diameter(\Delta))\le\alpha_{_{\sum_{j=1}^n\ord(P_j)}}(\diameter(\Delta))\to 0$, we conclude that $\cp$ is a Markov partition for $F$.
Finally, as $\{R>0\}=\bigcup_{P\in\cp}P$ and $F(P)=\Delta$ $\forall\,P\in\cp$, it follows that the Markov map $(F,\cp)$ is indeed an induced full Markov map.\cqd

Let $\mu$ be a $(\alpha,\delta)$-weak zooming measure with $\mu(X\setminus\Lambda)=0$ and let $U\subset X$ be an $\mu$ ergodic component.
Let ${A}$ be the attractor associated to
$U$ and ${A}_{\widetilde{\fz}}\subset{A}$ the compact
set such that ${\omega}_{f,{\widetilde{\fz}}}(x)={A}_{\widetilde{\fz}}$ for $\mu$-almost every point $x\in U$ (given by Proposition~\ref{PropositionErgodicComponents2} and by Lemma~\ref{LemmaOmegaLimitSets} applied to ${\cu}=\widetilde{\fz}$).

\begin{Lemma}\label{LemmaPartition}Let $(F,\cp)$ be as in Corollary~\ref{CorollaryDisjointness} and suppose that $\Delta\cap{A}_{\widetilde{\fz}}\ne\emptyset$. Then
$(F,\cp)$ is an induced full Markov map defined on $\Delta$ and it is compatible with $\mu|_{U}$.
\end{Lemma}
\dem
Let $p\in \Delta\cap{A}_{\widetilde{\fz}}$.
As $p\in{\omega}_{f,{\widetilde{\fz}}}(x)$ for $\mu$ almost every $x\in U$, we get $\mu|_{U}(U\setminus \bigcup_{n\ge0}f^{-n}(\Delta))=0$. Thus, as $\mu|_{U}\circ f^{-1}\ll\mu|_{U}$, $\mu|_{U}(\Delta)>0$.

By Corollary~\ref{CorollaryDisjointness}, we only need to show that $\mu|_{U}(\Delta\setminus\bigcup_{P\in\cp}P)$
$=$ $\mu\big((\Delta\setminus\bigcup_{P\in\cp}P)\cap U\big)=0$.
As  $p\in{\omega}_{f,{\widetilde{\fz}}}(x)$ for $\mu$ almost every $x\in U$, it follows that ${\Omega}(x)\ne\emptyset$ for
$\mu$ almost every $x\in \Delta$. Thus, $\mu|_{U}(\{R=0\})$
$=$ $\mu|_{U}(\Delta\setminus\bigcup_{P\in\cp}P)=0$.
\cqd

\begin{Theorem}\label{TheoremLocalinducedMaps}Suppose that for some $0<r_0<\delta/2$ and every $x$ the $(\alpha,\delta)$-zooming nested open ball $B_{r_0}^{\star}(x)$ is well defined and contains $B_{r_0/2}(x)$. Let $\Lambda$ $\subset$ $\limsup_{n}\tz(\alpha,\delta,f)$ be a positively invariant set and ${\widetilde{\fz}}=(\widetilde{\fz}(x))_{x\in\Lambda}$ a proper $(\alpha,\delta)$-zooming sub-collection. Let $\mu$ be a $(\alpha,\delta)$-weak zooming probability with bounded distortion and $\mu(\Lambda)=1$. Let $U\subset X$ an ergodic component for $\mu$ and ${A}_{\widetilde{\fz}}$ be the compact set such that ${\omega}_{f,{\widetilde{\fz}}}(x)={A}_{\widetilde{\fz}}$ for $\mu$-almost every point $x\in U$ (given by Theorem~\ref{PropositionFatZoomingErgodicAttractors}). Let $\Delta$ be a $(\alpha,\delta)$-zooming nested open ball with $\diameter(\Delta)\le r_0/2$ and such that $\Delta\cap{A}_{\widetilde{\fz}}\ne\emptyset$.

If $(F,\cp)$ is the induced Markov map associated to {\em ``the first $\widetilde{\ce}_{\cz}$-return time to $\Delta$}'' (as in Corollary~\ref{CorollaryDisjointness})
then $(F,\cp)$ is an induced full Markov map with $\mu$-bounded distortion. Furthermore, there exists $\nu\ll\mu$ an ergodic $F$-invariant probability with
$\log\frac{d\nu}{d\mu}\in L^{\infty}(\mu|_{\{\frac{d\nu}{d\mu}>0\}})$ and $\nu(\Delta)=1$.
\end{Theorem}
\dem
Let us show that, as $\mu$ has bounded distortion, $\mu|_U(\Delta)=\mu(\Delta)$. To prove this,
let $p\in \Delta\cap{A}_{\widetilde{\fz}}$. By Lemma~\ref{LemmaErgodicComponents}, $\mu(B_{{r_0}/2}(p)\cap U)=\mu(B_{{r_0}/2}(p))$. As $\diameter(\Delta)\le {r_0}/2$, $\Delta\subset B_{{r_0}/2}(p)$. So, $\mu|_U(\Delta)=\mu(\Delta)$.

As $\mu|_U(\Delta)=\mu(\Delta)$, Lemma~\ref{LemmaPartition} implies that $(F,\cp)$ is an induced full Markov map defined on $\Delta$ compatible with $\mu$.

Finally, as $\left|\log\frac{J_\mu F(x)}{J_\mu F(y)}\right| \leq \rho \dist(F(x),
F(y))$, $\forall\,x, y\in P$ and $\forall\,P\in\cp$ (because $P$ is contained in a zooming pre-ball of order $R(P)$ and $\mu$ has bounded distortion at the zooming times), we obtain that $(F,\cp)$ has $\mu$-bounded distortion.

Applying Proposition~\ref{FolkloreTheorem}, we obtain a $F$-invariant ergodic probability
$\nu\ll\mu$ with $\log\frac{d\nu}{d\mu}\in L^{\infty}(\mu|_{\{\frac{d\nu}{d\mu}>0\}})$ and, of course,
$\nu(\Delta)=1$.
\cqd

Given $\theta>0$ and $n\in\NN$, let $\cz_n(\alpha,\delta,\theta,f)$ be the set of points $x\in X$
such that $\#\{1\le j\le n$ $;$ $x\in\tz_j(\alpha,\delta,f)\}\ge\theta n$. Thus,
the set of points of $X$ with infinitely many moments with $\theta$-frequency of $(\alpha,\delta)$-zooming
times (with respect to $f$) is
$$\limsup_n\cz_n(\alpha,\delta,\theta,f)=\bigcap_{j=1}^{+\infty}\bigcup_{n\ge j}\cz_n(\alpha,\delta,\theta,f).$$

If $\mu$ is a $(\alpha,\delta)$-zooming measure with bounded distortion, $X$ can be decomposed in a finite collection of $\{U_1,...,U_s\}$ of $\mu$-ergodic components (Theorem~\ref{PropositionFatZoomingErgodicAttractors}). By ergodicity, $\exists\,\theta_i>0$ such that $$\limsup\frac{1}{n}\#\{1\le j\le n\,;\,\,x\in\tz_j(\alpha,\delta,f)\}\ge\theta_i$$ for $\mu$ almost every $x\in U_i$ $\forall\,i$. Furthermore, if ${\widetilde{\fz}}=(\widetilde{\fz}(x))_{x\in\Lambda}$ is a proper zooming sub-collection and $\mu(X\setminus\Lambda)=0$, there are also $\widetilde{\theta}_1,\cdots,\widetilde{\theta}_s>0$ such that
$$\limsup\frac{1}{n}\#\{1\le j\le n\,;\,\mbox{$j$ is a $\widetilde{\fz}$-time to }x\}\ge\widetilde{\theta}_i$$
for $\mu$ almost every $x\in U_i$ and all $1\le i\le s$.
Thus, we get the following remark.

\begin{Remark}\label{Remark7676}Let $\Lambda$ be a zooming set and ${\widetilde{\fz}}=(\widetilde{\fz}(x))_{x\in\Lambda}$ a proper zooming sub-collection. Let $\mu$ be a zooming measure with $\mu(X\setminus\Lambda)=0$. If $\mu$ has bounded distortion or, more in general, has a finite number of ergodic components then $\exists\,\widetilde{\theta}>0$ such that $$\limsup\frac{1}{n}\#\{j\le n\,;\,\mbox{ is a $\widetilde{\fz}$-time to }x\}\ge\widetilde{\theta}$$
for $\mu$ almost every $x\in X$. In particular, for every zooming measure $\mu$ with bounded distortion (or having a finite number of ergodic components) there is $\theta>0$ such that
$$\mu\big(X\setminus\limsup_m\cz_m(\alpha,\delta,\theta,f)\big)=0.$$

\end{Remark}

\begin{Theorem}\label{TheoremLocalLift}
Suppose that for some $0<{r_0}<\delta/2$ and every $x$ the $(\alpha,\delta)$-zooming nested
open ball $B_{r_0}^{\star}(x)$ is well defined and contains $B_{{r_0}/2}(x)$.
Let $\Lambda\subset X$ be a $(\alpha,\delta)$-zooming set and $\mu$ an ergodic
$f$-invariant zooming probability with $\mu(\Lambda)=1$.
Let ${\widetilde{\fz}}=(\widetilde{\fz}(x))_{x\in\Lambda}$ be a proper
$(\alpha,\delta)$-zooming sub-collection and ${A}_{+,{\widetilde{\fz}}}$ the compact
set such that ${\omega}_{+,f,{\widetilde{\fz}}}(x)={A}_{+,{\widetilde{\fz}}}$
for $\mu$-almost every point $x\in X$ (given by Lemma~\ref{LemmaOmegaLimitSets}
applied to ${\cu}={\widetilde{\fz}}$).
Let $\Delta$ be a $(\alpha,\delta)$-zooming nested open ball with $\diameter(\Delta)\le {r_0}/2$ and such that $\Delta\cap{A}_{+,{\widetilde{\fz}}}\ne\emptyset$.

If $R$ is {\em ``the first $\widetilde{\ce}_{\cz}$-return time to $\Delta$}'' and
$(F,\cp)$ is the induced Markov map associated to $R$ (as in Corollary~\ref{CorollaryDisjointness})
then $(F,\cp)$ is a full induced Markov map compatible with $\mu$ and there exists a
$F$-invariant finite measure $\nu\ll\mu$ (indeed, $\nu(Y)\le\mu(Y)$ for every Borel
set $Y\subset\Delta$) such that $\int R d\nu<+\infty$ and
$${\mu}=\frac{1}{\gamma}\sum_{j=0}^{+\infty}{f}\,^j_{*}(\nu|_{\{{{R}}>j\}}),$$
where $\gamma=\sum_{j=0}^{+\infty}{f}\,^j_{*}(\nu|_{\{{{R}}>j\}})(X)$.
\end{Theorem}
\dem
Let $A_{{{\widetilde{\fz}}}}$ be the compact set (given by  Lemma~\ref{LemmaOmegaLimitSets}) such that ${\omega}_{{f},{{\widetilde{\fz}}}}(x)=A_{{{\widetilde{\fz}}}}$
for ${\mu}$-almost every $x\in X$.
As ${A}_{+,{\widetilde{\fz}}}\subset A_{{{\widetilde{\fz}}}}$, we have $\Delta\cap{A}_{{\widetilde{\fz}}}\ne\emptyset$. Thus,
it follows from Lemma~\ref{LemmaPartition} that $(F,\cp)$ is an induced full Markov map defined on $\Delta$ and compatible with ${\mu}$  (see Definition~\ref{DefinitionMarkovMapCompatibleWithAMeasure}).

Let $\mathfrak{B}=\{x\in \Delta\,{;}\,F^j(x)\in \bigcup_{P\in{\cp}}P,\,\forall\,j\ge 0\}$. Because $\mu$ is $f$-invariant (in particular, $f$-non-singular), we get $\Delta=\mathfrak{B}$ (mod ${\mu}$).
As $\Delta\cap A_{_{+},{\widetilde{\fz}}}\ne \emptyset$ and ${\mu}$ is ${f}$-ergodic,
there is $\Theta>0$ such that
$$\limsup_{n\to\infty}\frac{1}{n}\#\{1\le j\le n\,;\,x\in G_j\mbox{ and }{f}\,^j(x)\in\Delta\}\ge\Theta$$ for ${\mu}$ almost every $x\in \Delta$, where $G_j=\{x\in\Lambda$ $;$ $j$ is a $\widetilde{\fz}$-time to $x\}$.
Thus, taking $B=\Delta$, $g=R$
and applying the first part of Lemma~\ref{LemmaIntegrabilidade} to ${f}$ we get
\begin{equation}\label{EquationLocalLift}
\limsup_{n\to\infty}\frac{1}{n}\#\{j\ge0\,;\,\sum_{k=0}^{j}R
\circ {{F}}^k(x)\le n\}\ge\Theta
\end{equation}
for ${\mu}$ almost every $x\in \Delta$.
Because $$\{j\ge0\,;\,\sum_{k=0}^{j}R
\circ {{F}}^k(x)<n\}=\{0\le j< n\,;\,{f}^j(x)\in\co_{{F}}^+(x)\},$$
it follows form (\ref{EquationLocalLift}) and Theorem~\ref{TheoremRogers} that
there exists a non trivial ${F}$-invariant measure such that $\nu(Y)\le\mu(Y)$ for every Borel set $Y\subset\Delta$ (in particular, $\nu\ll{\mu}$) with $\int R d \nu<+\infty$. Thus,
$\eta=\sum_{j=0}^{+\infty}{f}\,^j_{*}(\nu|_{\{{{R}}>j\}})$
is a ${f}$-invariant finite measure (see Remark~\ref{RemarkProjection}).
Note that, if $\eta(Y)>0$ for some Borel set $Y\subset X$ then $\nu(f^{-j}(Y))>0$ for some $j\ge0$ and, as $\nu\ll\mu$, $\mu(Y)=\mu(f^{-j}(Y))>0$. Thus, $\eta\ll\mu$.
As ${\mu}$ is ${f}$-ergodic probability, we get
$${\mu}=\frac{1}{\eta(X)}\,\eta=\frac{1}{\eta(X)}
\sum_{j=0}^{+\infty}{f}\,^j_{*}(\nu|_{\{{{R}}>j\}}).$$

\cqd

\begin{Lemma}\label{LemmaZuzun}
For every $k\ge1$,
$$\limsup_m\cz_m(\{\alpha_{n}\}_n,\delta,\theta,f)
\subset
\bigcup_{j=0}^{k-1}f^{-j}\big(\limsup_m\cz_{k\,m}(\{\alpha_{n}\}_n,
\delta,{\theta}/{k},f)\big)\subset$$
$$\subset\bigcup_{j=0}^{k-1}f^{-j}\big(\limsup_m\cz_m(\{\alpha_{k n}\}_n,\delta,{\theta}/{k},f^k)\big).$$
\end{Lemma}
\dem
Let ${k}\ge 1$. For each $x\in\limsup_m\tz_m(\{\alpha_{n}\}_n,\delta,f)$ and $0\le i<k$, set $\NN_x(i)=\{{k} j+i$ ${;}$ $j\in\NN$ and $x\in\tz_{{k} j+i}(\{\alpha_n\}_n,{\delta},f)\}$.
So, $x\in \tz_{m}(\{\alpha_n\}_n,{\delta},f)$ $\Leftrightarrow$ $m\in\bigcup_{i=0}^{{k}-1}\NN_x(i)$.
Note also that $\NN_x(j)\cap\NN_x(i)=\emptyset$, whenever $i\ne j$.

So, for each $x\in\limsup\cz_m(\{\alpha_{n}\}_n,\delta,\theta,f)$ $\subset$
$\limsup\tz_m(\{\alpha_{n}\}_n,\delta,f)$ one can choose
$\ell(x)\in\{0,...,{k}-1\}$ such that $\limsup_m\frac{1}{m}\#\{1\le j\le m$ ${;}$ $j\in\NN_x(\ell(x))\}\ge\theta/{k}$. Otherwise,
$\limsup_m\frac{1}{m}\#\{1\le j\le m$ ${;}$ $x\in
\tz_{j}(\{\alpha_n\}_n,{\delta},f)\}<\theta$, contradicting $x\in\limsup\cz_m(\{\alpha_{n}\}_n,$ $\delta,\theta,f)$.
As $j\in\NN_x(\ell(x))$ $\Leftrightarrow$ $x\in\tz_{j {k}+\ell(x)}(\{\alpha_n\}_n,{\delta},f)$
$\Leftrightarrow$ $f^{\ell(x)}(x)\in\tz_{j k}(\{\alpha_n\}_n,{\delta},f)$, it follows from Lemma~\ref{LemmaConcatenar} that
$$j\in\NN_x(\ell(x)) \Rightarrow f^{\ell(x)}(x)\in\tz_{k\,j}(\{\alpha_{ n}\}_n,{\delta},f)\subset\tz_{j}(\{\alpha_{k n}\}_n,{\delta},f^{k}).$$
Therefore,
$$\limsup_m\frac{1}{m}\#\{1\le j\le m\,{;}\,f^{\ell(x)}(x)\in\tz_{k\,j}(\{\alpha_{n}\}_n,{\delta},f)\}\ge$$
$$\ge\limsup_m\frac{1}{m}\#\{1\le j\le m\,{;}\,f^{\ell(x)}(x)\in\tz_{j}(\{\alpha_{k n}\}_n,{\delta},f^{k})\}\ge\theta/{k}.$$
As a consequence, if $x\in\limsup_m\tz_m(\{\alpha_{n}\}_n,\delta,f)$ then $f^{\ell(x)}(x)\in\limsup_m\mathbb\cz_m(\{\alpha_{k n}\}_n,$ ${\delta},\theta/{k},f^{k})$, with
$0\le \ell(x)<{k}$.
\cqd

\begin{Corollary}\label{CorollaryCuringa}
Let $\mu$ be a $f$-ergodic $(\alpha,\delta)$-zooming measure (not necessarily
invariant). For each $k>0$ there is a positively invariant set $E\subset X$ with $\mu$ positive measure and
such that $\frac{1}{\mu(E)}\mu|_{_{E}}$
is an $(\{\alpha_{k\,n}\},\delta)$-zooming ergodic probability with
respect to $f^k$.
Furthermore, if $\mu$ is $f$-invariant then $E$ is a $\mu$-ergodic component with respect to $f^k$, $\frac{1}{\mu(E)}\mu|_{_{E}}$ is $f^k$-invariant and $\mu=(\sum_{j=0}^{k-1}{f^j}_*){\mu|_E}$.
\end{Corollary}
\dem
It follows from Remark~\ref{Remark7676} that $\exists\,\theta>0$ such that  $\mu\big(X\setminus\limsup_{n}\cz_n({\alpha},\delta,{\theta},{f})\big)=0$.
By Lemma~\ref{LemmaZuzun}, there is $0\le j<k$ such that $\mu(f^{-j}(\widetilde{\cz}))>0$,
where $\widetilde{\cz}= \limsup_{n}\cz_n\big($ $\{\alpha_{k\,n}\}$, $\delta$, $\theta/k$, ${f^k}\,\big)$.
As $\mu$ is $f$-non-singular (by the definition of zooming measure),
we get $\mu(\widetilde{\cz})>0$.

Because $\mu$ has at most $k$ ergodic components with respect to ${f^k}$ (Lemma~\ref{LemmaReletingErgNumber}), there is one of these ergodic components ${E_0}\subset X$ such that $\mu({E_0}\cap\widetilde{\cz})>0$. Set $E={E_0}\cap\widetilde{\cz}$. As $\widetilde{\cz}\supset {f^k}(\widetilde{\cz})$ and ${f^k}^{-1}({E_0})={E_0}$, we have ${f^k}(E)$ $\subset$ $E$ $\subset$ $E_0$. Because $\mu|_{_{E_0}}$ is ${f^k}$-ergodic and ${f^k}$-non-singular, it follows from Lemma~\ref{LemmaContructingNon-SingErgMeasure} that
$\frac{1}{\mu(E)}\mu|_{_{E}}=\frac{1}{\mu(E)}(\mu|_{_{E_0}})|_{_E}$ is a probability ${f^k}$-non-singular and ${f^k}$-ergodic.
Of course $\frac{1}{\mu(E)}\mu|_{_{E}}(\widetilde{\cz})=1$ and so, $\frac{1}{\mu(E)}\mu|_{_{E}}$ is an $(\{\alpha_{k\,n}\},\delta)$-zooming ergodic probability with respect to ${f^k}$.

Suppose now that $\mu$ is $f$-invariant. In this case, as $E\subset{f}^{-k}(E)$,
it follows that $E={f}^{-k}(E)$ (mod $\mu$). Thus, changing $E$ by $\bigcap_{j\ge0}f^{-j\,k}(E)$,
it follows that $E$ is a $\mu$-ergodic component with respect to $f^k$.
So, it is easy to conclude that $\mu|_{_{E}}$ is ${f^k}$-invariant
and $\mu=(\sum_{j=0}^{k-1}{f^j}_*){\mu|_E}$.
\cqd

\begin{maintheorem}[Existence of invariant zooming measures]\label{TheoremZoomingInvariant}
If $\mu$ is a zooming measure with bounded
distortion then there exists a finite collection of ergodic
$f$-invariant probabilities absolutely continuous with respect to $\mu$
such that $\mu$-almost every
point in ${X}$ belongs to the basin of one of these probabilities.
\end{maintheorem}
\dem[{\bf Proof of Theorem~\ref{TheoremZoomingInvariant}}]
By Theorem~\ref{PropositionFatZoomingErgodicAttractors},
${X}$ can be partitioned in a finite collection of $\mu$-ergodic components with respect to $f$. Let $U$ be one of these ergodic components.
Choose any $0<r_0<\delta/2$ and let $k\ge 1$ be such that $\sum_{n=1}^{+\infty}\alpha_{k n}(\widetilde{r})<\widetilde{r}/4$ for $\widetilde{r}=r_0/4$ and for $\widetilde{r}=r_0$.

By Lemma~\ref{LemmaReletingErgNumber}, $U$ can be decomposed into a finite collection of disjoint
ergodic components with respect to $f^k$. As $U$ is invariant (in particular, $f(U)\subset U$), it follows from Lemma~\ref{LemmaContructingNon-SingErgMeasure} that $\mu|_{_U}$ is $f$-non-singular. Thus, $\mu|_{_U}$ is an ergodic $(\alpha,\delta)$-zooming measure. From Corollary~\ref{CorollaryCuringa}, there is $E\subset U\subset X$, with $f(E)\subset E$ and $\mu|_{_U}(E)>0$,
such that $\widetilde{\mu}=\frac{1}{\mu(E)}\mu|_{_E}=\frac{1}{\mu|_{_U}(E)}(\mu|_{_U})|_{_E}$ is a $(\widetilde{\alpha},\delta)$-zooming ergodic probability with respect to $\widetilde{f}=f^k$, where $\widetilde{\alpha}=\{\alpha_{k\,n}\}_n$.

Denote the set of $(\widetilde{\alpha},\delta)$-zooming images of $\widetilde{f}$ by ${\widetilde{\fz}}=({\widetilde{\fz}}(x))_{x\in\Lambda}$,
where ${\widetilde{\fz}}(x)$ $=$ $\{\widetilde{f}\,^n(x)\,;$ $n\in\NN$ and $x\in{\tz_n(\widetilde{\alpha},\delta,\widetilde{f})}\}$ is the set of
$(\widetilde{\alpha},\delta)$-zooming images of $x$ by $\widetilde{f}$ and $\Lambda$ $=$ $\limsup\tz_n(\widetilde{\alpha},\delta,\widetilde{f})$.

By Theorem~\ref{PropositionFatZoomingErgodicAttractors},
there exists a fat attractor ${A}$ (with respect to $\widetilde{f}$) such that
${{\omega}_{\widetilde{f}}}(x)=A$ for $\widetilde{\mu}$-almost every
point $x\in X$.
Moreover, there are compact sets
$A_{_{+},{\widetilde{\fz}}},A_{\widetilde{\fz}}\subset A$, with $A_{_{+},{\widetilde{\fz}}}\subset A_{\widetilde{\fz}}$,
such that ${\omega}_{{\widetilde{f}},{\widetilde{\fz}}}(x)=A_{\widetilde{\fz}}$ and ${\omega}_{+,{\widetilde{f}},{\widetilde{\fz}}}(x)=A_{_{+},{\widetilde{\fz}}}$
for $\widetilde{\mu}$-almost every point $x\in X$.

Let $r=r_0/4$ and choose any point $q\in A_{_{+},{\widetilde{\fz}}}$.
As $A_{_{+},{\widetilde{\fz}}}\subset A_{\widetilde{\fz}}$, we get $B^{\star}_r(q)\cap A_{\widetilde{\fz}}\supset B^{\star}_r(q)\cap A_{_{+},{\widetilde{\fz}}}\ne \emptyset$, where $B^{\star}_r(q)$ is the $(\widetilde{\alpha},\delta)$-zooming nested ball with respect to $\widetilde{f}$, radius $r$ and center on $q$ (see Definition~\ref{DefinitionZoomingNestedBalls} and Lemma~\ref{LemmaZoomingNestedBall}).

Taking $\Delta=B^{\star}_r(q)$ and $\widetilde{\ce}$ as the collection of all $(\widetilde{\alpha},\delta)$-zooming pre-balls with respect to $\widetilde{f}$ (see Notation~\ref{NotationCollectionZoomingPreBalls}), let $R$ be the first $\widetilde{\ce}$-return time to $\Delta$ (with respect to $\widetilde{f}$) given by (\ref{EquationInducing}), let $F=\widetilde{f}^R$ be induced map associated to $R$ and let $\cp$ be the Markov partition given by (\ref{EquationIndMArkovCollection}).

Applying Theorem~\ref{TheoremLocalinducedMaps} to $\widetilde{f}$, $\widetilde{\alpha}$, $\Delta=B^{\star}_r(q)$ (note that $\diameter(\Delta)=r_0/2$), $(F,\cp)$ and $\widetilde{\mu}$, we obtain a $F$-invariant measure $\nu\ll\widetilde{\mu}$.
To prove the existence of a $\widetilde{f}$-invariant ergodic probability $\widetilde{\eta}\ll\widetilde{\mu}$ we need only to show that the induced time $R$ is $\nu$ integrable
(see Proposition~\ref{FolkloreTheorem}).

Let $\mathfrak{B}=\{x\in B^{\star}_r(q)\,{;}\,F^j(x)\in \bigcup_{P\in\cp}P\,\forall\,j\ge 0\}$.
Note that $\widetilde{\mu}(B^{\star}_r(q)\setminus \mathfrak{B})=\nu(B^{\star}_r(q)\setminus \mathfrak{B})=0$.
Let $B$ be the set of $x\in \mathfrak{B}$ such that
$\limsup_n\frac{1}{n}\#\{1\le j\le n$ ${;}$ $x\in{\tz_j(\widetilde{\alpha},\delta,{\widetilde{f}})}$ and ${\widetilde{f}}^j(x)\in B\}>0$.
As $B^{\star}_r(q)\cap A_{_{+},{\widetilde{\fz}}}\ne \emptyset$, $\nu(B\setminus \mathfrak{B})=\widetilde{\mu}(B\setminus \mathfrak{B})=0$.

Taking $G_j={\tz_j(\widetilde{\alpha},\delta,{\widetilde{f}})}$, $g=R$ and $T=F$, it follows from Lemma~\ref{LemmaIntegrabilidade} that
$$\liminf_n\frac{1}{n}\sum_{j=0}^{n-1}R\circ F^j(x)<+\infty$$ for every $x\in B$. By Birkhoff's Ergodic Theorem,
$\int R d\nu=\lim_n\frac{1}{n}\sum_{j=0}^{n-1}R\circ F^j(x)$ for $\nu$-almost every $x\in B$.
Thus, $\int R d\nu<+\infty$.
As a consequence, the projection $\widetilde{\eta}=\sum_{P\in\cp}\sum_{j=0}^{R(P)-1}{{\widetilde{f}}^j}_*(\nu|_P)$ is a $\widetilde{\mu}$ absolutely continuous $\widetilde{f}$-invariant finite measure.

Taking $\eta=\frac{1}{k}\sum_{j=0}^{k-1}f^j_*\widetilde{\eta}$, it is easy to see that
$\eta$ is ${f}$-invariant finite measure and  ${\eta}\ll{\mu}$.
So, to finish the proof of the theorem we only need to verify that $U$ belongs to the basin of $\eta$.

By Birkhoff's Theorem, $\eta(\mathcal{B}(\eta)\cap U)=\eta(\mathcal{B}(\eta))>0$ and, as $\eta\ll\mu$, we get $\mu(\mathcal{B}(\eta)\cap U)>0$. As $\mathcal{B}(\eta)$ is a $f$ invariant set and $U$ is a $\mu$ ergodic component with respect to $f$, we conclude that $U=\mathcal{B}(\eta)$ (mod $\mu$).

\cqd

Before we begin the proof of Theorem~\ref{TeoremMarkovStructure} which gives the existence of a Markov structure for a zooming set,
we want to emphasize a difference between the proof of Theorem~\ref{TheoremZoomingInvariant} and proof of Theorem~\ref{TeoremMarkovStructure}.

In both proofs we begin with a reference measure $\mu$ and we need to show the existence of an induced invariant measure $\nu\ll\mu$ and also the $\nu$ integrability of the inducing time $R$.
In the hypothesis of Theorem~\ref{TheoremZoomingInvariant}, we have a zooming measure $\mu$
with bounded distortion, but we do not know if $\mu$ is
invariant. On the other hand, in Theorem~\ref{TeoremMarkovStructure}, we want to
study zooming measures for which we do not know anything about distortion,
but we know that they are invariant measures. In the proof of
Theorem~\ref{TheoremZoomingInvariant} the existence of $\nu$ is given by
Proposition~\ref{FolkloreTheorem} (this proposition is used to prove Theorem~\ref{TheoremLocalinducedMaps}) and in the proof of Theorem~\ref{TeoremMarkovStructure}
the existence of $\nu$ is assured by Theorem~\ref{TheoremRogers} (this theorem is central in the proof of Theorem~\ref{TheoremLocalLift}).
In both case, the  estimate
to get the integrability of the inducing time is given by Lemma~\ref{LemmaIntegrabilidade} (this lemma
appears in the proof of Theorem~\ref{TheoremLocalinducedMaps}~and~\ref{TheoremLocalLift}).

\begin{maintheorem}[Markov Structure for the zooming set]\label{TeoremMarkovStructure}Every zooming set $\Lambda$ admits a finite Markov structure $\mathfrak{F}=\{(F_1,\cp_1),...,(F_s,\cp_s)\}$. Furthermore, each $(F_i,\cp_i)\in\mathfrak{F}$ is a full Markov map defined on some topological open ball $U_i$ (also the elements of $\cp_i$ are topological open balls).
\end{maintheorem}

\dem[{\bf Proof of Theorem~\ref{TeoremMarkovStructure}}]If $f$ is backward separated and $\sup_{r>0}\sum_{n\ge1}\alpha_n(r)/r<\infty$,
choose any $0<r\le\frac{1}{4}r_0$, where $r_0$ is given by Lemma~\ref{LemmaZoomingNestedBall}, and set $m_{0}=1$. Otherwise, choose $0<r_0<\delta/2$, set $r=r_0/4$ and let $m_{0}$ be an integer big enough such that $\sum_{n\ge1}\alpha_{m_{0} n}(\widetilde{r})$ $<$ $\widetilde{r}/4$ for $\widetilde{r}=r_0$ and $\widetilde{r}=r$. Set also $\widetilde{f}=f^{m_{_{0}}}$, $\widetilde{\theta}=\theta/m_{_{0}}$, $\widetilde{\alpha_j}=\alpha_{m_{_{0}}\,j}$ and $\widetilde{\alpha}=\{\widetilde{\alpha_j}\}_j$.

Let $\widetilde{\fz}$ be the proper zooming sub-collection for $\widetilde{f}$ given by
$\widetilde{\fz}=(\widetilde{\fz}(x))_{x\in\widetilde{\Lambda}}$,
where $$\widetilde{\fz}(x)=\{\widetilde{f}\,^n(x)\,{;}\,n\in\NN\mbox{ and }x\in{\tz_{m_0\,n}(\alpha,\delta,f)}\}$$ and $$\widetilde{\Lambda}=\Lambda\cap\limsup_{n\to\infty}\tz_{m_0\,n}(\alpha,\delta,f)
\subset\Lambda\cap\limsup_{n\to\infty}\tz_n
(\widetilde{\alpha},\delta,\widetilde{f}).$$
One can easily check that $\widetilde{\fz}$ is indeed a proper $(\widetilde{\alpha},\delta)$-zooming collection with respect to $\widetilde{f}$ (see comments just below Definition~\ref{DefinitionPropreZooColl}).

As $X$ is compact, one can find $q_1,q_2,...,q_s\in X$ such that $\{B_r^{\star}(q_1),$ $...,$ $B_r^{\star}(q_s)\}$ is a finite cover of $X$ by $(\widetilde{\alpha},\delta)$-zooming nested balls, with respect to $\widetilde{f}$.

For each $1\le j\le s$, let $\widetilde{F}_j:B_r^{\star}(q_j)\to B_r^{\star}(q_j)$ be the induced map (with respect to $\widetilde{f}$) on
$B^{\star}_r(q_j)$ associated to {\em ``the first $\widetilde{\ce}$-return time to $B^{\star}_r(q_j)$}'' given by (\ref{EquationinducedMap}), where $\widetilde{\ce}$  is the collection of all $\widetilde{\fz}$-pre-balls $V_n(x,\widetilde{f})$ (with respect to $\widetilde{f}$) such that $x\in\Lambda$ and $n\ge1$ is a $\widetilde{\fz}$-time (with respect to $\widetilde{f}$) for $x$.
Note that $V_n(x,\widetilde{f})$ is a $(\widetilde{\alpha},\delta)$-zooming pre-ball of order $n$ with respect to $\widetilde{f}$ and also a $(\alpha,\delta)$-zooming pre-ball of order $m_0\,n$ with respect to $f$

Let $\widetilde{R}_j$ be the inducing time of $\widetilde{F}_j$
(with respect to $\widetilde{f}$\,)
and let $\cp_j$ be the Markov partition associated to
{\em ``the first $\widetilde{\ce}$-return time to $B_r^{\star}(q_j)$}''
(see (\ref{EquationIndMArkovCollection})).

Set, for every $x\in B_r^{\star}(q_j)$, $F_j(x)=f^{m_{_{0}}\widetilde{R}_j(x)}(x)$. As $\widetilde{f}=f^{m_{_{0}}}$, it is easy to see that $(F_j,\cp_j)$  is also an induced full Markov map with respect to $f$ with inducing time $R_j=m_{0} \widetilde{R}_j$. Of course, as it is a Markov map,
$(F_j,\cp_j)=(\widetilde{F}_j,\cp_j)$.

\begin{Remark}\label{RemarkAddtionalProperties}
Note that $\mathfrak{F}=\{({F}_1,\cp_1),...,({F}_s,\cp_s)\}$ is a finite collection of induced full Markov maps for ${f}$ and satisfies the following additional properties.
\begin{enumerate}
\item Each $F_j$ is defined on a topological open ball. Indeed, each $F_j$ is defined on a
$(\widetilde{\alpha},\delta)$-zooming nested ball with respect to $\widetilde{f}$.
\item For every $1\le j\le s$, each $P\in\cp_j$ is a topological open ball.
Furthermore, setting $n=\widetilde{R_j}|_P$,
we get $P\subset V_n(x,\widetilde{f})$ for some $x\in P\cap\tz_{m_{_{0}}\,n}(\alpha,\delta,f)\cap\Lambda$, where $V_n(x,\widetilde{f})$ is a $(\alpha,\delta)$-zooming pre-ball of order $m_0\,n$ with respect to $f$ (and also a $(\widetilde{\alpha},\delta)$-zooming pre-ball of order $n$ with respect to $\widetilde{f}$).
\end{enumerate}
\end{Remark}

Thus, to finish the proof we only need to show that $\mathfrak{F}$ is a Markov structure for $\Lambda$.

Let $\mu$ be an ergodic $f$-invariant probability such that $\mu(\Lambda)>0$.
By ergodicity, $\mu(\Lambda)=1$ (we are also using that $\Lambda$ is positively invariant and $\mu$ is invariant). As $\mu$ is $f$-non-singular (because $\mu$ is $f$-invariant) and $\Lambda$ is a $(\alpha,\delta)$-zooming set, it follows that $\mu$ is a $(\alpha,\delta)$-zooming measure.

It follows from Corollary~\ref{CorollaryCuringa}
that there is a $\mu$-ergodic component $U\subset X$ with respect to $\widetilde{f}$ such that $\widetilde{\mu}=\frac{1}{\mu(U)}\mu|_U$ is a $(\widetilde{\alpha},\delta)$-zooming ergodic invariant
probability with respect to $\widetilde{f}$ and $\mu=(\sum_{j=0}^{m_{_0}-1}{f^j}_*){\mu|_U}$.

By Proposition~\ref{PropositionErgodicComponents2}, there exists a
$\widetilde{f}$-attractor $A\subset {X}$ which attracts  $\widetilde{\mu}$-almost every point of $X$ and such that ${{\omega}_{\widetilde{f}}}(x)=A$ for $\widetilde{\mu}$-almost every $x\in X$ (indeed, as
$\widetilde{\mu}$ is $\widetilde{f}$-invariant, $A=\supp\widetilde{\mu}$).
By Lemma~\ref{LemmaOmegaLimitSets},
there are compact sets $A_{_{+},\widetilde{\fz}}$ and $A_{\widetilde{\fz}}$, with
$A_{_{+},\widetilde{\fz}}\subset A_{\widetilde{\fz}}\subset A$, such that
${\omega}_{\widetilde{f},\widetilde{\fz}}(x)=A_{\widetilde{\fz}}$ and
${\omega}_{_{+},\widetilde{f},\widetilde{\fz}}(x)= A_{_{+},\widetilde{\fz}}$ for $\widetilde{\mu}$-almost
every $x\in X$ (see Definition~\ref{DefinitionOmegaU} and \ref{DefinitionOmegaUTheta}).

Let $1\le j_{_{0}}\le s$ be such that $B_r^{\star}(q_{j_{_{0}}})\cap A_{_{+},\widetilde{\fz}}\ne\emptyset$. It follows from Theorem~\ref{TheoremLocalLift} that $(\widetilde{F}_{j_{_{0}}},\cp_{j_{_{0}}})$ is an induced full Markov map (with respect to $\widetilde{f}$) defined on $B_r^{\star}(q_{j_{_{0}}})$ and compatible with $\widetilde{\mu}$.
As a consequence, $(F_{j_{_{0}}},\cp_{j_{_{0}}})$ is an induced full Markov map with respect to $f$ (defined on $B_r^{\star}(q_{j_{_{0}}})$ and compatible with $\widetilde{\mu}$).
Also by Theorem~\ref{TheoremLocalLift},  there exists  a $\widetilde{F}_{j_{_{0}}}$-invariant measure $\nu\ll\widetilde{\mu}$ such that
$$\widetilde{\mu}=\frac{1}{\widetilde{\gamma}}
\sum_{j=0}^{+\infty}\widetilde{f}\,^j_{*}(\nu|_{\{{\widetilde{R}_{j_{_{0}}}}>j\}}),$$
where $\widetilde{\gamma}=\sum_{j=0}^{+\infty}\widetilde{f}\,^j_{*}
(\nu|_{\{{\widetilde{R}_{j_{_{0}}}}>j\}})(X)$.

It follows from the relation $R_{j_{_{0}}}=m_{_{0}} {\widetilde{R}_{j_{_{0}}}}$ that
\begin{equation}\label{OBS}
\{R_{j_{_{0}}}>m_{_{0}}\,j+k\}=\{R_{j_{_{0}}}>m_{_{0}}\,j\},
\end{equation}
$\forall$ $0\le k<m_{_{0}}$ and $\forall\,j\ge 0$.

Setting $\gamma=\widetilde{\gamma}/\mu(U)$, we get
$$\mu=(\sum_{k=0}^{m_{_{0}}-1}{f^k}_*){\mu|_U}=\mu(U)(\sum_{k=0}^{m_{_{0}}-1}{f^k}_*)\widetilde{\mu}=$$
$$=\frac{\mu(U)}{\widetilde{\gamma}}(\sum_{k=0}^{m_{_{0}}-1}{f^k}_*)\sum_{j=0}^{+\infty}\widetilde{f}\,^j_{*}
(\nu|_{\{{\widetilde{R}_{j_{_{0}}}}>j\}})
=\frac{1}{\gamma}(\sum_{k=0}^{m_{_{0}}-1}{f^k}_*)\sum_{j=0}^{+\infty}{f}^{m_{_{0}}\,j}_{*}
(\nu|_{\{{{R}_{j_{_{0}}}}/m_{_{0}}>j\}})=$$
$$=\frac{1}{\gamma}\sum_{k=0}^{m_{_{0}}-1}\sum_{j=0}^{+\infty}{f}^{m_{_{0}}\,j+k}_{*}
(\nu|_{\{{{R}_{j_{_{0}}}}>m_{_{0}}\,j\}})\underbrace{=}_{(\ref{OBS})} \frac{1}{\gamma}\sum_{j=0}^{+\infty}\sum_{k=0}^{m_{_{0}}-1}{f}^{m_{_{0}}\,j+k}_{*}
(\nu|_{\{{{R}_{j_{_{0}}}}>m_{_{0}}\,j+k\}})=$$
$$=\frac{1}{\gamma}\sum_{n=0}^{+\infty}{f}^{n}_{*}(\nu|_{\{{{R}_{j_{_{0}}}}>n\}}),$$
finishing the proof of Theorem~\ref{TeoremMarkovStructure}.

\cqd

The Corollary~\ref{CorollaryMarkovStructure} follows from the Remark~\ref{RemarkAddtionalProperties}.

\begin{Corollary}\label{CorollaryMarkovStructure}
If $\mathfrak{F}=\{(F_1,\cp_1),...,(F_s,\cp_s)\}$ is the induced Markov structure given by the Theorem~\ref{TeoremMarkovStructure} and, for each $1\le j\le s$, $R_j$ is the induced time of $F_j$
then $\mathfrak{F}$ has the following additional properties.
\begin{enumerate}
\item There is some $m_0\ge1$ such that each $F_j$ is defined on an $(\{\alpha_{m_0\,n}\}_n,\delta)$-zooming nested ball $B^*_j$ with respect to $f^{m_0}$.
\item Each $P\in\cp_j$ is a topological open ball, $\forall\,1\le j\le s$.
Furthermore, there is a $(\alpha,\delta)$-pre-ball $V_n$ with respect to $f$, where $n={R_j}|_P$,
such that $P=(f^n|_{V_n})^{-1}(\cp_j)\subset V_n$. In particular,
$$\dist(f^{\ell}(x),f^{\ell}(y))\le\bigg(\sum_{k>n-{\ell}}\alpha_k\bigg)\dist(F_j(x),F_j(y))\,\,\,\forall\,x,y\in P.$$
\item If $\mu$ is an ergodic (not necessarily invariant) zooming measure
with bounded distortion and $\mu(B^*_j)>0$ then
$\exists\,K>0$ such that
$$\left|\log\frac{J_\mu F_j(x)}{J_\mu F_j(y)}\right| \leq K \dist(F_j(x),
F_j(y)),$$ for $\mu$ almost every $\,x, y\in P$, all $\,P\in\cp_j$ and $1\le j\le s$.
\end{enumerate}
\end{Corollary}

%%%%%%%%%%%%%%%%%%%%%%%%%%%%%%%%%%%%%%%%%%%%%%%%%%%%%%%%%%%%%%%%%%%%%%%%%%%%%%

\section{A global induced Markov map}\label{SectionGlobalInducedMarkovMap}

In \cite{ALP}, Alves Luzzatto and Pinheiro study the decay of correlations of non-uniformly expanding maps using a local induced Markov map. Using a global induced map, Gouëzel \cite{G} could improve the results of \cite{ALP} to deal with decay faster then super-polynomially. The advantage of a global induced map is the possibility of dominating the induced time by the first hyperbolic time.
In this section we construct a global induced map (adapted to any zooming measure) with the induced time smaller or equal to the first zooming time with respect $f$ or, when we do not have enough backward contraction, with respect to a fixed iterate of the original map.

It was introduced in Section~\ref{SectionNestedSets} the notion of an open set being nested and this notion can be extended straightforwardly to {\em essentially open sets} (definition in Section~\ref{DefinitionEssentiallyOpenSets}). For this, consider a connected, compact, separable metric space $X$.

\begin{Definition}[Essentially open linked sets] We  say that two essentially open sets $U_1$ and $U_2$ are {\em linked} if their interior are linked.
\end{Definition}

Let $f:X\to X$ be a measurable map and let $\ce=(\ce_n)_n$ be a dynamically closed  family of pre-images.

Exactly as we have done to open sets, we say that a collection of essentially open sets $\ca$
is a {\em $\ce$-nested collection} if every
$A\in\ca$ is not linked with any $\ce$-pre-image of an element of $\ca$ with order bigger than zero.

Note that a collection $\ca=\{A_1,\dots,A_s\}$ of essentially open sets is $\ce$-nested if and only if the collection $\interior(\ca):=\{\interior(A_1),\dots,\interior(A_s)\}$ is a $\ce$-nested collection of sets as in Definition~\ref{DefinitionNestedCollection}. As a consequence, we get the following remark.
\begin{Remark}\label{RemarkMainNestedProperty}
Lemma~\ref{LemmaPreMainNestedProperty} is also valid for collections of essentially open sets.
That is, if $\ca$ is an $\ce$-nested collection of essentially open sets and $P_1$ and $P_2$
are $\ce$-pre-images of two elements of $\ca$ with $\ord(P_1)\ne\ord(P_2)$
then $P_1$ and $P_2$ are not linked.
\end{Remark}

Let $\delta>0$ and let $\alpha=\{\alpha_n\}_{1\le n\in\NN}$ be a zooming contraction
(Definition~\ref{DefinitionZoomingContration}).

\begin{maintheorem}[Global zooming induced Markov map]
\label{GlobalMarkovStructureForZooming} Given an $(\alpha,\delta)$-zooming set $\Lambda\subset X$ there are an induced Markov map $(F,\cp)$ defined on $X$ with induced time $R$, a finite partition $\cp_0$ of $X$ by essentially open sets and an integer $\ell\ge1$ satisfying the following properties.
\begin{enumerate}
\item\label{ItemGMSFZ1} For each $Q\in\cp$ there exists $P\in\cp_0$ such that $\interior(Q)\subset P$.
\item\label{ItemGMSFZ2} $F(P)\in\cp_0$ $\forall\,P\in\cp$(in particular, the elements of $P$ are essentially open sets).
\item\label{ItemGMSFZ3} Given $P\in\cp$ there is a zooming pre-ball $V_{R(P)}(x)$, $x\in \tz_{R(P)}(\alpha,\delta,f)\cap\Lambda$, such that $F|_P=\bigg(f^{R(P)}\big|_{V_{R(P)}(x)}\bigg)\bigg|_P$. In particular,
    \begin{enumerate}
    \item[(\ref{ItemGMSFZ3}.1)] $\dist(F(x),F(y))\ge 8\dist(x,y)$ $\forall\,x,y\in P$ and $\forall\,P\in\cp$;
    \item[(\ref{ItemGMSFZ3}.2)]\label{ItenZ4.2} for all $x,y\in P$, $P\in\cp$ and $0\le n< R(P)$, $$\dist(f^n(x),f^n(y))<\alpha_{_{R(P)-n}}\dist(F(x),F(y));$$
    \item[(\ref{ItemGMSFZ3}.3)] if $\mu$ is $(\alpha,\delta)$-zooming measure (not necessarily
        invariant) with bounded distortion (with respect to $f$) and $\mu(X\setminus\Lambda)=0$ then
        $\exists\,\rho>0$ such that
        $$\left|\log\frac{J_\mu F(x)}{J_\mu F(y)}\right|\le\rho\dist(F(x),
        F(y)),$$ for $\mu$ almost every $x, y\in P$ and $\forall\,P\in\cp$.
    \end{enumerate}
\item\label{ItemGMSFZ4} There is a good relationship between the tail of the partition and the tail of zooming times, i.e., $$\{R>n\}\cap\Lambda\subset \Lambda\setminus \bigcup_{j=1}^n \tz_{\ell\,j}(\alpha,\delta,f)\subset \Lambda\setminus \bigcup_{j=1}^n \tz_{j}(\widetilde{\alpha},\delta,f^{\ell}),$$
    where $\widetilde{\alpha}=\{\alpha_{\ell\,n}\}_n$.
\item\label{ItemGMSFZ5}If $\mu$ is a $f$-invariant measure with $\mu(\Lambda)>0$ then $\mu\big|_{\bigcap_{j\ge0}f^{-\ell\,j}(\{R>0\})}$ is an invariant measure with respect to $f^{\ell}$.
\item Every ergodic $f$-invariant zooming probability $\mu$ with $\mu(\Lambda)>0$ is liftable to $F$.
\end{enumerate}
\end{maintheorem}

\dem[{\bf Proof of Theorem~\ref{GlobalMarkovStructureForZooming}}]

Choose $0<r_{0}<\delta/2$ and set $r=r_0/4$. Let $\ell$ be an integer big enough such that $\sum_{j=1}^{\infty}\alpha_{\ell\,j}(\widetilde{r})\le \widetilde{r}/8$ for $\widetilde{r}=r_0$ and $\widetilde{r}=r$. Set $\widetilde{f}=f^{\ell}$, $\widetilde{\theta}=\theta/\ell$, $\widetilde{\alpha_j}=\alpha_{\ell\,j}$ and $\widetilde{\alpha}=\{\widetilde{\alpha_j}\}_j$.
As in the proof of Theorem~\ref{TeoremMarkovStructure}, let $\widetilde{\fz}$ be the proper zooming sub-collection for $\widetilde{f}$ given by
$\widetilde{\fz}=(\widetilde{\fz}(x))_{x\in\widetilde{\Lambda}}$,
where $$\widetilde{\fz}(x)=\{\widetilde{f}\,^n(x)\,{;}\,n\in\NN\mbox{ and }x\in{\tz_{\ell\,n}(\alpha,\delta,f)}\}$$ and $$\widetilde{\Lambda}=\Lambda\cap\limsup_{n\to\infty}\tz_{\ell\,n}(\alpha,\delta,f)
\subset\Lambda\cap\limsup_{n\to\infty}\tz_n
(\widetilde{\alpha},\delta,\widetilde{f}).$$

It is clear that $\widetilde{\Lambda}$ is $\widetilde{f}$-positively invariant. Also note that $\widetilde{\Lambda}$ is a large portion of $\Lambda$. Indeed, it follows from (\ref{EquationPeqZZ}) that $$\Lambda\subset\widetilde{\Lambda}\cup f^{-1}(\widetilde{\Lambda})\cup\cdots\cup f^{-(\ell-1)}(\widetilde{\Lambda}).$$ As a consequence $$\mu(\Lambda)>0\,\,\Longrightarrow\,\,\mu(\widetilde{\Lambda})>0$$ for every $f$-non-singular measure $\mu$.

Denote by  $\ce_{\cz,\widetilde{f}}=(\ce_{\cz,\widetilde{f},n})_n$, where $\ce_{\cz,\widetilde{f},n}=\{V_n(x,\widetilde{f})$
${;}$ $x\in\tz_n(\widetilde{\alpha},\delta,\widetilde{f})\}$ is the collection of all $(\widetilde{\alpha},\delta)$-zooming pre-balls with respect $\widetilde{f}$ of order $n$.
Let us denote  the order with respect to $\widetilde{f}$ by $\ord_{\widetilde{f}}$.

Let $\widetilde{\ce}\subset\ce_{\cz,\widetilde{f}}$ be the collection of all $\widetilde{\fz}$-pre-ball $V_n(x,\widetilde{f})$ (with respect to $\widetilde{f}$) for all $x\in\Lambda$ and all $\widetilde{\fz}$-time $n$ for $x$ with respect to $\widetilde{f}$.

Let $\{q_1,...,q_s\}\subset X$ be a maximal ${r}/2$ -separated set and consider the collection $\ca=\{B_{{r}}(q_1),...,B_{{r}}(q_s)\}$ of open balls. As we have contraction in zooming times, the elements of $\ca$ are not contained in any ${\ce_{\cz,\widetilde{f}}}$~-pre-image of an element of $\ca$ with  order bigger than zero (see Notation~\ref{NotationCollectionZoomingPreBalls}). As $\sum_n\widetilde{\alpha}_n\le 1/8$, every chain of ${\ce_{\cz,\widetilde{f}}}$~-pre-images of $\ca$ has diameter smaller that ${r}/4$, that is, if $(P_0,...,P_n)\in ch_{{\ce_{\cz,\widetilde{f}}}}(B_{r}(q_i))$ then $$\diameter(\bigcup_j P_j)<\sum_j\diameter(P_j)<\sum_j\alpha_j \diameter(B_{{r}}(q_i))<{r}/4.$$
Thus, $$(B_{{r}}(q_i))^{\star}=B_{r}(q_i)\setminus\Bigg(\overline{\bigcup_{(P_j)_j\in ch_{{\ce_{\cz,\widetilde{f}}}}(B_{r}(q_i))}\bigcup_j P_j}\Bigg)\supset B_{(3/4)r}(q_i),$$ for all $1\le i\le s$.

Let $\ca'=\{\Delta_1,...,\Delta_s\}$, where $\Delta_i$ is the connected component of $(B_{{r}}(q_i))^{\star}$ containing $B_{(3/4){r}}(q_i)$.
It follows from Proposition~\ref{PropositionExisteNested2} that $\ca'$ is a $\ce_{\cz,\widetilde{f}}$~-nested collection of sets.
Moreover, as $\{q_1,...,q_s\}\subset X$ is maximal ${r}/2$ -separated, $\ca'$ is a cover of $X$ by opens sets.

Setting $\cp^s=\{{\Delta_1}\cap\cdots\cap{\Delta_s}\}$ and, for $1\le\wp<s$,
$$\cp^\wp=\bigg\{{\Delta_{i_1}}\cap\cdots\cap {\Delta_{i_\wp}}\setminus
\bigg(\bigcup_{1\le k_1<\cdots<k_{\wp+1}\le s}{\Delta_{k_1}}\cap\cdots\cap {\Delta_{k_{\wp+1}}}
\bigg)\,\mbox{\LARGE{;}}\,1\le i_1<\cdots<i_\wp\le s\bigg\},$$ it follows that $\cp_0:=\cp^1\cup\cdots\cup\cp^{s}$ is a partition of $X$ by essentially open sets. Note that, $\bigcup_{P\in\cp_0}\partial P\subset\bigcup_{j=1}^s\partial \Delta_j$.

\begin{Claim}\label{ClaimOKJJUHG6}
Let $Q$ be a $\ce_{\cz,\widetilde{f}}$~-pre-image (with respect to $\widetilde{f}$) of some $\Delta_i\in\ca'$ with $\ord_{\widetilde{f}}(Q)>0$ and let $P\in\cp_0$. If $Q\cap P\ne\emptyset$  then $Q\subset \interior(P)$.
\end{Claim}
\dem[Proof of Claim~\ref{ClaimOKJJUHG6}]
Suppose that $Q\cap P\ne\emptyset$ and $Q\not\subset \interior(P)$. As $Q$ is a connected open set, $Q\cap\partial P\ne\emptyset$. Thus, there exists $\Delta_k\in\ca'$ such that $Q\cap\partial \Delta_k\ne\emptyset$. As $\ca'$ is  $\ce_{\cz,\widetilde{f}}$~-nested collection of open sets and $\ord_{\widetilde{f}}(Q)>0$, $\Delta_k\subset Q$. As $\diameter(Q)<\big(\sum_{j=1}^{\ord_{\widetilde{f}}(Q)}\widetilde{\alpha}_j\big)\diameter(\Delta_{i})<{r}/8$.
But this leads to a contradiction because $B_{(3/4){r}}(p_{k})\subset \Delta_{k}$.
\cqd

\begin{Claim}\label{Claim878jhjh}
If $Q_1$ and $Q_2$
are $\ce_{\cz,\widetilde{f}}$~-pre-images (with respect to $\widetilde{f}$) of respectively $P_1,P_2\in\cp_0$ then $Q_1$ and $Q_2$ are not linked. In particular, $\cp_0$ is an $\ce_{\cz,\widetilde{f}}$~-nested collection (with respect to $\widetilde{f}$) of essentially open sets.
Furthermore,
\begin{enumerate}
\item if $Q_1\cap Q_2\ne\emptyset$ then $\ord_{\widetilde{f}}(Q_1)\ne\ord_{\widetilde{f}}(Q_2)$;
\item if $Q_1\subsetneqq Q_2$ then $\ord_{\widetilde{f}}(Q_1)>\ord_{\widetilde{f}}(Q_2)$.
\end{enumerate}
\end{Claim}
\dem[Proof of Claim~\ref{Claim878jhjh}]
Note that if $Q_1\cap Q_2\ne\emptyset$ then $\ord_{\widetilde{f}}(Q_1)\ne\ord_{\widetilde{f}}(Q_2)$. Otherwise $P_1\cap P_2=f^{\ord_{\widetilde{f}}(Q_1)}(Q_1\cap Q_2)\ne\emptyset$. This shows the first item.

Suppose that  $\interior(Q_1\cap Q_2)\ne\emptyset$, with $Q_1\ne Q_2$. We may assume that $P_1\ne P_2$ (if $P_1= P_2$, the claim follows from Corollary~\ref{CorollaryMainNestedProperty}). Set $\wp_j=\ord_{\widetilde{f}}(P_j)$ (with respect to $\widetilde{f}$) for $j=1,2$.
By the first item, assume for instance that $\wp_1<\wp_2$.

Write $Q_2=(f^{\wp_2}|{V_{\wp_2}(x,\widetilde{f})})^{-1}(P_2)$, with  $x\in\tz_{\wp_2}(\widetilde{\alpha},\delta,\widetilde{f})$. Let $\Delta_{j_2}\in\ca'$ be such that $\interior(P_2)\subset \Delta_{j_2}$ and set $\widetilde{\Delta}_{j_2}=(f^{\wp_2}|{V_{\wp_2}(x,\widetilde{f})})^{-1}(\Delta_{j_2})$.
Thus,
\begin{equation}\label{EquationGlobalPartition1}
\interior(P_1\cap \widetilde{f}^{\wp_1}(\widetilde{\Delta}_{j_2}))\supset \interior(P_1\cap
\widetilde{f}^{\wp_1}(Q_2))=\widetilde{f}^{\wp_1}(\interior(Q_1\cap Q_2))\ne\emptyset.
\end{equation}

As $\widetilde{f}^{\wp_1}(\widetilde{\Delta}_{j_2})$ is a
$\ce_{\cz,\widetilde{f}}$~-pre-image of $\Delta_{j_2}$ with order $\wp_2-\wp_1>0$, it follows from Claim~\ref{ClaimOKJJUHG6}
that $\interior(P_1)\supset \widetilde{f}^{\wp_1}(\widetilde{\Delta}_{j_2})$.
So, $\interior(P_1)\supset \interior(\widetilde{f}^{\wp_1}(Q_2))$. Using (\ref{EquationGlobalPartition1}), we get $$\widetilde{f}^{\wp_1}(\interior(Q_1))=\interior(P_1)\supset \widetilde{f}^{\wp_1}(\interior(Q_1\cap Q_2)).$$ As a consequence, $\interior(Q_1)\supset\interior(Q_2)$.

So, we obtain that $\interior(Q_1\cap Q_2)\ne\emptyset$ and $\wp_1<\wp_2$ implies that $\interior(Q_1)\supset\interior(Q_2)$ (or, if $\wp_1>\wp_2$, $\interior(Q_1)\subset\interior(Q_2)$). From this we conclude that $Q_1$ and $Q_2$
are not linked and also the second item of the claim.
\cqd

Define an inducing time $\widetilde{R}:X\to\{0,1,2,\dots\}$ on $X$
as follows.
Given
$x\in X$, let ${\Omega}(x)$ be the collection of
all $\widetilde{\ce}$-pre-images $Q$ of any  $P\in\cp_0$ such that $x\in Q$. That is, $Q\in\Omega(x)$ if $x\in Q$ and
there are $n\in\NN$, $y\in\widetilde{\Lambda}$ and $P\in\cp_0$ such that  $Q=(\widetilde{f}\,^{n}|_{V_{n}(y,\widetilde{f})})^{-1}(P)$, where $n\ge1$ is a $\widetilde{\fz}$-time (with respect to $\widetilde{f}$) for $y$. Note that $V_n(y,\widetilde{f})$ is both a $(\widetilde{\alpha},\delta)$-zooming pre-ball of order $n$ with respect to $\widetilde{f}$ and a $(\alpha,\delta)$-zooming pre-ball of order $\ell\,n$ with respect to $f$.
If ${\Omega}(x)\ne\emptyset$ let
$\widetilde{R}(x)=\min\{\ord_{\widetilde{f}}(V)\,{;}$ $V\in{\Omega}(x)\}$ and let
$\widetilde{R}(x)=0$ whenever ${\Omega}(x)=\emptyset$.

Note that if $x\in\widetilde{\Lambda}$ then $\widetilde{R}(x)$ is smaller than or equal to the first $\widetilde{\fz}$-time of $x$, i.e., $\widetilde{R}(x)$ $\le$ $\min\{n$ ${;}$ $n$ is a $\widetilde{\fz}$-time (with respect to $\widetilde{f}$) to $x\}$ $=$ $\min\{n$ ${;}$ $x\in{\tz_{\ell\,n}(\alpha,\delta,f)}\}$. Thus,
\begin{equation}\label{EquationinducedGlobalTime}
\{\widetilde{R}>n\}\cap\Lambda\subset\Lambda\setminus \bigcup_{j=1}^n \tz_j(\widetilde{\alpha},\delta,\widetilde{f})\subset\Lambda\setminus \bigcup_{j=1}^{n} \tz_{\ell\,j}(\alpha,\delta,f).
\end{equation}

Define the induced map $\widetilde{F}$ on
$X$ associated to {\em the first $\ce_{\widetilde{\cz}}$ time} by
\begin{equation}\label{EquationinducedGlobalMap}
\widetilde{F}(x)=\widetilde{f}\,^{\widetilde{R}(x)}(x),\,\forall x\in X.
\end{equation}

If $\Omega(x)\ne\emptyset$, it follows from Claim~\ref{Claim878jhjh} that the collection of sets ${\Omega}(x)$ is totally ordered by inclusion. Moreover, there is a unique $Q\in{\Omega}(x)$ such that $\ord_{\widetilde{f}}(Q)=\widetilde{R}(x)$. In this case, set $I(x)=Q$.

Also by Claim~\ref{Claim878jhjh}, $\ord_{\widetilde{f}}(I(x))<\ord_{\widetilde{f}}(J)$ $\forall\,I(x)\ne J\in\Omega(x)$ and $\forall\,x\in X$. Furthermore, if $I(x)\cap I(y)\ne\emptyset$ then $I(x)=I(y)$ (see the proof of Lemma~\ref{LemmaDisjointness} which is analogous).

Proceeding as in the proof of Corollary~\ref{CorollaryDisjointness}, one can easily conclude that $$\mbox{$\cp:=\{I(x)$ $;$ $x\in X$ and $\Omega(x)\ne\emptyset\}$}$$ is a Markov partition for $\widetilde{F}$. Besides, defining $R(x)=\ell\,\widetilde{R}(x)$ and $F(x)=f^{R(x)}(x)$ $=$ $\widetilde{f}\,^{\widetilde{R}(x)}(x)$ $=$ $\widetilde{F}(x)$, one can see that $\cp_0,\cp,F$ and $R$ satisfy the first four items of the theorem.

\begin{Remark}\label{RemarkIndMarPar78}
For future references we note that $\cp|_{\widetilde{\Lambda}}$ $=$ $\{P\cap\widetilde{\Lambda}$ $;$ $P\in\cp\}$ is an induced Markov partition of $\widetilde{\Lambda}$ with respect to $\widetilde{f}=f^{\ell}$ (this follows from Claim~\ref{Claim878jhjh}).
\end{Remark}

Let $\mu$ be a $f$-invariant measure with $\mu(\Lambda)>0$. To check the item (\ref{ItemGMSFZ5}) set $E=\bigcap_{j}\widetilde{f}^{-j}(\{R>0\})$. As $\widetilde{f}^{-1}(E)\supset E\supset\widetilde{\Lambda}$, we get $\widetilde{f}^{-1}(E)=E$ (mod $\mu$) and $\mu(E)\ge\mu(\widetilde{\Lambda})>0$. Thus $\mu|_E$ is $\widetilde{f}$-invariant.

To prove the last item we will construct a local Markov map induced from $\widetilde{F}$.

\subsubsection*{Constructing a local induced map from the global one}Let $\mu$ be a $f$-invariant ergodic probability with $\mu(\Lambda)>0$.
By ergodicity, $\mu(\Lambda)=1$ (we are also using that $\Lambda$ is positively invariant and $\mu$ is invariant).
As $\widetilde{\Lambda}$ is $\widetilde{f}$-positively invariant and $\mu$ is also $\widetilde{f}$-invariant, $\mu|_{\widetilde{\Lambda}}$ is $\widetilde{f}$-invariant.
On the other hand, as $X$ can be decomposed into at most $\ell$
$\mu$-ergodic components with respect to $\widetilde{f}$,
there is a $\mu$-ergodic component $U\subset\widetilde{\Lambda}$.

Thus, $\widetilde{\mu}=\frac{1}{\mu(U)}\mu\big|_U$ is a $(\widetilde{\alpha},\delta)$-zooming ergodic invariant
probability with respect to $\widetilde{f}$, $\widetilde{\mu}(\widetilde{\Lambda})=1$ and $\mu=(\sum_{j=0}^{\ell-1}{f^j}_*)\mu|_U$.

By Proposition~\ref{PropositionErgodicComponents2}, there exists a
$\widetilde{f}$-attractor $A\subset {X}$ which attracts  $\widetilde{\mu}$ almost every point of $X$
(indeed $A=\supp\widetilde{\mu}$ because
$\widetilde{\mu}$ is $\widetilde{f}$-invariant).
By Lemma~\ref{LemmaOmegaLimitSets},
there are compact sets $A_{_{+},\widetilde{\fz}}$ and $A_{\widetilde{\fz}}$, with
$A_{_{+},\widetilde{\fz}}\subset A_{\widetilde{\fz}}\subset A$, such that
${\omega}_{\widetilde{f},\widetilde{\fz}}(x)=A_{\widetilde{\fz}}$ and
${\omega}_{_{+},\widetilde{f},\widetilde{\fz}}(x)= A_{_{+},\widetilde{\fz}}$ for $\widetilde{\mu}$-almost
every $x\in X$ (see Definition~\ref{DefinitionOmegaU} and \ref{DefinitionOmegaUTheta}).
Let $1\le j_{_{0}}\le s$ be such that
$\Delta_{j_{_{0}}}\cap A_{_{+},\widetilde{\fz}}\ne\emptyset$.

As in (\ref{EquationInducing}), we define the first $\widetilde{\ce}$-return time to $\Delta_{j_{_{0}}}$ (with respect to $\widetilde{f}$). Precisely, given
$x\in \Delta_{j_{_{0}}}$, let ${\Omega}_0(x)$ be the collection of
$\widetilde{\ce}$-pre-images $V$ of $\Delta_{j_{_{0}}}$ such that $x\in V$ (i.e.,
$x\in V=(\widetilde{f}\,^n|_{V_n(y)})^{-1}(\Delta_{j_{_{0}}})$ for some $y\in\tz_{\ell\,n}(\alpha,\delta,f)\cap\Lambda$ and $n\in\NN$) and define the first $\widetilde{\ce}$-return time to $\Delta_{j_{_{0}}}$ as the map $\widetilde{R}_0:\Delta_{j_{_{0}}}\to\NN$ given by
\begin{equation}\label{EquationInducing2}
\widetilde{R}_0(x)=\begin{cases}
\min\{\ord(V)\,;\,V\in{\Omega_0}(x)\} & \text{ if }{\Omega_0}(x)\ne\emptyset\\
0 & \text{ if }{\Omega_0}(x)=\emptyset
\end{cases}.
\end{equation}
Thus, the induced map $\widetilde{F}_0$ on
$\Delta_{j_{_{0}}}$ (with respect to $\widetilde{f}$) associated to {\em ``the first $\widetilde{\ce}$-return time to $\Delta_{j_{_{0}}}$}'' is given by
\begin{equation}\label{EquationinducedMap2}
\widetilde{F}_0(x)=\widetilde{f}\,^{\widetilde{R}_0(x)}(x),\,\forall x\in \Delta_{j_{_{0}}}.
\end{equation}

We claim that $\widetilde{F}_0$ is also an induced map with respect to $\widetilde{F}$.
\begin{Claim}\label{ClaimGlobalInd1}
For each $x\in\Delta_{j_{_{0}}}$ there is $\widetilde{k}(x)\in\NN$
such that $\widetilde{F}_0(x)=\widetilde{F}\,^{\widetilde{k}(x)}(x)$.
\end{Claim}
\dem[Proof of Claim~\ref{ClaimGlobalInd1}]
By definition, $V\in \Omega_0(x)$ if and only if $\exists\,y\in\tz_{\ell\,n}(\alpha,\delta,f)\cap\Lambda$ such that $x\in V=(\widetilde{f}\,^n|_{V_n(y,\widetilde{f})})^{-1}(\Delta_{j_{_{0}}})$, where $n=\ord_{\widetilde{f}}(V)$ and $V_n(y,\widetilde{f})$ is the $(\alpha,\delta)$-zooming pre-ball with respect to $f$ of order $\ell\,n$ ``centered on $y$'' (as noted before, $V_n(y,\widetilde{f})$ is also a $(\widetilde{\alpha},\delta)$-zooming pre-ball with respect to $\widetilde{f}$ of order $n$ ``centered on $y$'').
If $P$ is the element of $\cp_0$ that contains $\widetilde{f}\,^n(x)$, we get
$P\subset B_{\delta}(\widetilde{f}\,^n(y))=\widetilde{f}\,^n(V_n(y,\widetilde{f}))$  (because $\diameter(P)\le r_0/2<\delta/4$). Thus, $V':=(\widetilde{f}\,^n|_{V_n(y,\widetilde{f})})^{-1}(P)\in\Omega(x)$ and $\ord_{\widetilde{f}}(V')=n$. As a consequence, $0\le \widetilde{R}(x)\le \widetilde{R}_0(x)$ $\forall\,x\in\Delta_{j_{_{0}}}$.

Let $x\in\Delta_{j_{_{0}}}$ be such that $m:=\widetilde{R}_0(x)>0$. In this case
set $s=\sum_{n=0}^{\wp_x-1} \widetilde{R}\circ \widetilde{F}^n(x)$, where $$\wp_x=\max\{j\ge1\,;\,\sum_{n=0}^{j-1} \widetilde{R}\circ \widetilde{F}^n(x)< m\}.$$

Let $P\in\cp$ be such that $\widetilde{F}_0(x)\in P$ and let $y$ $\in$ $\tz_{\ell\,m}(\alpha,\delta,f)$ $\cap$ $\Lambda$ be such that $I_0(x)$ $=$ $(\widetilde{f}\,^m|_{V_{m}(y,\widetilde{f})})^{-1}(\Delta_{j_0})$, where $I_0(x)$ $\in$ $\Omega_0(x)$ is the unique element of $\Omega_0(x)$ such that $\ord_{\widetilde{f}}(I_0(x))$ $=$ $m$ $=$ $\widetilde{R}_0(x)$ (see the comment just above Lemma~\ref{LemmaDisjointness}). As $P$ $\subset$ $\overline{\Delta_{j_0}}$ $\subset$ $B_{\delta}(\widetilde{f}\,^m(y))$ $=$ $\widetilde{f}\,^m(V_m(y,\widetilde{f}))$,  $I:=(\widetilde{f}\,^m|_{V_{m}(y,\widetilde{f})})^{-1}(P)\in\Omega(x)$. Thus, $\widetilde{f}\,^s(I)\in\Omega(\widetilde{f}\,^{s}(x))$ and, as a consequence, $\widetilde{R}(\widetilde{f}\,^s(x))\le m-s$ (because $\ord_{\widetilde{f}}(\widetilde{f}\,^s(I))$ $=$ $m-s$). On the other hand, as $\widetilde{F}\,^{\wp_x}(x)$ $=$ $\widetilde{f}\,^s(x)$ and $\widetilde{R}(\widetilde{F}\,^{\wp_x}(x))+s$ $=$ $\sum_{n=0}^{\wp_x} \widetilde{R}\circ \widetilde{F}\,^n(x)$ $\ge$ $m$, we get $\widetilde{R}(\widetilde{f}\,^s(x))$ $\ge$ $m-s$.
Thus, $\widetilde{R}(\widetilde{f}\,^s(x))$ $=$ $m-s$. This implies that $\widetilde{R}_0(x)$ $=$ $m$ $=$ $\widetilde{R}(\widetilde{f}\,^s(x))+s$ $=$ $\sum_{n=0}^{\wp_x} \widetilde{R}\circ \widetilde{F}\,^n(x)$, i.e., $\widetilde{F}_0(x)$ $=$ $\widetilde{f}\,^{\widetilde{R}_0(x)}(x)$ $=$ $\widetilde{F}\,^{\widetilde{k}(x)}$, where $\widetilde{k}(x)$ $=$ $\wp_x+1$.\cqd

Now, to finish the proof, we will show the existence of
a $F$-invariant finite measure $\nu\ll\mu$ such that $\mu=\sum_{j=0}^{+\infty}f\,^j_{*}(\nu|_{\{R>j\}})$.

Because $\sum_{n\ge1}\widetilde{\alpha}_n(r_0)\le r_0/8$, the $(\widetilde{\alpha},\delta)$-zooming nested open ball $B_{r_0}^{\star}(x)$ (with respect to $\widetilde{f}$) is well defined and contains $B_{r_0/2}(x)$ for every $x\in X$ (Lemma~\ref{LemmaZoomingNestedBall}). As $\diameter(\Delta_{j_0})\le\frac{1}{2} r_0$, $\Delta_{j_0}\cap A_{_{+},{\widetilde{\fz}}}\ne\emptyset$ and $\widetilde{\mu}$ is an ergodic $\widetilde{f}$-invariant $(\widetilde{\alpha},\delta)$-zooming probability, we can apply Theorem~\ref{TheoremLocalLift} and obtain a finite $\widetilde{F}_0$-invariant finite measure $\nu_0\ll\widetilde{\mu}$ with $\int \widetilde{R}_0 d \nu_0<+\infty$ and such that $$\widetilde{\mu}=\frac{1}{\widetilde{\gamma}}
\sum_{j=0}^{+\infty}\widetilde{f}\,^j_{*}({\nu_0}|_{\{\widetilde{R}_0>j\}}),$$
where $\widetilde{\gamma}=\sum_{j=0}^{+\infty}\widetilde{f}\,^j_{*}(\nu_0|_{\{\widetilde{R}_0>j\}})(X)$.
As $\widetilde{F}_0(x)=\widetilde{F}^{\widetilde{k}(x)}(x)$ and $\int \widetilde{k}\,d\nu_0\le\int \widetilde{R}_0 d\nu_0$, it follows from Remark~\ref{RemarkProjection} that $\nu=\frac{1}{\widetilde{\gamma}}\sum_{j=0}^{+\infty}
\widetilde{F}\,^j_{*}(\nu_0|_{\{{{\widetilde{R}_0}}>j\}})$ is a $\widetilde{F}$-invariant finite measure (note that $\nu$ is not necessarily a probability).  Moreover, it is straightforward to check that $\sum_{n=0}^{+\infty}
\widetilde{f}\,^n_{*}(\nu|_{\{\widetilde{R}>n\}})$ $=$ $\frac{1}{\widetilde{\gamma}}$ $\sum_{n=0}^{+\infty}\widetilde{f}\,^n_{*}\big(\big($
$\sum_{j=0}^{+\infty}
\widetilde{F}\,^j_{*}({\nu_0}|_{\{{{\widetilde{R}_0}}>j\}})\big)
|_{\{\widetilde{R}>n\}}\big)$
$=$ $\frac{1}{\widetilde{\gamma}}$ $\sum_{n=0}^{+\infty}\widetilde{f}\,^n_{*}({\nu_0}
|_{\{\widetilde{R}_0>n\}})$ $=$ $\widetilde{\mu}$ (see, for instance,  Lemma~4.1 of \cite{Zw}).
That is $$\widetilde{\mu}=\sum_{n=0}^{+\infty}
\widetilde{f}\,^n_{*}(\nu|_{\{\widetilde{R}>n\}}).$$

Proceeding as in the end of the proof of Theorem~\ref{TeoremMarkovStructure} (or alternatively, using Lemma~4.1 of \cite{Zw}) we get $$\mu=\frac{1}{\gamma}\sum_{n=0}^{+\infty}{f}^{n}_{*}(\nu|_{\{{R}>n\}}),$$
where $\gamma=1/\mu(U)$.
\cqd

%%%%%%%%%%%%%%%%%%%%%%%%%%%%%%%%%%%%%%%%%%%%%%%%%%%%%%%%%%%%%%%%%%%%%%%%%%%%%%

\section{Expanding measures}
\label{ApplicationsExpandingMeasures}

Let $M$ be a compact Riemannian manifold of dimension $d\ge1$ and $f:M\to M$ is a non-flat map with a critical/singular set $\cc\subset M$.

{\bf Hyperbolic Times.}
The idea of hyperbolic times is a key notion
on the study of non-uniformly hyperbolic dynamics and was introduced
by Alves et al~\cite{Al,ABV}.
Let us fix $0<b<\frac{1}{2}\min\{1,1/\beta\}$. Given $0<\sigma<1$ and
$\varepsilon>0$, we will say that $n$ is a $(\sigma,\varepsilon)$-{\em
hyperbolic time} for a point $x\in M$ (with respect to the non-flat map $f$ with a
$\beta$-non-degenerate critical/singular set $\cc$) if for all $1\le k\le n$ we
have $\prod_{j=n-k}^{n-1}\|(Df\circ f^j(x))^{-1}\|\le
{\sigma}^k\mbox{ and }$ $\dist_{\varepsilon}(f^{n-k}(x),\mathcal{C})\ge
{\sigma}^{b k}$. We denote the set of points of $M$ such that $n\in\NN$ is a
$(\sigma,\varepsilon)$-{\em hyperbolic time} by $\th_n(\sigma,\varepsilon,f)$.

\begin{Proposition}\cite{ABV}\label{PropositionABV}
Let $U\subset M$ be a set
satisfying the slow approximation condition. Given $\lambda>0$ there exist
$\theta>0$ and $\varepsilon>0$
such that, for every $x\in U$, $$\#\{1\le j\le n\,{;}\,x\in
\th_j(e^{-\lambda/4},\varepsilon,f)\}\ge\theta\,n,$$ whenever
$\sum_{i=0}^{n-1}
\log \|(Df(f^{i}(x)))^{-1}\|^{-1}\ge\lambda\,n$.
\end{Proposition}

It follows from Proposition~\ref{PropositionABV} that the points
of an expanding set (recall
Definition~\ref{DefinitionExpandingSet}) have infinitely many
moments with positive frequency of hyperbolic times. In
particular, they have infinitely many hyperbolic times.

The proposition below assures that the hyperbolic times are indeed zooming times.

\begin{Proposition}\cite{ABV}
\label{PropositionHyperbolicBalls}
Given $\sigma\in(0,1)$ and
$\varepsilon>0$, there is $\delta>0$,
depending only on $\sigma,\varepsilon$ and on the map $f$, such
that if $x\in {\th}_n(\sigma,\varepsilon,f)$ then there exists a
neighborhood $V_n(x)$ of $x$ with the following
properties:
\begin{enumerate}
\item $f^{n}$ maps $\overline{V_n(x)}$ diffeomorphically onto the ball
$\overline{B_{\delta}(f^{n}(x))}$;
\item $dist(f^{n-j}(y),f^{n-j}(z)) \le
\sigma^{j/2}\dist(f^{n}(y),f^{n}(z))$ $\forall y, z\in V_n(x)$ and $1\le j<n$;
\end{enumerate}
\end{Proposition}

The sets $V_n(x)$ are called \emph{hyperbolic
pre-balls} and their images $f^{n}(V_n(x))=B_{\delta}(f^n(x))$,
\emph{hyperbolic balls}.

Given $\sigma\in(0,1)$, $\varepsilon>0$ and $\theta\in(0,1]$,
define ${\ch}_n(\sigma,\varepsilon,\theta,f)$ as the set $x\in
{\th}_n(\sigma,\varepsilon,f)$ such that $\#\{1\le j\le
n\,{;}\,x\in \th_j(\sigma,\varepsilon,f)\}\ge\theta\,n$.
\begin{Remark}\label{RemarkHH}
It follows from Proposition~\ref{PropositionABV} that if $x$ is a
$\lambda$-expanding point then there are $\varepsilon>0$ and
$\theta\in(0,1)$ such that
$x\in\limsup{\ch}_n(e^{-\lambda/4},\varepsilon,\theta,f)$. That
is, every $\lambda$-expanding point $x$ belongs not only to
$\limsup {\th}_n(e^{-\lambda/4},$ $\varepsilon,f)$ but also  to
$\limsup{\ch}_n(e^{-\lambda/4},\varepsilon,\theta,f)$. In
particular, if $\mu$ is a $\lambda$-expanding measure then there
exists $\varepsilon>0$ and $\theta\in(0,1)$ such that
$$\mu(M\setminus\limsup{\ch}_n(e^{-\lambda/4},\varepsilon,\theta,f))=0.$$
\end{Remark}

The proof of Lemma~\ref{LemmaHyperbolicBalls} just below is easy
and straightforward. For instance, replacing $\det Df$ by
$J_{\mu}f$, the proof proceeds exactly as in the Lebesgue case of
Proposition~2.5. of \cite{Pi}.

\begin{Lemma}
\label{LemmaHyperbolicBalls}
If $\mu$ is a $f$-non-flat measure then there is $\rho>0$ such that
$$
\bigg|\log \frac{J_\mu f^n(p)}{J_\mu f^n(q)}\bigg|\le \rho\,\dist(f^n(p),f^n(q))
$$
for every $x\in {\th}_n(\sigma,\varepsilon,f)$ and $\mu$-almost every $p$ and $q\in V_n(x)$.
\end{Lemma}

By Proposition~\ref{PropositionHyperbolicBalls}, if $n\in\NN$ is
$(\sigma,\varepsilon)$-hyperbolic time for $x\in {M}$ then $n$ is
an $(\{\alpha_n\}_n,\delta)$-zooming time to $x$, where $\alpha_n(r)=\sigma^{n/2}r$. Thus, using
Proposition~\ref{PropositionABV} and
\ref{PropositionHyperbolicBalls}, it follows that if $\lambda>0$
and $\ch$ is a $\lambda$-expanding set then there is $\delta>0$
such that $\ch$ is an $(\{\alpha_n\}_n,\delta)$-zooming
set, where $\alpha_n(r)=e^{-(\lambda/8)n}r$ (in particular, every $\lambda$-expanding measure is an
$(\{\alpha_n\}_n,\delta)$-zooming measure). Furthermore, using
Lemma~\ref{LemmaHyperbolicBalls}, we easily get the following
lemma.

\begin{Lemma}\label{LemmaZooExp}
Given $\lambda>0$ there is $\delta>0$ (depending only on $\lambda$
and $f$) such that every $f$-non-flat $\lambda$-expanding measure
is an $(\{\alpha_n\}_n,\delta)$-zooming measure with
bounded distortion at the zooming times, where $\alpha_n(r)=e^{-(\lambda/8)n}r$.
\end{Lemma}

\dem[{\bf Proof of Theorem~\ref{TheoremSRB}}]
 The proof of Theorem~\ref{TheoremSRB} follows
 straightforwardly from Lemma~\ref{LemmaZooExp}
 and Theorem~\ref{TheoremZoomingInvariant}.
\cqd

\begin{Remark}
If in Theorem~\ref{TheoremSRB} we set $\lambda=0$, then the
results will be the same with one difference only: the collection
of measures is not finite  but countable.
\end{Remark}
To prove the remark above, let $M'$ be the set of points $y\in M$
such that Equation~(\ref{EquationExpanding}) holds for every $x\in
\bigcup_{k=0}^{+\infty}f^{-k}(y)$. As $\mu\circ f^{-1}\ll\mu$,
$\mu(M\setminus M')=0$. For each $0<n\in\NN$, let $M_n$ be the set
of $x\in M'$ such that
$$
\limsup_{n\to\infty}\frac{1}{n}\sum_{i=0}^{n-1}
\log \|(Df(f^{i}(x)))^{-1}\|^{-1}\in \bigg(\frac{1}{n+1},\frac{1}{n}\bigg].
$$
Note that $M_n$ is an invariant set $\forall\,n\in\NN$, i.e,
$f^{-1}(M_n)=M_n$. Let $\mathcal{N}\subset\NN$ the set of
$n\in\NN$ such that $\mu(M_n)>0$. For each $n\in\mathcal{N}$, we
can apply Theorem~\ref{TheoremSRB} to $\mu|_{M_n}$. Thus, we only have
 to consider the collection of all measure $\nu$ such that
$\nu$ is $\mu|_{M_n}$ absolutely continuous ergodic $f$-invariant
probabilities for some $n\in\mathcal{N}$.

\dem[{\bf Proof of
Theorem~\ref{TeoremMarkovStructureForExpanding}}] Observe that
Theorem~\ref{TeoremExpandingMarkovStructure} below is a more
complete result than
Theorem~\ref{TeoremMarkovStructureForExpanding}. On the other
hand, the proof of Theorem~\ref{TeoremExpandingMarkovStructure} is
a direct consequence of Theorem~\ref{TeoremMarkovStructure},
Corollary~\ref{CorollaryMarkovStructure} and the fact that the
$\lambda$-expanding set $\ch$ is an $(\alpha,\delta)$-zooming set,
where $\alpha=\{\alpha_n\}_n$, $\alpha_n(r)=e^{-(\lambda/8)n}r$, $\delta$ is given by
Proposition~\ref{PropositionHyperbolicBalls} and $\theta$ is given
by Proposition~\ref{PropositionABV}. The case $\lambda=0$ follows
directly from the case $\lambda>0$ by taking any sequence
$\lambda_n\searrow0$ and setting the Markov structure as
$$\mathfrak{F}=\{(F,\cp)\,{;}\,(F,\cp)\in \mathfrak{F}(\lambda_n)\mbox{ and }n\in\NN\},$$
where $\mathfrak{F}(\lambda_n)$ is the Markov structure for $\lambda_n$.

\cqd

\begin{Theorem}[Markov Structure for an expanding set]
\label{TeoremExpandingMarkovStructure}
Let $M$ be a compact Riemannian manifold and $f:M\to M$ a non-flat
map. Let $\lambda\ge 0$ and $\ch$ be a $\lambda$-expanding set. Then
there is a Markov structure $\mathfrak{F}=\{(F_i,\cp_i)\}_i$ for
$\ch$. Furthermore, denoting the domain of $F_i$ by $U_i$ and its
inducing time by $R_i$, we have the following additional
properties.
\begin{enumerate}
\item If $\lambda>0$ then $\mathfrak{F}=\{(F_i,\cp_i)\}_i$ is a
finite collection; \item $U_i$ is a topological open ball
$\forall\,i$; \item Each $(F_i,\cp_i)\in\mathfrak{F}$ is a full
Markov map, i.e., is a Markov map with $F_i(P)=U_i$
$\forall\,P\in\cp_i$. In particular, every $P\in\cp_i$ is a
topological open ball $\forall\,i$; \item For each
$(F_i,\cp_i)\in\mathfrak{F}$ there is $\lambda_i>0$ such that
$$\log\|(DF_i(x))^{-1}\|^{-1}>\lambda_i,\,\,\forall\,x\in\bigcup_{P\in\cp_i}P.$$
\item Every $P\in\cp_i$  is contained in a hyperbolic pre-ball of
order $R_i(P)$, i.e., there is $x$ $\in$ ${\th}_{R_i(P)}$
$(e^{-\lambda_i/4},$ $\varepsilon,$ $f)$ such that $P\subset
V_{R_i(P)}(x)$, where $\varepsilon$ is given by
Proposition~\ref{PropositionABV}.
\end{enumerate}
\end{Theorem}

\dem[{\bf Proof of Theorem~\ref{GlobalExpandingInducedMarkovMap}}]
Let $\ell=\min\{k\in\NN\,;k\ge(16\log3)/\lambda\}$. Similar to the proof of
Theorem~\ref{TeoremMarkovStructureForExpanding}, the proof of
Theorem~\ref{GlobalExpandingInducedMarkovMap} is a consequence of
the zooming version of the theorem (in this case,
Theorem~\ref{GlobalMarkovStructureForZooming}) and the fact that
the $\lambda$-expanding set $\ch$ is an $(\alpha,\delta)$-zooming
set for any $\delta\in(0,\varepsilon]$, where
$\alpha=\{\alpha_n\}_n$, $\alpha_n(r)=e^{(-\lambda/8)n}r$ and
$\varepsilon>0$ is given by
Proposition~\ref{PropositionHyperbolicBalls}. Indeed, as
$\ell=\min\{k\in\NN\,;k\ge(16\log3)/\lambda\}$, we get
$\ell=\min\{k\in\NN$ $;$ $\sum_{n=1}^{\infty}\alpha_{k\,j}(r)\le
r/8\}$ $\forall\,r\ge0$. Thus, follows from
Theorem~\ref{GlobalMarkovStructureForZooming} the first nine items
of Theorem~\ref{GlobalExpandingInducedMarkovMap}.

Let us show the last item. If $\ch$ is uniformly expanding, then there are
$1<\tau<\rho<\infty$ such that $\tau<\|Df^n(x)^{-1}\|^{-1}<\rho$
$\forall\,x\in\ch$. This implies that both  $R$ and $DF$ is
bounded from above. Thus there is $\gamma>0$ such that
$\mbox{Leb}(P)>\gamma$ $\forall\,P\in\cp$ (because $F(P)\in\cp_0$
$\forall\,P\in\cp$ and $\cp_0$ is a finite collection of sets with
non-empty interior), where $\mbox{Leb}$ is a volume measure
associated to the Riemannian structure of M. From this we get
$\#\cp\le\mbox{Leb}(M)/\gamma<\infty$.
\cqd

%%%%%%%%%%%%%%%%%%%%%%%%%%%%%%%%%%%%%%%%%%%%%%%%%%%%%%%%%%%%%%%%%%%%%%%%%%%%%%

\section{Examples and applications}
\label{SectionExamplesAndAplications}

The purpose of the current section is to give examples of
expanding and zooming measures and also to give some illustrative
applications of the theorems and ideas previously developed.

For now, consider a  compact Riemannian manifold $M$ of dimension $d\ge1$.
Let $f:M\to M$ be a non-flat map  and $\cc\subset M$ its critical/singular set.

We say that a point $x\in M$ has all Lyapunov exponents positive if
\begin{equation}\label{EquationLYPFLD}
\limsup\frac{1}{n}\log|Df^n(x)\,v|>0\hspace{0.4cm}\forall\,0\ne v\in T_x M.
\end{equation}

A periodic  repeller is a periodic point $p$ such that $Df^n(p)$
is well defined and all eigenvalues of $Df^n(p)$ is bigger that
$1$, where $n$ is the period of $p$. In this case, this condition
is equivalent to all Lyapunov exponents of $p$ being positive.
Furthermore, as
\begin{equation}\label{EquationExpPerPoints}
\lim_{\ell\to\infty}\big\|\big((Df^{n}(p))^{-1}\big)^{\ell}\big\|^{\frac{1}{\ell}}
=\min\{\lambda^{-1}\,;\,\lambda\mbox{ is an eigenvalue of }Df^n(p)\},
\end{equation}
the periodic point $p$ is a repeller if and only if there is
$\ell>0$ such that $p$ is a periodic point for
$\widetilde{f}=f^{\ell}$ with period $n$ and such that
$\log\big\|\big(D\widetilde{f}^n(p)\big)^{-1}\big\|^{-1}>0$ (for
this take any prime $\ell\in\NN$ big enough).

\begin{Lemma}\label{LemmaExpandingPeriodicPreOrbit}If $p$
is a periodic repeller of period $n$ and
$\co_f^-(p)\cap\cc=\emptyset$ then there exists $\ell>0$ such that
$p$ is a periodic point of period $n$ with respect to
$\widetilde{f}=f^{\ell}$ and the pre-orbit of $p$ with respect to
$f$ is a $(24\,\log2)$-expanding set for $\widetilde{f}$.
\end{Lemma}
\dem
Let $m\ge1$ and set $\cc_{f^m}=\bigcup_{j=0}^{m-1}f^{-j}(\cc)$, the critical set of $f^m$.
As $\co_f^-(p)\cap\cc=\emptyset$, we get $\co_f^-(p)\cap\cc_{f^m}=\emptyset$.
From this,
it follows that, for every $0<\delta<\dist(\co_f^+(p),\cc_{f^m})$ and all $y\in\co_f^-(p)$,
$$\lim_{j\to+\infty}
\frac{1}{j} \sum_{i=0}^{j-1}-\log \mbox{dist}_{\delta}((f^m)^i(x),\cc_{f^m})
=\lim_{j\to+\infty}
\frac{1}{j} \sum_{i=0}^{j-1}-\log \mbox{dist}_{\delta}((f^m)^i(p),\cc_{f^m})=0.
$$
Thus, $\co_f^{-}(p)$ satisfies the slow approximation condition
with respect to $f^m$ (and the critical set of $f^m$), for every
$m\ge1$.

Let $n_0$ be such that $\lambda_0$ $:=$
$\log\big\|\big(Df^{n\,n_0}(p)\big)^{-1}\big\|^{-1}$ $>$ $0$. Let
$n_1$ be such that $e^{\lambda_0\,n_1/8}$ $>$ $8$. Let $\ell$ be a
prime number such that $\ell\ge n_0\,n_1$. Thus
$\log\big\|\big(Df^{n\,\ell}(p)\big)^{-1}\big\|^{-1}$ $>$
$24\,\log2$ and $p$ is a periodic point with period $n$ for the
map $\widetilde{f}=f^\ell$.

As $p$ is a periodic point, $\co_f^{-}(p)$ is a positively invariant
set (indeed $\co_f^{-}(p)$ is an invariant set) with respect to
$f$ and, as a consequence, it is a positively invariant set with
respect to $\widetilde{f}$.

Finally, as $$\lim_{j\to\infty}\frac{1}{j}\log\|(Df^{n\,j}(y))^{-1}\|=\lim_{j\to\infty}\frac{1}{j}\log\|(Df^{n\,j}(p))^{-1}\|$$
for all $y\in\co_f^-(p)$, it follows that $\co_f^{-}(p)$ is a $(24\,\log2)$-expanding set with respect to $\widetilde{f}$.

\cqd

Recall that the $\alpha$-limit set of a point $x$, denoted by
$\alpha_f(x)$, is the set
of accumulating points of
$\co_f^{-}(x)=\bigcup_{j=0}^{\infty}f^{-j}(x)$, the pre-orbit of
$x$.
As there is a substantial number of example of dynamics exhibiting
periodic repellers whose $\alpha$-limits have non-empty interior,
the proposition below show an easy way to find expanding measures
with big support in the topological sense (non-empty interior).

\begin{Proposition}\label{PropositionBacana}
Let $f:M\to M$ be a non-flat map and $\cc$ its critical/singular set. If
there is a periodic repeller $p$ contained in the interior of its
$\alpha$-limit set and such that $\co_f^-(p)\cap\cc=\emptyset$
then  there is an open neighborhood $\Delta$ of $p$ and an
uncountable collection $\mathcal{M}$ of ergodic invariant
probabilities such that if $\mu\in\mathcal{M}$ then all of its
Lyapunov exponents are positive and the support of $\mu$ contains
$\Delta$. Furthermore, $\exists\,\ell\ge1$ such that every
$\mu\in\mathcal{M}$ is an invariant ergodic zooming probability
for $\widetilde{f}=f^{\ell}$ and
$$\lim\frac{1}{n}\sum_{j=0}^{n-1}\log\|D\widetilde{f}(\widetilde{f}^j(x))^{-1}\|^{-1}\ge8$$
for $\mu$ almost every $x\in M$.
\end{Proposition}
\dem Let $\ell$ be as in
Lemma~\ref{LemmaExpandingPeriodicPreOrbit} and
$\widetilde{f}=f^{\ell}$. It follows from
Proposition~\ref{PropositionABV} that there are $\delta$ and
$\theta_0>0$ such that
$$\limsup_{n\to\infty}\frac{1}{n}\#\{1\le j\le n\,;\,x\in\th_j(e^{-6\log 2},\delta,\widetilde{f})\}\ge\theta_0$$
for every $x\in\co_f^-(p)$.

From Proposition~\ref{PropositionHyperbolicBalls} it follows that
each $(e^{-6\log 2},\delta)$-hyperbolic time is a
$(\alpha,\delta)$-zooming time and every $(e^{-6\log
2},\delta)$-hyperbolic pre-ball is a $(\alpha,\delta)$-zooming
pre-ball (all with respect to $\widetilde{f}$), where
$\alpha=\{\alpha_n\}$ and $\alpha_n(r)=(1/8)^{n}(r)$. In particular,
$\co_f^-(p)$ is a $(\alpha,\delta)$-zooming set.

Noting that $\sum_n\alpha_n(r)<r/4$ $\forall\,r>0$, let $0<r<\delta/4$ be small such that
$\overline{B_r(p)}\subset\interior(\alpha_f(p))$.
Thus $\Delta:=B_r^{\star}(p)$ is a $(\alpha,\delta)$-zooming nested ball (with respect to $\widetilde{f}$) containing $B_{r/2}(p)$.

Let $\widetilde{\fz}=(\fz(x))_{x\in\limsup\th_n(e^{-6\log
2},\delta,\widetilde{f})}$ be the collection of zooming images of
$\widetilde{f}$ that are $(e^{-6\log 2},\delta)$-expanding images
and let $\widetilde{\ce}$ be the collection of all
$\widetilde{\fz}$-pre-balls. Let $R$ be the ``first
$\widetilde{\ce}$-return time to $\Delta$'', $F$ be the induced
map associated to the  ``first $\widetilde{\ce}$-return time to
$\Delta$'' and $\cp$ be the Markov partition associated to the
``first $\widetilde{\ce}$-return time to $\Delta$'' as in
Definition~\ref{DefinitionRetIndTime},~\ref{DefinitionIndMapToRetIndTime}~and~\ref{DefinitionMarkCollecFirstRetTime}
(all this definitions applied to $\widetilde{f}$ instead to $f$).
By Corollary~\ref{CorollaryDisjointness}, $(F,\cp)$ is an induced
full Markov map for $\widetilde{f}$ on $\Delta$ with inducing time
$R$.

As the zooming points are dense on $\Delta$ (because $\co_f^-(p)$
is dense), $\{R>0\}=\bigcup_{P\in\cp}P$ is an open and dense
subset of $\Delta$. Let
$\mathfrak{B}=\bigcap_{j\ge0}F^{-j}(\{R>0\})$, that is,
$\mathfrak{B}$ is the set of points $x\in\Delta$ such that
$F^j(x)\in\Delta$ $\forall\,j\ge0$. Of course, $\mathfrak{B}$ is a
residual set of $\Delta$. Furthermore, $\mathfrak{B}$ is a metric
space with the distance induced by the distance of $M$ and its
topology is the induced topology.

Let $\cw $ be the collection of subsets of $\mathfrak{B}$ formed
by the empty set $\emptyset$ and all $Y\subset \mathfrak{B}$ such
that
$Y=(F|_{P_1})^{-1}\circ\dots\circ(F|_{P_s})^{-1}(\mathfrak{B})$
for some sequence of $P_1,...,P_s\in\cp$. Note that $\cw $
generates all open sets of $\mathfrak{B}$.

Let $\mathcal{A}$ be the collection of all sequence of numbers
$\{a_P\}_{P\in\cp}$ satisfying $a_P\in(0,1)$,
$\sum_{P\in\cp}a_P=1$ and $\sum_{P\in\cp}a_P\,R(P)<\infty$.

Choose any $\{a_P\}_{P\in\cp}\in\mathcal{A}$. Given any
$Y\in\cw\setminus\{\emptyset,\mathfrak{B}\}$, there is a unique
sequence $P_1,...,P_s\in\cp$ such that
$Y=(F|_{P_1})^{-1}\circ\dots\circ(F|_{P_s})^{-1}(\mathfrak{B})$.
In this case, define $\nu(Y)=\prod_{j=1}^s a_{P_j}$. Set also
$\nu(\emptyset)=0$ and $\nu(\mathfrak{B})=1$. It easy to see that
$\nu(A\cup B)=\nu(A)+\nu(B)$ for every $A,B\in\cw$ with $A\cap
B=\emptyset$. Moreover, $\nu(A)\le\nu(B)$ for every $A,B\in\cw$
with $A\subset B$. As $\cw$ generates the Borel algebra of
$\mathfrak{B}$, $\nu$ can be extended as a measure on the Borel
set of $\mathfrak{B}$. Furthermore, $\nu$ is $F$-invariant.
Indeed, given $Y\in\cw$, say
$Y=(F|_{P_1})^{-1}\circ\dots\circ(F|_{P_s})^{-1}(\mathfrak{B})$,
we get
$$\nu\big(F^{-1}(Y)\big)=\nu\bigg(F^{-1}\bigg((F|_{P_1})^{-1}\circ\dots\circ(F|_{P_s})^{-1}(\mathfrak{B})\bigg)\bigg)=$$
$$=\sum_{P\in\cp}\nu\bigg((F|_P)^{-1}\bigg((F|_{P_1})^{-1}\circ\dots\circ(F|_{P_s})^{-1}(\mathfrak{B})\bigg)\bigg)=$$
$$=\sum_{P\in\cp}a_P\,a_{P_1}\,\cdots\,a_{P_s}=a_{P_1}\,\cdots\,a_{P_s}\,\underbrace{\sum_{P\in\cp}a_P}_{1}=\nu(Y).$$

Note that if $x\in\mathfrak{B}$ and $P_0,P_1,P_2,\cdots$ is the
itinerary of $x$ by $F$ (i.e., $P_j$ is the element of $\cp$ that
contains $F^j(x)$ $\forall\,j\ge0$) then
$x=\bigcap_{j=0}^{\infty}C_n(x)$, where
$C_n(x)=(F|_{P_0})^{-1}\circ\dots\circ(F|_{P_j})^{-1}(\mathfrak{B})$
is the $n$-th cylinder containing  $x$. As
$$\frac{\mu(F(C_n(x))}{\mu(C_n(x))}=\frac{a_{P_1}\,a_{P_2}\,\cdots\,a_{P_n}}
{a_{P_0}\,a_{P_1}\,a_{P_2}\,\cdots\,a_{P_n}}=\frac{1}{a_{P_0}}\,\,\forall\,n>0,$$
one can prove that the Jacobian of $F$ with respect to $\mu$,
$J_{\mu}F(x)$ is well defined and it is constant on every
$P\in\cp$ (indeed, $J_{\mu}F|_P=a_P$). This implies that
$\frac{J_{\mu}F^n(x)}{J_{\mu}F^n(y)}=1$ for all $y\in C_n(x)$. As
a consequence, one can easily conclude that $\mu$ is ergodic with
respect to $F$.

As $\int R\,d\nu=\sum_{P\in\cp}R(P)\nu(P)=\sum_{P\in\cp}R(P)\,a_P<\infty$, it follows from Remark~\ref{RemarkProjection} that
$$
\widetilde{\mu}=\sum_{P\in\cp}\sum_{j=0}^{R(P)-1}\widetilde{f}_{\ast}^j\left(\nu|_{P}\right)\,
$$
is a $\widetilde{f}$-invariant finite measure. Furthermore, as
$\nu$ is $F$-ergodic, it follows that $\widetilde{\mu}$ is
$\widetilde{f}$-ergodic.

Applying Corollary~\ref{CorollaryRogers}, we conclude that
\begin{equation}\label{EquationAQZ1}
\limsup_n\frac{1}{n}\#\{0\le j<n\,;\,\widetilde{f}^{j}(x)\in\co_F^+(x)\}>0
\end{equation}
for $\widetilde{\mu}$ almost every $x\in\Delta$. If for some $n$
and $i$ we have $\widetilde{f}^n(x)=F^i(x)$ then, by construction,
$x$ belongs to some $(e^{-6\log 2},\delta)$-hyperbolic pre-ball
$V_n(y,\widetilde{f})$ with respect to $\widetilde{f}$. As a
consequence, $n$ is a zooming time for $x$ (with respect to
$\widetilde{f}$) and
$\frac{1}{n}\sum_{j=0}^{n-1}\log\|D\widetilde{f}(\widetilde{f}^j(x))^{-1}\|^{-1}$
$\ge$ $\bigg(\big(e^{-6\log 2})^{1/2}\bigg)^{-1}=8$, whenever
$\widetilde{f}^{n}\in\co_F^+(x)$. Thus
\begin{equation}
\label{EquationTYR4E}
\limsup_n\frac{1}{n}\#\{0\le j<n\,;\,x\in\tz_j(\alpha,\delta,\widetilde{f})\}>0
\end{equation}
and
\begin{equation}
\label{EquationJKJKL}
\limsup\frac{1}{n}\sum_{j=0}^{n-1}\log\|D\widetilde{f}(\widetilde{f}^j(x))^{-1}\|^{-1}\ge8
\end{equation}
for $\widetilde{\mu}$ almost every $x\in\Delta$. As
$\widetilde{\mu}(\Delta)>0$, it follows by ergodicity that
(\ref{EquationTYR4E}) and (\ref{EquationJKJKL}) holds for
$\widetilde{\mu}$ almost every point $x\in M$.

One can easy check that $\mu=\sum_{j=0}^{\ell-1}f_{\ast}^j\widetilde{\mu}$
is an ergodic $f$-invariant measure such that
(\ref{EquationTYR4E}) and (\ref{EquationJKJKL}) holds for $\mu$
almost every point $x\in M$. Therefore, $\mu$ is a zooming measure
(with respect to $\widetilde{f}$) and it follows from (\ref{EquationJKJKL})
that all Lyapunov exponent are positive for $\mu$ almost every point of $M$.

Of course distinct elements of $\mathcal{A}$ give rise to
distinct ergodic $F$-invariant probabilities. By ergodicity these
probabilities are mutually singular and so the $f$-invariant
measures generated from them are also mutually singular. So, as
$\mathcal{A}$ are uncountable, this process gives an uncountable
collection of $f$-invariant measures.

To finalize the proof, note that $\mu(U)\ge\nu(U)>0$ for every
open subset of $\Delta$, because every open subset of $\Delta$
contains a non-empty $Y\in\cw$, and $\nu(Y)>0$ for all
$\emptyset\ne Y\in\cw$. This implies that $\supp\mu\supset\Delta$.

\cqd

\begin{Definition}
A map $f:M\to M$ is called {\em strongly topologically transitive}
if we get $\bigcup_{n\ge0}f^n(U)=M$ for all open set $U\subset M$.
\end{Definition}

\begin{Theorem}\label{TheoremBacana}
Let $f:M\to M$ be a $C^{1+}$ map, possible with a critical region
$\cc$. If $f$ is strongly topologically transitive and it
has a periodic repeller  then some iterate of $f$ admits an
uncountable number of ergodic invariant expanding probabilities
whose supports are the whole manifold.
\end{Theorem}
\dem As $f$ is strongly topologically transitive,
$\alpha_f(x)=M$ for every point $x\in M$. In particular, if $p$ is
a periodic repeller, we get $\alpha_f(p)=M\ni p$. Thus, we can
apply Proposition~\ref{PropositionBacana}. Let $\Delta$, $\ell$
and $\cm$ be given by Proposition~\ref{PropositionBacana}.

Let $\mu\in\cm$. As $f$ is strongly topologically transitive, given any open
set $U\subset M$ there is $n\ge0$ such that
$f^{-n}(U)\cap\Delta\ne\emptyset$. Thus
$\mu(U)=\mu(f^{-n}(U))\ge\mu(f^{-n}(U)\cap\Delta)>0$ for every
$\mu\in\cm$ (because $\Delta\subset\supp\mu$). This implies that
$\supp\mu=M$ $\forall\,\mu\in\cm$.

Given $\mu\in\cm$, we have that
$$\lim\frac{1}{n}\sum_{j=0}^{n-1}\log\|D\widetilde{f}(\widetilde{f}^j(x))^{-1}\|^{-1}\ge8$$
for $\mu$ almost every $x\in M$. As  the equation above implies
that there are no negative Lyapunov exponents with respect to any
iterate of $f$ and for $\mu$-almost all point, it follows from
Lemma~\ref{LemmaNonFlatToSlowRecurrence} that $\mu$ is an
expanding measure with respect to $\widetilde{f}$. \cqd

\begin{Corollary}\label{CorollaryBacana}
If a $C^{1+}$ map $f:M\to M$, possibly with a critical region
$\cc$, is strongly topologically transitive and it has a
periodic repeller then the set of periodic repeller is dense on
$M$.
\end{Corollary}
\dem
This corollary follows from Theorem~\ref{TheoremBacana} and the fact that the support of any
expanding invariant measure is contained in the closure of the periodic repellers
(see Lemma~\ref{LemmaRepeller} of the appendix).
\cqd

\begin{Example}[See figure~\ref{FIGUREExampleAppl1}]
\label{ExampleAppl1}
Let $f:[0,1]\to[0,1]$ be given by
$$f(x)=\left\{\begin{array}{cc}
    g(x) & \mbox{if }x<1/2 \\
    1-g(1-x) & \mbox{if }x\ge1/2
  \end{array}\right.$$ where $g(x)=x+2x^2$.

The map $f$ can be seen as a $C^{\infty}$ map of the circle
$S^1=\RR/\ZZ$ and this map is topologically conjugated to the
uniformly expanding map $h(x)=2x$ (mod $\ZZ$). Thus, $f$ is
strongly topologically transitive. Note that $f$ has expanding periodic points.
Indeed, $f$ has a periodic point $p\in(0,1)$ of period two
(because $h$ does) and, as $Df^2(x)>1$ $\forall\,x\in(0,1)$, it
follows that $p$ is an expanding periodic point. Thus, follows
from Theorem~\ref{TheoremBacana} that some iterate of $f$ admits
an uncountable number of ergodic invariant expanding probabilities
whose supports are the whole circle.
%%%%%%%%%%%%%%%%%%%%%%%%%%%%%%%%%%%%%%%%%%%%%%
\begin{figure}\label{FIGUREExampleAppl1}
  \includegraphics{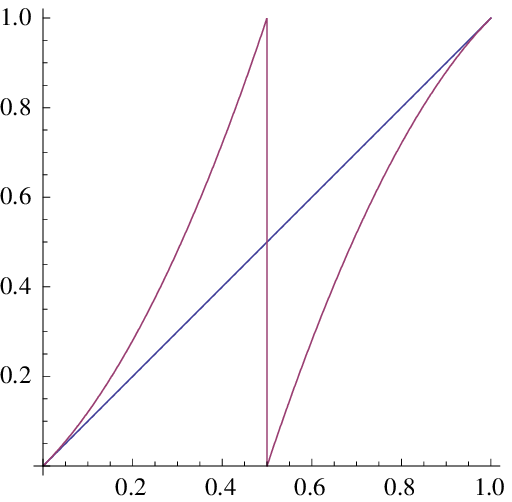}\\
  \caption{\label{FIGUREExampleAppl1}}
\end{figure}
%%%%%%%%%%%%%%%%%%%%%%%%%%%%%%%%%%%%%%%%%%%%%%
\end{Example}

In \cite{Y2} Young shows that for maps like $f$ of
Example~\ref{ExampleAppl1}
$\frac{1}{n}\sum_{i=0}^{n-1}\delta_{f^i(x)}$ converges weakly to
the Dirac measure at $0$ for Lebesgue almost every $x$. In particular, $f$
admits no invariant measures that is absolutely continuous with
respect to the Lebesgue measure. Furthermore, for Lebesgue almost
every point the Lyapunov exponent is zero. In contrast, it follows
from Theorem~\ref{TheoremBacana} that $f$ admits an uncountable
number of ergodic invariant probabilities whose supports are the
whole manifold and whose Lyapunov exponents are positive.

\begin{Example}\label{ExampleAppl2}
Let $F:[0,1]^2\to[0,1]^2$ be the skew product given by
$$F(x,y)=\big(f(x),(1+x)\phi(y)\big)$$
where $f$ is as in Example~\ref{ExampleAppl1} and
$\phi(y)=1/2-|y-1/2|$ is the ``tent'' map of slope one. Taking any
periodic point $p\in(0,1)$ for $f$, it is easy to see that
$\psi(y)=G(p,y)$ is an uniformly expanding map, where
$(f^{n}(x),G(x,y))=F^n(x,y)$ and $n$ is the period of $p$. Thus
taking any periodic point $q\in(0,1)$ with respect to $\psi$, it
follows that $(p,q)\in(0,1)^2$ is an expanding periodic point of
$F$. It is not difficult to check that $F\big|_{F([0,1]^2)}$ is
strongly topologically transitive and so it follows from
Theorem~\ref{TheoremBacana} that  some iterate of $F$ admits an
uncountable number of ergodic invariant expanding probabilities
whose supports are $F([0,1]^2)$.

\end{Example}

As in Example~\ref{ExampleAppl1}, the scenario of the expanding
invariant measures of Example~\ref{ExampleAppl2} is much richer
than the Lebesgue measure scenario. Indeed, as
$$DF^n(x,y)=\left(
                \begin{array}{cc}
                  (f^n)'(x) & 0 \\
                  * & ({+/-})\prod_{j=0}^{n-1}(1+f^j(x)) \\
                \end{array}
              \right),$$
it follows that the Lyapunov exponents of a point $(x,y)$ are
$\limsup\frac{1}{n}\log|(f^n)'(x)|$, the Lyapunov exponent of $x$
with respect to $f$, and
$\limsup\frac{1}{n}\log\prod_{j=0}^{n-1}(1+f^j(x))$. As
$\limsup\frac{1}{n}\log$ $|(f^n)'(x)|=0$ and $0$ $\le$
$\limsup\frac{1}{n}\log$ $\prod_{j=0}^{n-1}(1+f^j(x))$ $\le$
$\limsup\frac{1}{n}\sum_{j=0}^{n-1}f^j(x)$ $=$ $0$ for
Lebesgue-a.e. $x\in[0,1]$ (Theorem~5 of \cite{Y2}), we conclude
that the {\em all Lyapunov exponents for Lebesgue almost every
point are zero}.

One can find other examples of expanding measures in, for instance, \cite{BDV04} and \cite{Do08a}.

Let us apply Theorem~\ref{TheoremBacana} to unimodal maps. For this, we note that
every non-flat $S$-unimodal map $f:[0,1]\to[0,1]$
without a periodic attractor has an expanding periodic point
$p\in(0,1)$. Moreover, one can show that a non-flat $S$-unimodal map $f$ has an
expanding periodic $p\in(0,1)$ with dense pre-orbit if and only if
$f$ is not an infinitely renormalizable map and $f$ does not have
a periodic attractor. Thus, we get from Theorem~\ref{TheoremBacana} the following corollary.

\begin{Corollary}
If $f:[0,1]\to[0,1]$ is a non-flat $S$-unimodal map then one and
only one of the following alternatives can occur.
\begin{enumerate}
\item $f$ has a periodic attractor. \item $f$ is an infinitely
renormalizable map. \item $f$ admits an expanding invariant
probability whose support has non-empty interior (indeed an
uncountable number of these probabilities).
\end{enumerate}
\end{Corollary}

The corollary above shows that the dynamic of any $S$-unimodal map in the complement of the Axiom A and the infinitely
renormalizable maps exhibits uncountable many non trivial expanding measures, even when there are not SRB measures.

\begin{Theorem}[Markov structure for expanding sets of local diffeomorphisms] \label{TheoremaLOCALDIFFEO}
If $f:M\to M$ a $C^{1+}$ map is a local diffeomorphism then the
set of points with all Lyapunov exponents positive admits a Markov
structure. Furthermore, given any finite collection of ergodic
invariant probabilities $\{\mu_1,\cdots,\mu_s\}$ with all Lyapunov
exponents positive there is an induced map $(F,\cp)$ as in
Theorem~\ref{GlobalExpandingInducedMarkovMap} such that every
$\mu_j$ is liftable to this map.
\end{Theorem}
\dem
Let $\lambda>0$ and, for each $\ell\in\NN$, let $\Lambda_{\ell}$ be the set of $\lambda$-expanding points of $M$ with respect to $f^{\ell}$ (as $f$ is a local diffeomorphism there is no slow approximation condition). It follows from Lemma~\ref{LemmaExpPosImpliesNUE} that $\bigcup_{\ell\in\NN}\Lambda_{\ell}$ contains the set of points with all Lyapunov exponents positive (indeed, it is equal). As each $\Lambda_{\ell}$ has a Markov structure with respect to $f^{\ell}$ (and so, a Markov structure with respect to $f$), it follows that $\bigcup_{\ell\in\NN}\Lambda_{\ell}$ has a Markov structure with respect to $f$.

Given any finite collection of ergodic invariant probabilities $\{\mu_1,\cdots,\mu_s\}$ with all Lyapunov
exponents positive, there is some $\ell\in\NN$ such that $\mu_j(\Lambda_{\ell})=1$ for all $j=1,\cdots,s$.
Thus, applying Theorem~\ref{GlobalExpandingInducedMarkovMap}, there is an induced map $(F,\cp)$
such that every $\mu_j$ is liftable to this map, $\forall\,j=1,\cdots,s$.
\cqd

\subsection{Maps with a dense expanding set}

Besides the previous examples there are many examples of maps with
a dense expanding set. Indeed, most of the results of the so
called ``non-uniformly expanding maps'' was done with the
hypothesis of an expanding set of full Lebesgue measure, in
particular, a dense expanding set. This is, for instance, the case
of Viana maps (see Example~\ref{Viana_maps} below).

The crucial property used in the proof of
Proposition~\ref{PropositionBacana} and
Theorem~\ref{TheoremBacana} is indeed the existence of an expanding (or zooming)
set that is dense and also some condition to spread open sets to the whole
manifold. In the theorem below the hypotheses are chosen to obtain
these properties again.

\begin{Theorem}\label{TheoremBACANAdenso}
Let $f:M\to M$ be a transitive non-flat map with $\#f^{-1}(x)<\infty$ $\forall\,x\in M$.
If $f$ has a dense
$\lambda$-expanding set, $\lambda>0$, then there is an uncountable
collection of ergodic invariant $\lambda'$-expanding
probabilities, $\lambda'\ge\lambda/8$, whose support are the whole manifold.
\end{Theorem}
\dem
Given any $x\in\th_{n}(e^{-\lambda/4},\varepsilon,f)$, let $V_n(x)$ is the $(e^{-\lambda/4},\delta)$-hyperbolic pre-ball of center $x$ and order $n$, where $\varepsilon,\delta>0$ follows from Proposition~\ref{PropositionABV} and \ref{PropositionHyperbolicBalls}.
Note that $$W:=\bigcap_{j=0}^{\infty}\bigcup_{n\ge j}\bigg(\bigcup_{x\in\th_{n}}V_n(x)\bigg)$$
is a residual set, where $\th_{n}=\th_{n}(e^{-\lambda/4},\varepsilon,f)$.
Thus, the set of points $x\in W$ that are transitive ($\omega(x)=M$) is also a residual set (because the set of transitive points is residual). Choose a transitive point $q\in W$. As $q\in W$, there are sequences $n_k\to\infty$ and $x_k\in \th_{n_k}$ such that $q\in V_{n_k}(x_k)$ $\forall\,k\in\NN$ and $\lim_{k\to\infty}f^{n_k}(x_k)=p$, for some $p\in M$.
Of course that $x_{k}\to q$, indeed, $\dist(x_{k},q)$ $\le$ $e^{-(\lambda/8)n_k}$ $\delta$, $\forall\,k$.

Let $\alpha=\{\alpha_n\}_n$, where $\alpha_n(r)=e^{-(\lambda/8)n}r$. As $f$ is backward separated (because $\#f^{-1}(x)<\infty$ $\forall\,x\in M$) and as $\sup_{r>0}\sum_{n\ge1}\alpha_n(r)/r<+\infty$, we can choose any $0<r<r_0$ and consider the $(\alpha,\delta)$-zooming nested ball $B_{r}^{\star}(p)$, where $0<r_0<\delta/2$ is given by Lemma~\ref{LemmaZoomingNestedBall}.

We claim that there is $\Lambda\subset B_{r}^{\star}(p)$ dense in $B_{r}^{\star}(p)$ and such that every $x\in\Lambda$ has as hyperbolic return to $B_{r}^{\star}(p)$, that is, given $x\in\Lambda$ there is $s\ge1$ such that $x\in\th_{s}$ and $f^{s}(x)\in B_{r}^{\star}(p)$. Indeed, for each $y\in B_{r}^{\star}(p)$ and $\gamma>0$ one can find $\widetilde{y}\in\co^+(q)$, say $\widetilde{y}=f^{i}(q)$, so that $\dist(\widetilde{y},y)<\gamma/2$. Taking $k>i$ big enough so that $\dist(f^{i}(x_{k}),\widetilde{y})$ $=$ $\dist(f^{i}(x_{k}),f^{i}(q))$ $<$ $\gamma/2$, it follows that $\dist(f^{i}(x_{k}),y)<\gamma$, $f^{i}(x_{k})\in\th_{n_k-i}$ and $f^{n_k-i}(x_k)\in B_{r}^{\star}(p)$.

Now, the proof follows as the proof of Proposition~\ref{PropositionBacana} with a single difference. Here we do not need to consider an iterate $\widetilde{f}=f^{\ell}$ of $f$. Taking $\widetilde{f}=f$ and $\Delta=B_{r}^{\star}(p)$, construct the induced map $F$ and everything else as in the proof of Proposition~\ref{PropositionBacana}.
\cqd

\subsection{Decay of correlation and the Central Limit Theorem}

In \cite{ALP,ALP1} Alves, Luzzatto and Pinheiro study the decay of correlation associated to the decay of the tail of expanding moments. There it was proved that a polynomial decay of the tail of expanding moments, measured by the Lebesgue measure, implies a polynomial decay of correlation for the absolutely continuous invariant measure with respect to the Lebesgue measure (the SRB measure). It was also proved that the Central Limit Theorem holds for the SRB whenever the tail of expanding moments decays more then quadratically. In \cite{G} Gouëzel complemented this study for Lebesgue measure by showing that an exponential (or a stretched exponential) decay of the tail of expanding moments, measured by the Lebesgue measure, implies an exponential (or a stretched exponential) decay of correlation for the SRB measure. Here our construction permits to extend the results of \cite{ALP,G} for general expanding measures.

\subsubsection{Tail estimate for Maps with bounded derivative}
Follows from Lemma~\ref{LemmaSlowTofm} of the appendix that if $f:M\to M$ is a local
$C^{1+}$ diffeomorphism in the whole manifold except  in  a {\em
non-degenerate critical/singular set} $\mathcal{C}\subset {M}$ and
$\sup\{\|Df(x)\|$ $;$ $\,x\in M\setminus\cc\}$ $<$ $+\infty$ then
every set $U$ satisfying the slow approximation condition with
respect to $f$ also satisfies the slow approximation condition
with respect to $f^m$, $\forall m\ge 1$. On the other hand, a
direct calculation shows that for every big $n$
$$\frac{1}{n}\sum_{i=0}^{n-1} \log
\|\big(Df(f^{i}(x))\big)^{-1}\|^{-1}>\lambda$$
\begin{equation}\label{Equationkjkjj1}\Longrightarrow\end{equation}
$$\frac{1}{[n/m]}\sum_{i=0}^{[n/m]-1} \log
\|\big(Df^m((f^m)^{i}(x))\big)^{-1}\|^{-1}>\lambda,$$
where $[n/m]=\max\{j\in\ZZ; j\le n\}$ is the integer part of $n/m$.
Thus, it is easy to obtain the following result.

\begin{Lemma}\label{LemmaExpftofm}
If $f:M\to M$ is a local $C^{1+}$ diffeomorphism in the whole
manifold except  in  a {\em non-degenerate critical/singular set}
$\mathcal{C}\subset {M}$ and $\sup\{\|Df(x)\|$ $;$ $\,x\in
M\setminus\cc\}$ $<$ $+\infty$ then every $\lambda$-expanding set
with respect to $f$ is $(m \lambda)$-expanding with respect to
$f^m$ for every $m\ge 1$.
\end{Lemma}

As a consequence of the lemma above, every $\lambda$-expanding measure
with respect to $f$ is a $(m \lambda)$-expanding measure with
respect to $f^m$ $\forall\,m\ge 1$, whenever $\|Df\|$ is bounded.

\begin{Theorem}[Global expanding induced Markov map for
maps with bounded derivative]
\label{GlobalExpandingInducedMarkovMapForMapsWithBV}
Let $f:M\to M$
be a local $C^{1+}$ diffeomorphism in the whole manifold except in
a {\em non-degenerate critical/singular set} $\mathcal{C}\subset
{M}$ and suppose that $\sup\{\|Df(x)\|$ $;$ $\,x\in
M\setminus\cc\}$ $<$ $+\infty$. Given $\lambda>0$, let $\ell=\min\{k\in\NN$ $;$ $k\ge(16\log 3)/\lambda\}$.
If $\ch$ is a $\lambda$-expanding
set with $\lambda>0$ then there exist an induced
Markov map $(F,\cp)$  with inducing time $R$, and a finite
partition $\cp_0$ of $M$ by essentially open sets satisfying
all items of Theorem~\ref{GlobalExpandingInducedMarkovMap} and also the following ones.
\begin{itemize}
\item $\ch\subset\{R>0\}$.
\item $\cp|_{\ch}$ $=$ $\{P\cap\ch$ $;$ $P\in\cp\}$ is an induced Markov partition of $\ch$ with respect to $f^{\ell}$.
\item\label{ItemE6b} There are $r>0$, $\varepsilon>0$ and $n_0\in\NN$ such that
$$\{R>n\}\cap\ch\subset\{h>n/\ell\}\,\,\,\forall\,n\ge n_0,$$
where
$$h(x)=\inf\{j>0\,;\,\frac{1}{\,j}\sum_{k=0}^{j-1}\log\|Df(f^{k}(x))^{-1}\|^{-1}\ge\lambda\mbox{\;
and}$$
$$\frac{1}{\,j}\sum_{k=0}^{\,j-1}-\log\dist_{\varepsilon}(f^k(x),\cc)
\le r\}.$$
\end{itemize}
\end{Theorem}
\dem
As $Df$ is bounded, $f^{\ell}$ is also a non-flat map with a critical region $\cc_{\ell}=\bigcup_{j=0}^{\ell-1}$.
In this case, follows from Lemma~\ref{LemmaExpftofm}, Proposition~\ref{PropositionABV} and \ref{PropositionHyperbolicBalls} that (1) $\ch\subset\limsup_n\th_{n}(1/9,\varepsilon_{\ell},f^{\ell})$ for some $\varepsilon_{\ell}>0$; (2) the whole $\ch$ is a $(\alpha,\delta_{\ell})$-zooming set (with respect to $f^{\ell}$) for some $\delta_{\ell}>0$ and with $\alpha=\{\alpha_n\}_n$ and $\alpha_n(r)=(1/9)^n r$.
Furthermore, one can take $0<\varepsilon<\varepsilon_{\ell}$ so that
\begin{equation}
\th_{\ell\,j}(e^{-\lambda/4},\varepsilon,f)\subset\th_{j}(1/9,\varepsilon_{\ell},f^{\ell}),\,\,\forall\,j\ge1.
\end{equation}
This produces a proper zooming sub-collection ${\widetilde{\fz}}=(\widetilde{\fz}(x))_{x\in\ch}$  associated to the $(e^{-\lambda/4},\varepsilon)$-hyperbolic time (both with respect to $f$). The important point here is that ${\widetilde{\fz}}$ is defined on the whole $\ch$ (in the proofs of Theorem~\ref{GlobalMarkovStructureForZooming} and \ref{GlobalExpandingInducedMarkovMap} we didn't have any warranty that the zooming sub-collection covered all the given zooming set when $\ell>1$).
Thus, applying Theorem~\ref{GlobalMarkovStructureForZooming}, we can get all items of Theorem~\ref{GlobalExpandingInducedMarkovMap}. Furthermore, in this case $\ch$ $\subset$ $\{R>0\}$ $=$ $\bigcup_{P\in\cp}P$ and, by Remark~\ref{RemarkIndMarPar78}, we get that $\cp|_{\ch}$ is an induced Markov partition of $\ch$ with respect to $f^{\ell}$.

Let $r=\frac{b}{4}\lambda$ and $b=\frac{1}{3}\min\{1,1/\beta\}$, with $\beta>0$ being the constant that appears in the definition of non-flatness (see Section~\ref{ApplicationsExpandingMeasures}). It follows from (\ref{Equationkjkjj1}), (\ref{Equationkjkjj2}) and from Propositions~\ref{PropositionABV} that $\th_{\ell\,j}(e^{-\lambda/4},\varepsilon,f)$ can be estimate by $h(\ell\,j+s)$ for any $0\le s<\ell$. That is, $\th_{\ell\, j}(e^{-\lambda/4},\varepsilon,f)\subset\{h=n\}$ if $j=[n/\ell]$ and $n$ is big enough. Thus, the last item of the theorem follows from item (\ref{ItemE6}) of Theorem~\ref{GlobalExpandingInducedMarkovMap}.

\cqd

\begin{Theorem}[Decay of correlation]\label{TheoremDacayOfCorrelation}
Let $f:M\to M$ a $C^{1+}$ map with a {\em non-degenerate critical set} $\mathcal{C}\subset
{M}$. Suppose that $\mu$ is a non-flat $\lambda$-expanding probability, $\lambda>0$, and that $f^n\big|_{\supp\mu}$ is transitive for all $n\ge1$. Then there is a unique absolutely continuous invariant probability $\nu\ll\mu$. The probability $\nu$ is mixing and it basin contains $\mu$ almost every point of $M$. Furthermore, there exist $r,\epsilon>0$ such that, if we set $$h(x)=\inf\{j>0\,;\,\frac{1}{\,j}\sum_{k=0}^{j-1}\log\|Df(f^{k}(x))^{-1}\|^{-1}\ge\lambda\mbox{\;
and}$$
$$\frac{1}{\,j}\sum_{k=0}^{\,j-1}-\log\dist_{\varepsilon}(f^k(x),\cc)
\le r\},$$ then we have the following estimates for the decay of correlation
$$\mathrm{Cor}(\phi,\psi\circ f^n)=\bigg|\int \phi\,\psi\circ f^n d\nu-\int\phi d\nu\,\int\psi d\nu\bigg|$$ for any given functions $\phi,\psi:M\to\RR$ with $\phi$ H\"{o}lder and $\psi$ bounded.
\begin{enumerate}
\item If $\mu\{h>n\}=O(n^{-\gamma})$ for some $\gamma>0$, then $\mathrm{Cor}(\phi,\psi\circ f^n)=O(n^{-\gamma/\ell})$,
    where $\ell=\min\{k\in\NN\,;k\ge(16\log3)/\lambda\}$.
\item If $\mu\{h>n\}=O(\exp(-\rho\,n^\gamma))$ for some $\rho,\gamma>0$, then there exist $\widetilde{\rho},\widetilde{\gamma}>0$ such that $\mathrm{Cor}(\phi,\psi\circ f^n)=O(\exp(-\widetilde{\rho}\,n^{\widetilde{\gamma}}))$.
\end{enumerate}
\end{Theorem}
\dem
Let $(F,\cp)$, $R$ and $\cp_0$ be the global induced Markov map, the inducing time and the finite partition of $M$ by essentially open sets given by Theorem~\ref{GlobalExpandingInducedMarkovMapForMapsWithBV}. Furthermore, we can assume that $(F,\cp)$, $R$ and $\cp_0$ are also given by Theorem~\ref{GlobalMarkovStructureForZooming}. Indeed,  for some $\delta>0$,   $\alpha=\{\alpha_n\}_n$ and $\alpha_n(r)=e^{-(\lambda/4)n} r$, $\mu$ is a $(\alpha,\delta)$-zooming measure and so, Theorem~\ref{GlobalExpandingInducedMarkovMap} and \ref{GlobalExpandingInducedMarkovMapForMapsWithBV} are straightforward application of Theorem~\ref{GlobalMarkovStructureForZooming} (considering for instance that $\Lambda\subset M$ is the set of $(\alpha,\delta)$-zooming points).

Of course there is a natural identification of an induced full Markov map with $\mu$-bounded distortion (see Definition~\ref{DefinitionMU-boundDist} and \ref{DefinitionInducdMMap}) with a Young Tower. To construct the Young Tower \cite{Y1,Y2} (or equivalently, an induced full Markov map with $\mu$-bounded distortion) we can proceed as in the proof of Theorem~4.1 of \cite{G}. In this theorem a global induced Markov map as in Theorem~\ref{GlobalExpandingInducedMarkovMapForMapsWithBV} induces a Young Tower with essentially the same estimates of the tail of the partition. We note that the Lebesgue measure is not important to the proof. Indeed the fundamental ingredient of Theorem~4.1 of \cite{G} is the Lemma~\ref{LemmaGouezel} below.

Let $\nu\ll\mu$ be the $f$-invariant probability given by Theorem~\ref{TheoremSRB}.
The construction of the local induced Markov map associated to the global one was already  done in the proof of Theorem~\ref{GlobalMarkovStructureForZooming}, this is precisely the induced map $\widetilde{F}_0$ given by (\ref{EquationinducedMap2}) and defined in the topological open ball $\Delta_{j_0}$.

To emphasize $f$ instead of its iterate $\widetilde{f}=f^{\ell}$, set $\tau(x)=\ell\,\widetilde{R}_0$ and  $F_0(x)=f^{\tau(x)}(x)=\widetilde{F}_0(x)$. Recall that as $\mu$ is a $\lambda$-expanding measure, $\ell=\min\{k\in\NN$ $;$ $(16\log 3)/\lambda\}$ (see Theorem~\ref{GlobalExpandingInducedMarkovMapForMapsWithBV}).

Let $k(x)=\ell\,\widetilde{k}(x)$, where $\widetilde{k}(x)$ is given by Claim~\ref{ClaimGlobalInd1}. Thus, $F_0(x)=F^{k(x)}(x)$. Set $t_j(x)=\sum_{i=0}^{j-1}R(F^{i}(x))$, for every $j\le k(x)$. Of course that $\tau(x)=t_{k(x)}(x)$.

Let $\cp^{n}$ be the partition of $M$ given by $\cp^n=\bigcap_{i=0}^{n-1}F^{-i}(\cp_0)$. As $\Delta_{j_0}$ is contained a $\nu$-ergodic component $U\subset M$ with respect to $f^{\ell}$ (see the proof of Theorem~\ref{GlobalMarkovStructureForZooming}), there is some $L>0$ such that every element of $\cp_0$ contains an element of $\cp^n$ whose image under $F^n$ is $\Delta_{j_0}$. From the distortion control it follows that there exist constants  $C_0,\varepsilon>0$ such that
$$\mu\big(\{\tau=t_j\mbox{ or ... or }\tau=t_{j+L-1}\,;\,t_1,\dots,t_{j-1},\tau>t_{j-1}\}\big)\ge\varepsilon$$
$$\mbox{and}$$
$$\mu\big(\{t_{j+1}-t_{j}>n\,;\,t_1,\dots,t_{j}\}\big)\le C_0\mu(\{R>n\}).$$
Thus, applying Lemma~\ref{LemmaGouezel} and Theorem~3 of \cite{Y2}  the theorem follows.
\subsubsection*{Remark} {\em  Although there is no explicit reference to the stretched exponential decay in the statement of Theorem~3 of \cite{Y2}, Young proofs can be adapted to the general case (see the comments in the proof Lemma 4.2. of \cite{G} and also the comments in the begging of Section~4 of the same paper).}

\cqd

\begin{Lemma}[Gouëzel \cite{G}]
\label{LemmaGouezel}
Let $(X,\mu)$ be a space endowed with a finite measure and $k:X
\to \NN$ and $t_0,t_1, t_2,\ldots : X\to \NN$ measurable functions such
that $0=t_0<t_1<t_2<\dots$ almost everywhere. Set $\tau(x)=t_{k(x)}(x)$,
and assume that there exist $L>0$ and $\epsilon>0$ such that
  \begin{equation}
  \label{young1prime}
  \mu\{\tau=t_j \text{ or }\ldots \text{ or }\tau=t_{j+L-1}\,;\,
  t_1,\ldots,t_{j-1}, \tau>t_{j-1}\} \geq
  \epsilon.
  \end{equation}
Assume moreover that there exist a positive sequence $u_n$ and a
constant $C_0$ such that
  \begin{equation}
  \label{young1seconde}
  \mu\{ t_{j+1}-t_j>n\,;\,t_1,\ldots,t_j\}
  \le C_0 u_n.
  \end{equation}
Then
  \begin{enumerate}
  \item  If $u_n$ has polynomial decay, $\mu\{ \tau > n\}=O(u_n)$.
  \item If $u_n=e^{-cn^\eta}$ with $c>0$ and $\eta\in (0,1]$, then
  there exists $c'>0$ such that $\mu\{ \tau >
  n\}=O(e^{-c'n^\eta})$.
  \end{enumerate}
\end{Lemma}

\begin{Theorem}[Central Limit Theorem]\label{TheoremCentralLimitTheorem}
Let $f$, $\mu$, $\nu$ and $h$ being as in Theorem~\ref{TheoremDacayOfCorrelation}.
If $\mu\{h>n\}=O(n^{-(1+\gamma)})$ for some $\gamma>(16\log3)/\lambda$ then the {\em Central Limit Theorem} holds for any H\"{o}lder function $\phi:M\to\RR$ such that $\phi\circ f\ne \psi\circ f-\psi$ for any $\psi$.
\end{Theorem}
\dem
As Theorem~\ref{TheoremDacayOfCorrelation}, this result follows from  Theorem~\ref{GlobalExpandingInducedMarkovMapForMapsWithBV}, Lemma~\ref{LemmaGouezel} and from \cite{Y2} (Theorem~4).
\cqd

\begin{Example}[Viana maps]\label{Viana_maps}
An important class of non-uniformly expanding dynamical
systems (with critical sets) in dimension greater than one was
introduced by Viana in \cite{V}.
This class of  maps can be described as
follows. Let $a_0\in(1,2)$ be such that the critical point $x=0$
is pre-periodic for the quadratic map $Q(x)=a_0-x^2$. Let
$S^1=\RR/\ZZ$ and $b:S^1\rightarrow \RR$ be a Morse function, for
instance, $b(s)=\sin(2\pi s)$. For fixed small $\alpha>0$,
consider the map
 \[ \begin{array}{rccc} \hat f: & S^1\times\RR
&\longrightarrow & S^1\times \RR\\
 & (s, x) &\longmapsto & \big(\hat g(s),\hat q(s,x)\big)
\end{array}
 \]
 where  $\hat q(s,x)=a(s)-x^2$ with
$a(s)=a_0+\alpha b(s)$, and $\hat g$ is the uniformly expanding
map of the circle defined by $\hat{g}(s)=ds$ (mod $\ZZ$) for some
integer $d\ge16$. It is
easy to check that for $\alpha>0$ small enough there is an
interval $I\subset (-2,2)$ for which $\hat f(S^1\times I)$ is
contained in the interior of $S^1\times I$. Thus, any map $f$
sufficiently close to $\hat f$ in the $C^0$ topology has
$S^1\times I$ as a forward invariant region. We consider from here
on these maps restricted to $S^1\times I$.

Most of the results for $f\in \cn$ are summarized below:
\begin{enumerate}
\item $\exists\,\ch\subset S^1\times I $, with full Lebesgue
measure on $S^1\times I$, and $\lambda>0$ such that $\ch$ is a
$\lambda$-expanding set with respect to $f$ \cite{V} (indeed, to
be coherent with the estimate (\ref{Equation_h}) we may assume
that $\ch$ is a $2\lambda$-expanding set); \item for each
$0<c<1/4$ and $\varepsilon>0$ there are constants
$C(c,\varepsilon)$ and $\delta(x)>0$ such that
$$\leb(\{x\,;\,h_{\varepsilon}(x)>n\})\le
C(c,\varepsilon)e^{-c\sqrt n}$$ for every $n\ge1$ \cite{AA,V},
where
$$
h_{\varepsilon}(x)=\inf\{j>0\,;\,\frac{1}{\,j}
\sum_{k=0}^{j-1}\log\|Df(f^{k}(x))^{-1}\|^{-1}\ge\lambda\mbox{\;
and}
$$
\begin{equation}\label{Equation_h}\frac{1}{\,j}\sum_{k=0}^{\,j-1}
-\log\dist_{\delta(\varepsilon)}(f^k(x),\cc) \le\varepsilon\};
\end{equation}
\item $f\big|_{f(S^1\times I)}$ is strongly topologically transitive and has a unique ergodic absolutely continuous
invariant (thus SRB) measure whose support is $f(S^1\times I)$
\cite{Al}; \item the Central Limit Theorem holds for $f$
\cite{ALP}; \item the correlations of Hölder functions decay at
least like $e^{-c'\sqrt{n}}$, for some $c'>0$ \cite{G}.
\end{enumerate}
\end{Example}

Of course Viana maps satisfies most of hypothesis of the theorems
in this section (Section~\ref{SectionExamplesAndAplications}). In
particular, it follows from Theorem~\ref{TheoremBACANAdenso} that
there is an uncountable collection $\mathcal{M}$ of ergodic
invariant probabilities such that all  Lyapunov exponents of every
$\mu\in\cm$ are positive and the support of any $\mu\in\cm$ is the
whole manifold. Furthermore, every $\mu\in\mathcal{M}$ is an $f$
invariant ergodic zooming probability  and
$$\lim\frac{1}{n}\sum_{j=0}^{n-1}\log\|D{f}({f}^j(x))^{-1}\|^{-1}\ge\lambda/2$$
for $\mu$ almost every $x\in M$ and every $\mu\in\cm$.

We can also apply
Theorem~\ref{TheoremSRB},~\ref{TeoremMarkovStructureForExpanding}
and ~\ref{GlobalExpandingInducedMarkovMap} to the Viana maps. From
Theorem~\ref{TheoremSRB}, we conclude that all non-flat expanding
measure admits an absolutely continuous invariant measure (in
particular one can apply this theorem to obtain the SRB
measure). Furthermore, we can apply Theorem~\ref{TheoremDacayOfCorrelation} (or Theorem~\ref{TheoremCentralLimitTheorem}) to study the decay of correlation (or the Central Limit Theorem) of any expanding invariant measure with bounded distortion (in particular, this gives a new proof of the decay of correlation of the SRB measure for the Viana maps). We can also study the decay for the measures that come from Theorem~\ref{TheoremBacana} or Theorem~\ref{TheoremBACANAdenso}.

\begin{Theorem}
If $f$ is a Viana map then there exist an uncountable number of ergodic invariant expanding measures with exponential decay of correlation and whose support is the whole $f(S^1\times I)$.
\end{Theorem}
\dem
The construction of the collection of expanding measures given by  Theorem~\ref{TheoremBACANAdenso} or Theorem~\ref{TheoremBacana} comes from that proof of Proposition~\ref{PropositionBacana}. This measures are associated to an induced full Markov map $(F,\cp)$ defined on a topological ball $\Delta$ and to the collection $\ca$ of all sequence $\{a_P\}_{P_\cp}$ satisfying $\sum_{P\in\cp}a_P=1$ and $\sum_{P\in\cp}a_P R(P)<\infty$, where $R$ is the induced time of $F$.
As one can see in the proof of Proposition~\ref{PropositionBacana}, each $a=\{a_P\}_{P_\cp}\in\ca$ generates a $F$-invariant measure $\nu_a$ and also a $f$-invariant measure $\mu_a$, with $\nu_a\ll\mu_a$. Moreover, we have a very good distortion control of $J_{\nu_a}F^n$ in every cylinder $C_n$. Indeed $\frac{J_{\nu_a}F^n(x)}{J_{\nu_a}F^n(y)}=1$ $\forall y\in C_n(x)$ (see details in the proof of Proposition~\ref{PropositionBacana}).
Let $a=\{a_P\}_{P_\cp}\in\ca$ be any sequence satisfying  $\lim_{n}\frac{1}{n}\log\big(\sum_{P\in\cp_n}a_P\big)=\gamma<0$, where $\cp_n=\{P\in\cp\,;\,R(P)=n\}$.
Thus, $\nu_a(\{R>n\})=\sum_{j>n}\nu_a(\sum_{P\in\cp_j}a_p)=O(e^{-\gamma\,n})$ and it follows from \cite{Y2} that $\mu_a$ has exponential decay of correlation.

\cqd

\subsection{Expanding measures on metric spaces}

In Section~\ref{SettingAndStatementOfMainResults} and \ref{ApplicationsExpandingMeasures}
we deal with expanding
sets and measures on Riemannian manifold because the standard way to define
these objects is using the derivative of the map. Precisely, using $\|(Df)^{-1}\|^{-1}$.
So, to extend the notion of expanding sets or measures we need to rewrite
this expression in terms of the distance.
For this, note that if $f:M\to M$ is differentiable at a point $p\in M$
then $$\|(Df(p))^{-1}\|^{-1}=\liminf_{x\to p}\frac{\dist(f(x),f(p))}{\dist(x,p)}.$$
Thus, given a metric spaces $X$ and $Y$ and a map $f:X\to Y$
define
$$\DD^{-}f(p)=\liminf_{x\to p}\frac{\dist(f(x),f(p))}{\dist(x,p)},$$
where we are using the notation $\dist$ to assign the distance on both spaces.
Define also
$$\DD^{+}f(p)=\limsup_{x\to p}\frac{\dist(f(x),f(p))}{\dist(x,p)}.$$

Of course one can rewrite the expanding
condition (\ref{EquationExpanding}) in terms of $\DD^- f$, that is,
\begin{equation}\label{EquationSecExpCond}
\limsup_{n\to\infty}\frac{1}{n}\sum_{j=0}^{\infty}\log(\DD^{-}f)\circ f^j(x)>0,
\end{equation}
and use this condition to define
the expanding condition on a metric space.
The critical/singular set $\cc$ can be defined as the set of points $x\in X$ having $\DD^{-}f(x)=0$ or $\DD^+f(x)=\infty$. In the condition of non-degenerateness we only need to replace the expressions (1) and (2)
by
\begin{enumerate}
\item[(1)]
\quad $\displaystyle{\frac{1}{B}dist(x,\mathcal{C})^{\beta}\le
\DD^{-}f(x)\le\DD^+f(x)\le B\,dist(x,\mathcal{C})^{-\beta} }$.
\end{enumerate}
and
\begin{enumerate}
\item[(2)] \quad $\displaystyle{\left|\log\DD^{-}f(x)-
\log\DD^{-}f(x)\:\right|\leq
\frac{B}{dist(x,\mathcal{C})^{\beta}}dist(x,y)}$.
\end{enumerate}

It is straightforward to check that, Proposition~\ref{PropositionABV},
Proposition~\ref{PropositionHyperbolicBalls} et cetera remain true if we
change $\|(Df)^{-1}\|^{-1}$ by $\DD^{-}f$. In particular, the
expanding sets and measures with this definition are zooming sets
and measures. As a consequence, if $X$ is a connected, compact,
separable metric space there are results analogues to
Theorem~\ref{TheoremSRB}, \ref{TeoremMarkovStructureForExpanding}
and \ref{GlobalExpandingInducedMarkovMap} for this context
(see Remark~\ref{RemarkCONNECTED} if $X$ is not connected).

\begin{Definition}The map
$f$ is called conformal at $p\in X$ if $\DD^{+}f(p)=\DD^{-}f(p)$.
In this case the conformal derivative of $f$ at $p$ is
$$\DD f(p)=\lim_{x\to p}\frac{\dist(f(x),f(p))}{\dist(x,p)}.$$
\end{Definition}
It is easy to check that the {\em chain rule} holds for the conformal derivative.
Moreover, it is obvious that one can rewrite the expanding
condition (\ref{EquationExpanding}) or (\ref{EquationSecExpCond}) in terms,
if it exists, of the conformal derivative $\DD f$.

An example of a conformal in this definition is the shift with the usual metric.

\begin{Example}[Expanding sets on a metric space]\label{ExpandingShift}Consider the {\em one-side shift}
$\sigma:\sum_2^+\to\sum_2^+$ with its usual metric, that is,
$$\dist(x,y)=\sum_{n=1}^{+\infty}\frac{|x_n-y_n|}{2^n},$$
where $x=\{x_n\}_n$ and $y=\{y_n\}_n$. It is easy to verify that $\sigma$
is a conformal map and that $\DD\sigma(x)=2$ $\forall\,x\in\sum_2^+$.
\end{Example}

As we could have expected, every positively invariant set (in particular the whole $\sum_2^+$) and all
invariant measure for the map $\sigma$ of the Example~\ref{ExpandingShift} are expanding.

In this paper we are basically interested in zooming and
expanding measures. As we saw, the set of zooming measures
contains the expanding measures.
Now we will give examples of zooming sets and measures that are not
expanding, i.e., examples of sets and invariant
measures that are zooming with only a polynomial backward contraction.

Note that if $f:X\to X$ is a conformal map defined on a compact metric space $X$
and $\DD f\le1$ then it follows by compactness that given any $\varepsilon>0$ there is
$\delta>0$ such that $$\dist(f(x),f(y))\le(1+\varepsilon)\dist(x,y)$$
$\forall\,x,y\in X$ satisfying $\dist(x,y)<\delta$. So, given any $\varepsilon>0$ there is $\delta>0$ such that
if $x,y\in X$, $n\ge1$ and
$\dist(f^j(x),f^j(y))<\delta$ $\forall 0\le j<n$ then
$$\dist(f^n(x),f^n(y))\le(1+\varepsilon)^n\dist(x,y),$$
that is,
$\DD f\le1$ prohibits any exponential backward contraction.
In particular, it does not admit any expanding set or measure.

\begin{Example}[Zooming but not expanding]\label{ShiftPolinomial}Consider the {\em one-side shift}
$\sigma:\sum_2^+\to\sum_2^+$ with its usual topology.
Consider the compatible metric given by
$$\dist(x,y)=\left\{\begin{array}{cc}
               0 & \mbox{if }x=y \\
               \big(\phi(x,y)\big)^{-2} & \mbox{if }x\ne y
             \end{array}.\right.
$$
where $x=\{x_n\}_n$, $y=\{y_n\}_n$ and $\phi(x,y)=\min\{n\ge1$ $;$ $x_n\ne y_n\}$. It is easy to verify that $\sigma$
is a conformal map and that $\DD\sigma(x)=1$ $\forall\,x\in\sum_2^+$.
In particular, $$\lim_{n\to\infty}\frac{1}{n}\log\DD {\sigma}^n(x)=\lim_n\frac{1}{n}\sum_{n=0}^{\infty}\log\DD\sigma({\sigma}^n(x))=0,\,\,\,\,\forall\,x\in{\sum}_2^+.$$
So, $\sigma$ does not admit any expanding set or measure. In contrast, given any $p\in\sum_2^+$ and $x,y$ $\in$ $C_n(p)$ $=$ $\{q\in\sum_2^+$ $;$ $p_1=q_1\dots p_n=q_n\}$, we have $\phi(\sigma^j x,\sigma^j y)=\phi(\sigma^n x,\sigma^n y)+(n-j)$, for $0\le j\le n$. Thus,
$$\sqrt{\dist(\sigma^j x,\sigma^j y)}=\frac{1}{\phi(\sigma^n x,\sigma^n y)+(n-j)}=\frac{\sqrt{\dist(\sigma^n x,\sigma^n y)}}{1+(n-j)\sqrt{\dist(\sigma^n x,\sigma^n y)}}$$
and so,
$$\dist(\sigma^j x,\sigma^j y)=\bigg(\frac{1}{1+(n-j)\sqrt{\dist(\sigma^n x,\sigma^n y)}}\bigg)^2\dist(\sigma^n x,\sigma^n y).$$
As a consequence, the cylinder $C_n(p)$ is a
$(\alpha,1)$-zooming pre-ball for $p$, where  $\alpha=\{\alpha_n\}_n$ and $\alpha_n(r):=(\frac{1}{1+n\sqrt{r}})^2 r$ (one can check that $\alpha_n\circ\alpha_j(r)=\alpha_{n+j}(r)$).
This implies that every positively invariant set of $\sum_2^+$ is an $(\alpha,1)$-zooming set and any $\sigma$-invariant
measure is $(\alpha,1)$-zooming.
\end{Example}

\subsubsection*{Future applications} Recently there was an increasing  development of the study of
the thermodynamic formalism beyond the uniformly hyperbolic
context (including countable Markov shift) by several authors
(this list is certainly not complete):
Araujo~\cite{Ara07},
Arbieto, Matheus, Oliveira, Varandas, Viana~\cite{AMO,Oli03,OV07,Var07,VaVi},
Bruin, Keller, Todd~\cite{BruTo06,BrK98,BT06,BruTo,BruTo07},
Buzzi, Paccaut, Sarig, Schmitt~\cite{Bu99,BPS01,BS03,Sar99,Sar01,Sar03},
Dobbs~\cite{Do08b}, Denker, Keller, Nitecki, Przytycki, Rivera-Letelier, Urba{\'n}ski~\cite{DKU90, DU91a,DNU95,DPU96,DU91f,DU92b,PrRLe,Ur98},
Leplaideur, Rios~\cite{LR06},
Pesin, Senti, Zhang~\cite{PS05,PS06,PZ07,PSZ07},
Wang, Young~\cite{WY},
Yuri~\cite{Yur97,Yur99,Yur03}.

In many cases, a natural place to look for an equilibrium state is the set of expanding measures.
Thus, as we have shown the existence of an expanding invariant measure that is absolutely continuous with respect to any given reference expanding measures with a ``nice'' Jacobian (the Jacobian of conformal measures associated to a H\"{o}lder potential has a ``nice'' Jacobian), one can expect to get among  the invariant expanding measures a (non-lacunary) Gibbs measure that is a candidate to the equilibrium state. Furthermore, as we have constructed an induced Markov map such that every ergodic expanding invariant measure is liftable to it (in particular, allowing a symbolic study of these measures), we believe that the results presented in this paper can be useful to the program of extending the thermodynamic formalism to the general non-uniformly hyperbolic setting.

%%%%%%%%%%%%%%%%%%%%%%%%%%%%%%%%%%%%%%%%%%%%%%%%%%%%%%%%%%%%%%%%%%%%%%%%%%%%%%

\section{Appendix}

Let $\mu$ be a finite
measure defined on the Borel sets of the compact,
separable metric space ${X}$.
 The proof of the following fact can be found in, for instance,
 \cite{KaHa} (Lemma~A.6.8.).

\begin{Lemma}\label{LemmaSmallPartition}
Let ${U}\subset {X}$ be a measurable set. Given $\varepsilon>0$
there exists a finite partition {\em (mod $\mu$)} $\cp$ of ${U}$
by open sets of $U$ (in the induced topology) such that $\diam
(\cp)<\varepsilon$.
\end{Lemma}

Now, let $M$ be a  compact Riemannian manifold  of dimension
$d\ge1$.

\begin{Lemma}\label{LemmaNonFlatToSlowRecurrence}
Let $f:M\to M$ be a $C^{1+}$ map. If $\mu$ is a $f$-invariant
ergodic probability with all of its Lyapunov exponent finite then
$\mu$ satisfies the slow approximation condition, that is, for each $\varepsilon>0$ there is a $\delta>0$
such that
$$
\limsup_{n\to+\infty}
\frac{1}{n} \sum_{j=0}^{n-1}-\log \mbox{dist}_{\delta}(f^j(x),\cc)
\le\varepsilon,
$$
for $\mu$ almost every $x\in M$.
\end{Lemma}
\dem

Let $\cc$ be the critical region of $f$ (of course we may assume
that $\cc\ne\emptyset$). As $f$ is $C^{1+}$,
$\cc$ is a compact set and also $\det Df$ is Holder. That is, $\exists\,k_0,k_1>0$ such that  $|\det Df(x) - \det Df(y)|$ $\le$ $k_0\dist(x,y)^{k_1}$ $\forall\,x,y\in M$. Given $x\in M$ there is $y_x\in\cc$ such that $\dist(x,y_x)=\dist(x,\cc)$. Thus, we get $|\det Df(x)|$ $=$ $|\det Df(x) - \det Df(y_x)|$ $\le$ $k_0\dist(x,y_x)^{k_1}$ $=$ $k_0\dist(x,\cc)^{k_1}$.
That is, $\log|\det Df(x)|$ $\le$ $\log k_0+k_1\log\dist(x,\cc)$.
If $\int\log|\det Df|\,d\mu=-\infty$, it follows from Birkhoff that $\sum_{i}\lambda_i=\lim\frac{1}{n}\sum_{j=0}^{n-1}\log|\det Df(f^j(x))|=-\infty$
for $\mu$-almost every $x$, contradicting our hypothesis as $\sum_{i}\lambda_i$ is the sum of the Lyapunov exponents. So,  $-\infty$ $<$ $\int\log|\det Df|\,d\mu-\log k_0$ $\le$
 $k_1\int\log\dist(x,\cc)\,d\mu(x)$ $\le$ $k_1\log\diameter(M)$. As the logarithm of the distance to the critical set is integrable, it follows that $$\int\log\dist_{e^{-n}}(x,\cc)\,d\mu(x)=
 \int_{\{x\,;\,\log\dist(x,\cc)<-n\}}\log\dist(x,\cc)\,d\mu\rightarrow0$$
 when $n\rightarrow\infty$. This implies (by Birkhoff)  the slow approximation condition.
\cqd

From Lemma~\ref{LemmaNonFlatToSlowRecurrence} follows the Corollary~\ref{CorollaryNonFlatToSlowRecurrence}.

\begin{Corollary}\label{CorollaryNonFlatToSlowRecurrence}
Let $f:M\to M$ be a $C^{1+}$ map. An
ergodic invariant probability $\mu$ is an expanding measure if and only if
(\ref{EquationExpanding0}) holds for $\mu$ almost every $x\in M$
\end{Corollary}

The lemma below is a remark that appears in Section~1.1 of
\cite{ABV}.

\begin{Lemma}\label{LemmaExpPosImpliesNUE}
Let $f:M\to M$ be a $C^1$ local diffeomorphism and let $\mu$ be a $f$-invariant
probability. If for
$\mu$-almost every $x\in m$ we have
\begin{equation}\label{EquationY7ygy}
\lim_{n\to\infty}\log|Df^n(x)\,v|>0,\,\forall\,|v|=1,
\end{equation}
then there exist an iterate $\tilde{f}=f^{\ell}$ of $f$
such that
\begin{equation}\label{EquationY7ygyhh}\lim_{n\to\infty}\frac{1}{n}\sum_{j=0}^{n-1}\log
\|D\widetilde{f}(\widetilde{f}\,^j(x))^{-1}\|^{-1}>0
\end{equation}
for
$\mu$-almost every $x\in M$.
\end{Lemma}
\dem By the compactness of $M$, (\ref{EquationY7ygy}) implies that there is $\lambda>0$ such that for each $x\in M$  $\exists\,n_x\in\NN$ satisfying $\log|Df^n(x)\,v|\ge2\lambda$ $\forall\,|v|=1$ and $\forall\,n\ge n_x$, that is, $$\log\|(Df^{n}(x))^{-1}\|^{-1}\ge2\lambda\,\,\forall\,n\ge n_x.$$
Let $K=|\min_{x\in M}\log\|(Df^n(x))^{-1}\|^{-1}|$ and let $\varepsilon>0$ be such that $\varepsilon(1+K/\lambda)<1$. Let $\ell\ge1$ be so that $\mu(U)>1-\varepsilon$, where $U=\{x\in M$ $;$ $\log\|(Df^{\ell}(x))^{-1}\|^{-1}>\lambda\}$. Thus,
$$\int\log\|(Df^{\ell})^{-1}\|^{-1}d\mu>\lambda\mu(U)-K(1-\mu(U))=\lambda(1-\varepsilon(1+K/\lambda))>0$$
and the proof of the lemma follows from Birkhoff.
\cqd

\begin{Lemma}[Transporting the slow approximation from $f$ to $f^m$]\label{LemmaSlowTofm}
Let $f:M\to M$ be a local $C^{1+}$ diffeomorphism in the whole
manifold except  in  a {\em non-degenerate critical/singular set}
$\mathcal{C}\subset {M}$. Suppose that $K=\sup\{\|Df(x)\|;\,x\in
M\setminus\cc\}<+\infty$. Given $m>1$, let $F=f^m$ and let $\cc_F$
be the critical region of $F$, i.e.,
$\cc_F=\bigcup_{j=0}^{m-1}f^{-j}(\cc)$. If $\delta<{K^{-m}}$ then,
for every $n\in\NN$,
\begin{equation}
\label{Equationkjkjj2}\sum_{j=0}^{n-1}-\log\dist_{\delta/K^m}(F^j(x),\cc_F)\le 2
\sum_{j=0}^{m  n-1}-\log\dist_{\delta}(f^j(x),\cc).
\end{equation}
\end{Lemma}
\dem Given $\delta<{K^{-m}}$, $x\in M$ and $n\in\NN$,
set$$J_n=\{0\le j<n\,{;}\,\dist(F^j(x),\cc_F)<\delta/K^m\}.$$ Of
course we may assume that  $\sum_{j=0}^{m
n}-\log\dist_{\delta}(f^j(x),\cc)<+\infty$.

Given $j\in J_n$, let $0\le s_j<m$ be such that
$\dist(F^j(x),\cc_F)=\dist(F^j(x),f^{-s_j}(\cc))$ and let
$\gamma:[0,1]\to M$ a continuous path from $F^j(x)$ to
$f^{-s_j}(\cc)$ with length equal to
$\dist(F^j(x),f^{-s_j}(\cc))$. Note that $\gamma([0,1))\subset
M\setminus\bigcup_{j=0}^{m-1}f^{-j}(\cc)$ and so,
$f^i(\gamma([0,1)))\cap \cc=\emptyset$ $\forall\,0\le i<m$. Thus,
$\dist(f^{j m+s_j}(x),\cc)$ $=$
$\dist(f^s(F^j(x)),f^{s_j}(f^{-s_j}(\cc)))$ $\le$
$\mbox{length}(f^{s_j}\circ\gamma)$ $\le$
$K^{s_j}\mbox{length}(\gamma)$ $\le$ $K^m\dist(F^j(x),\cc_F)$.
That is,
\begin{equation}\label{EquationRecorrencia1}
\dist(f^{j m+s_j}(x),\cc)\le K^m\dist(F^j(x),\cc_F).
\end{equation}
As $j\in J_n$, we get $K^m\dist(F^j(x),\cc_F)<\delta<K^{-m}$ and
then
\begin{equation}\label{EquationRecorrencia2}
\dist(f^{j m+s_j}(x),\cc)<\delta<K^{-m}.
\end{equation}
Multiplying (\ref{EquationRecorrencia2}) by $\dist(f^{j
m+s_j}(x),\cc)$,
\begin{equation}\label{EquationRecorrencia3}
(\dist(f^{j m+s_j}(x),\cc))^2<K^{-m}\dist(f^{j m+s_j}(x),\cc).
\end{equation}
Combining  (\ref{EquationRecorrencia1}) and
(\ref{EquationRecorrencia3}),
\begin{equation}\label{EquationRecorrencia4}
(\dist(f^{j m+s_j}(x),\cc))^2<\dist(F^j(x),\cc_F),\;\forall j\in
J_n.
\end{equation}
From (\ref{EquationRecorrencia4}), we get  for every $j\in J_n$
$$-\log\dist(F^j(x),\cc_F)<-2\log\dist(f^{j m+s_j}(x),\cc)=$$
$$\underbrace{=}_{\mbox{using (\ref{EquationRecorrencia2})}}-2\log\dist_{\delta}(f^{j m+s_j}(x),\cc)
\le 2\sum_{i=0}^{m-1}-\log\dist_{\delta}(f^{j m+i}(x),\cc).$$

Therefore,
$$\sum_{j=0}^{n-1}-\log\dist_{\delta/K^m}(F^j(x),\cc_F)
=\sum_{j\in J_n}-\log\dist(F^j(x),\cc_F)\le$$
$$\le 2 \sum_{j\in J_n}\sum_{i=0}^{m-1}-\log\dist_{\delta}(f^{j m+i}(x),\cc)\le$$
$$\le 2 \sum_{j=0}^{n-1}\sum_{i=0}^{m-1}-\log\dist_{\delta}(f^{j m+i}(x),\cc)=
2\sum_{j=0}^{m  n-1}-\log\dist_{\delta}(f^j(x),\cc)$$ \cqd

\begin{Lemma}\label{LemmaRepeller}
The support of any ergodic invariant expanding measure $\mu$,
with respect to a non-flat map $f:M\to M$ (possibly with a critical/singular region $\cc$),
is contained in the closure of the set of periodic repellers of $f$.
Furthermore, for each $\varepsilon>0$ there is a periodic repeller whose orbit is $\varepsilon$-dense on the support of $\mu$.
\end{Lemma}
\dem
Let $\mu$ be an ergodic invariant expanding measure and $p\in M$ be a $\mu$-generic point.
Thus, $\omega(p)=\supp\mu$.
By Proposition~\ref{PropositionHyperbolicBalls}, there is a sequence $n_j\to\infty$
of hyperbolic times for $p$ and a sequence of
hyperbolic pre-balls $V_{n_j}(p)$ with $f^{n_j}$ mapping
$\overline{V_{n_j}(p)}$ diffeomorphically onto the ball
$\overline{B_{\delta}(f^{n_j}(p))}$.

Let $m\ge1$ be big enough so
that $\{p,f(p),\cdots,f^m(p)\}$ is $\delta/10$ dense on $\supp\mu$.
Given any $\varepsilon>0$, let $k_0$
be big enough so that $\{f^m(p),$ $f^{m+1}(p),$ $\cdots,$ $f^{k_0}(p)\}$ is $\varepsilon/2$-dense on $\supp\mu$.
Let $0<r_0$ be small so that $f^m\big|_{{B_{r_0}(p)}}$
is a diffeomorphism and $\diameter\big(f^j(B_{r_0}(p))\big)<\delta/10$
$\forall\,0\le j\le m$ (as $x$ is an expanding point, note that
$\{p,f(p),\cdots,f^m(p)\}\cap\cc=\emptyset$).
Choose $0<r<\varepsilon/3$ so that $\overline{B_r(f^m(p))}\subset f^m(B_{r_0}(p))$ and let
$U=\big(f^m\big|_{{B_{r_0}(p)}}\big)^{-1}\big(B_r(f^m(p))\big)$.
Note that every ball of radius $\delta/2$ and center on a point
of $\supp\mu$ contain at least one of
the pre-images $U,f(U),\cdots,f^m(U)=B_r(f^m(p))$
(because $\{p,f(p),\cdots,f^m(p)\}$ be $\delta/10$ is dense).

Let $k>>k_0$ be a very big hyperbolic time for $f^m(p)$.
Thus, the diameter of the associated pre-ball $V_k(f^m(p))$ is smaller
then $r/2$ and so, $\overline{V_k(f^m(p))}\subset B_r(f^m(p))$.
As noted before, $B_{\delta/2}(f^{k+m}(p))$ contains the closure of some $f^s(U)$. So $f^k(V_k(f^m(p)))=B_{\delta}(f^{k+m}(p))\supset\overline{f^s(U)}$.
Let $W=\big(f^k\big|_{V_k(f^m(p))}\big)^{-1}(f^s(U))\subset V_k(f^m(p))\subset B_r(f^m(p))$.
Thus, $f^{k+m-s}$ maps $\overline{W}\subset U$ diffeomorphically onto $\overline{U}$. Furthermore,
as we can choose $k$ as big as we want, the expansion
of $f^k|_W$ is as big as we want. On the other hand,
we can loose expansion only on the transport of $f^s(U)$ to $f^m(U)=B_r(f^m(p))$
and this is at most $m$ steps. Therefore, it follows that $g=(f^{k+m-s}|_W)^{-1}$
is a contraction. In particular, $f^{k+m-s}$ has a repeller fixed point $\widetilde{q}\in W$.
Of course, $\widetilde{q}$ is a periodic repeller for $f$ and, as the diameter of $W$ is as small as $k$ is big, $\widetilde{q}$ is as close of $f^m(p)$ as we want.
From this follows that $\{q,$ $f(q),$ $\cdots,$ $f^{k_0}(q)\}$ is $\varepsilon$-dense.
\cqd

%%%%%%%%%%%%%%%%%%%%%%%%%%%%%%%%%%%%%%%%%%%%%%%%%%%%%%%%%%%%%%%%%%%%%%%%%%%%%%

\end{document}